\documentclass[11pt,reqno]{amsart}

\usepackage[T1]{fontenc}
\usepackage[utf8]{inputenc}
\usepackage[english]{babel}

\usepackage{lmodern}

\usepackage[a4paper,margin=1.15in]{geometry}
\usepackage{microtype}
\usepackage{mathtools}
\usepackage{amssymb}
\usepackage{bm}
\usepackage{braket}
\usepackage{extpfeil}
\SetSymbolFont{stmry}{bold}{U}{stmry}{m}{n}
\allowdisplaybreaks
\numberwithin{equation}{section}

\usepackage{array}
\usepackage{tabularx}

\usepackage[shortlabels]{enumitem}
\usepackage[title]{appendix}
\usepackage[all,poly]{xy}
\usepackage{tikz}
\usetikzlibrary{arrows.meta,positioning}
\usepackage{capt-of}

\usepackage{xcolor}
\definecolor{linkblue}{RGB}{0,72,128}
\usepackage[
colorlinks=true,
linkcolor=linkblue,
citecolor=linkblue,
urlcolor=linkblue
]{hyperref}
\usepackage[nameinlink,noabbrev]{cleveref}

\theoremstyle{plain}
\newtheorem{thm}{Theorem}[section]
\newtheorem{theorem}[thm]{Theorem}
\newtheorem{lemma}[thm]{Lemma}

\newtheorem{corollary*}[thm]{Corollary*}
\newtheorem{proposition}[thm]{Proposition}
\newtheorem{proposition*}[thm]{Proposition*}

\theoremstyle{definition}

\newtheorem{definition}[thm]{Definition}

\newtheorem{defn-thm}[thm]{Definition-Theorem}

\theoremstyle{remark}
\newtheorem{remark}[thm]{Remark}

\newcommand{\Rmnum}[1]{\uppercase\expandafter{\romannumeral #1}}

\tikzset{
	solid black dot/.style={
		circle,
		fill=black,
		inner sep=1.5pt
	}
}

\def\d{\mathrm{d}}

\newcommand{\Ga}{\Gamma}

\newcommand{\bC}{\mathbb{C}}

\newcommand{\bP}{\mathbb{P}}
\newcommand{\bQ}{\mathbb{Q}}
\newcommand{\bR}{\mathbb{R}}

\newcommand{\bT}{\mathbb{T}}

\newcommand{\cH}{\mathcal{H}}

\newcommand{\cM}{\mathcal{M}}
\newcommand{\cO}{\mathcal{O}}

\newcommand{\fp}{\mathfrak{p}}
\newcommand{\fl}{\mathfrak{l}}

\newcommand{\fQ}{\mathfrak{Q}}

\newcommand{{\inv} }{\mathrm{inv}}
\newcommand{\ev}{\mathrm{ev}}

\newcommand{\vir}{{\mathrm{vir}}}
\newcommand{\Jac}{ {\mathrm{Jac}} }

\begin{document}
\raggedbottom
	
	 \title[$\mathbb P^1$ and resolved conifold]{All-Genus Equivalence Between $\mathbb{P}^1$ and Resolved Conifold with Diagonal and Anti-Diagonal Actions via Topological Recursion}
	\author{Dingkang Xiong}
	\address{Dingkang Xiong, Department of Mathematical Sciences, Tsinghua University, Haidian District, Beijing 100084, China}
	\email{xdk22@mails.tsinghua.edu.cn}
	\author{Jinghao Yu}
	\address{Jinghao Yu, School of Mathematical Sciences, Peking University, Haidian District, Beijing 100871, China}
	\email{yujinghao@math.pku.edu.cn}
	\author{Jian Zhou}
	\address{Jian Zhou, Department of Mathematical Sciences, Tsinghua University, Haidian District, Beijing 100084, China}
	\email{jianzhou@mail.tsinghua.edu.cn}
	
	\begin{abstract}
		We study the Liu--Yang--Zhang--Zhou (LYZZ) conjectures
		\cite{liuEquivariantGromovWitten2022} relating the descendant Gromov--Witten partition functions of \(\mathbb P^1\) and the equivariant resolved conifold with two torus actions through linear changes of dilaton-shifted coordinates on their big phase spaces.
		We prove the transformation for the anti-diagonal torus action and establish a Hodge-twisted version for the diagonal torus action via Eynard--Orantin topological recursion. In particular, we provide the spectral curve for the resolved conifold with the anti-diagonal torus action and the diagonal torus action with Hodge twists.
	\end{abstract}
	\maketitle
	
	\setcounter{tocdepth}{1}
	\tableofcontents

\section{Introduction}
\label{sec:introduction}

Understanding relations between the Gromov--Witten theories of
different targets is an important problem in enumerative geometry.
At genus zero, quantum cohomology organizes the primary invariants into a Frobenius manifold
\cite{dubrovinGeometry2DTopological1996,
kontsevichGromovWittenClassesQuantum1994}, from which the all-genus theory can be reconstructed through several complementary frameworks.
The Dubrovin--Zhang approach encodes higher-genus descendant potentials into tau-functions of integrable hierarchies obtained by deforming the genus-zero principal hierarchy
\cite{dubrovinNormalFormsHierarchies2001}.
In the semisimple setting, the Givental--Teleman reconstruction recovers higher-genus descendant and ancestor potentials from genus-zero data via quantized $S$- and $R$-matrix actions
\cite{giventalSemisimpleFrobeniusStructures2001,
giventalGromovWittenInvariantsQuantization2001,
telemanStructure2DSemisimple2012}.
Alternatively, Eynard--Orantin topological recursion governs higher-genus invariants through a spectral curve $(\Sigma, x, y, B)$, constructing symmetric multidifferentials $\omega_{g,n}$ via residue calculus at the simple critical points of $x$
\cite{chekhovFreeEnergyTopological2006,chekhovHermiteanMatrixModel2006,eynardInvariantsAlgebraicCurves2007}.
When the spectral curve is compatible with a semisimple cohomological field theory (CohFT), the Dunin-Barkowski--Orantin--Shadrin--Spitz (DOSS) correspondence identifies $\omega_{g,n}$ with the all-genus stable range invariants of the CohFT \cite{dunin-barkowskiIdentificationGiventalFormula2014}; see also \cite{milanovEynardOrantinRecursion2014,xiongEYNARDORANTINTOPOLOGICALRECURSION,chekhovTopologicalRecursionHard2019}.

These frameworks suggest a natural question: when two enumerative problems are linked at genus zero via transformations of their Frobenius structures or principal hierarchies, does this relation lift to their all-genus descendant potentials? This paper addresses this question for the Gromov-Witten invariants of the projective line $\mathbb P^1$ and the equivariant Gromov-Witten invariants of the resolved conifold
\[
X=\operatorname{Tot}\!\left(
\mathcal O_{\mathbb P^1}(-1)\oplus
\mathcal O_{\mathbb P^1}(-1)
\right),
\]
equipped with the anti-diagonal and diagonal $\mathbb C^*$-actions having trivial base action and fiber weights $(1,-1)$ and $(1,1)$, respectively. To distinguish the notations,
we  denote  by $X^{\mathrm{ad}}$ and $X^{\mathrm{di}}$ the resolved conifold $X$ with these two torus actions.

The interplay between the Gromov--Witten theories of $\mathbb P^1$ and $X$ grew out of closely related developments in integrable systems and Frobenius manifolds.
For $\mathbb P^1$, initiated by the matrix-model formulation and Toda conjecture of Eguchi and Yang \cite{eguchiTOPOLOGICALCP1MODEL1994a}, the integrable hierarchy description was developed in \cite{getzlerTODACONJECTURE2001,zhangCP1TopologicalSigma2002} and formalized via the extended Toda hierarchy by Carlet, Dubrovin, and Zhang \cite{carletExtendedTodaHierarchy2004}, with its Virasoro symmetries and connection to descendant potentials established in \cite{dubrovinVirasoroSymmetriesExtended2004}.
In the equivariant setting, Okounkov and Pandharipande described the descendant theory of $\mathbb P^1$ via the $2$-Toda hierarchy and Hurwitz theory with completed cycles \cite{okounkovEquivariantGromovWitten2006,okounkovGromovWittenTheory2006,okounkovGromovWittenTheoryHurwitz2009}.

For the anti-diagonal resolved conifold $X^{\mathrm{ad}}$, Brini linked the genus-zero theory to the dispersionless Ablowitz--Ladik hierarchy and conjectured an all-genus integrable description \cite{briniLocalGromovWitten2012}.
Brini, Carlet, and Rossi then realized the Ablowitz--Ladik hierarchy as a two-component reduction of the $2$-Toda hierarchy and showed that its Frobenius structure (with non-flat unity) is almost dual, in the sense of Dubrovin \cite{dubrovinAlmostDualityFrobenius2004}, to that of $X^{\mathrm{ad}}$ \cite{briniIntegrableHierarchiesMirror2012}.
For further developments concerning melting-crystal models, higher-rank $2$-Toda reductions, difference equations, and local bihamiltonian structures of the Ablowitz--Ladik hierarchy, see \cite{takasakiModifiedMeltingCrystal2013,briniRationalReductions2017,alimIntegrableHierarchyResolved2022,liTriHamiltonianAblowitzLadik2022}.

Brini \cite{briniLocalGromovWitten2012} also considered the genus zero potential of
the resolved conifold in the diagonal action case. 
The link to the $\bP^1$ was established through the work of several groups.
Riley and Strachan \cite{rileyNoteRelationshipRational2007} studied the almost dual of the $\bP^1$-Frobenius manifold
and the NLS-Frobenius manifold.
As is well-known from Dubrovin's work \cite{dubrovinGeometry2DTopological1996},
these latter two Frobenius manifolds are related to each other by a Legendre transformation
for Frobenius manifolds.
As noted in \cite{rileyNoteRelationshipRational2007},
their almost duals are related to each other by what they called a twisted Legendre transformation.  
Stedman and Strachan \cite{strachanGeneralizedLegendre2017} generalized this to generalized Legendre transformations
of generalized Frobenius manifolds. (See also \cite{liuGeneralizedFrobeniusManifolds2025,liuLegendreTransformationsClass2025} for further studies.) 
Liu, Yang, Zhang, and Zhou (LYZZ)  pointed out that the AL-Frobenius manifold and the $\bP^1$-Frobenius 
manifold are related to each other via generalized Legendre transformation. For more details, see also \cite{liuLegendreTransformationsClass2025}.  
Furthermore,
they observed almost duality between the $\bP^1$-Frobenius manifold and the generalized Frobenius manifold
defined by the diagonal potential function  in \cite{briniLocalGromovWitten2012}.
Note the potential function in \cite{rileyNoteRelationshipRational2007} has an equivalent but different form.

Motivated also by Brini's conjecture and the Hodge--GUE correspondence \cite{dubrovinHodgeGUECorrespondence2020}, LYZZ proposed exact identities relating the all-genus descendant partition functions of $\mathbb P^1$ to those of $X^{\mathrm{di}}$ and $X^{\mathrm{ad}}$ under explicit affine changes of variables, up to prescribed genus-zero quadratic terms. 
They proved both identities at genus zero via the principal hierarchies and provided evidence in genera one and two.

On the integrable-systems side, these genus-zero structures were recently placed into a systematic framework through the theory of generalized Frobenius manifolds with non-flat unity and their topological deformations \cite{liuGeneralizedFrobeniusManifolds2025}.
Within this program, the extended Ablowitz--Ladik hierarchy was realized as the topological deformation of the corresponding principal hierarchy \cite{liuExtendedAblowitzLadik2026} and shown to be related to the extended Toda hierarchy via generalized Legendre and Miura transformations \cite{dubrovinGeometry2DTopological1996,strachanGeneralizedLegendre2017,liuLegendreTransformationsClass2025}.
Nevertheless, directly extending the genus-zero proof of LYZZ along the Dubrovin--Zhang approach remains out of reach: the Frobenius manifolds of $X^{\mathrm{di}}$ and $X^{\mathrm{ad}}$ do not admit Euler vector fields, and a general Dubrovin--Zhang theory of higher-genus topological deformations for Frobenius manifolds without Euler vector fields (or an all-genus lifting of almost duality) has not yet been fully developed.

To bypass this obstacle, we turn to Eynard--Orantin topological recursion, which provides a direct route to compare the all-genus descendant potentials through their spectral curves without relying on Euler vector fields.

For $\mathbb P^1$, Norbury and Scott proposed that stationary Gromov--Witten invariants are governed by the spectral curve
\begin{equation}
x(\ell)=\ell+\ell^{-1},
\qquad y(\ell)=\log\ell,
\qquad
B(\ell_1,\ell_2)
=\frac{d\ell_1\,d\ell_2}{(\ell_1-\ell_2)^2},
\label{eq:intro-P1-spectral-curve}
\end{equation}
which was proved in all genera by DOSS \cite{norburyGromovWittenInvariants2014,dunin-barkowskiIdentificationGiventalFormula2014}.
Borot and Norbury further extracted both stationary and non-stationary descendant invariants by applying linear functionals to the multidifferentials $\omega_{g,n}$ of \eqref{eq:intro-P1-spectral-curve} \cite{borotLoopEquationsGromovWitten2019}, and Fang, Liu, and Zong extended this spectral-curve description to equivariant $\mathbb P^1$ \cite{fangEynardOrantinRecursionEquivariant2017}.

For toric Calabi--Yau threefolds, building on Mari\~no \cite{marinoOpenStringAmplitudes2008}, Bouchard, Klemm, Mari\~no, and Pasquetti (BKMP) formulated the remodeling conjecture, identifying open and closed Gromov--Witten invariants with the Eynard--Orantin topological recursion of the mirror curve \cite{bouchardRemodelingBModel2009,bouchardTopologicalOpenStrings2010}.
For $\mathbb C^3$, the open-string conjecture was proved for the one-legged topological vertex by Chen and Zhou independently \cite{chenBKMPConjecture2018,zhouLocalMirrorOneLegged2009} and for the three-legged vertex by Zhou \cite{zhouLocalMirrorTopologicalVertex2009}, alongside results for closed-string free energies \cite{bouchardTopologicalRecursionMirror2012,bouchardRemodelingFaberPandharipande2013,zhuBouchardSulkowski2015}.
The general conjecture was proved by Eynard and Orantin for smooth semi-projective toric Calabi--Yau threefolds \cite{eynardComputationOpenGromov2015}, and extended by Fang, Liu, and Zong to general semi-projective toric Calabi--Yau three-orbifolds \cite{liuAllgenusOpenclosedMirror2020,fangRemodelingConjectureToric2019}.
Fang, Liu, Yu, and Zong recently established a descendant version of the remodeling conjecture \cite{fangRemodelingConjectureDescendants2025}, providing an all-genus spectral-curve framework for equivariant descendant invariants of toric Calabi--Yau threefolds.

In this paper, we prove the all-genus conjectures of LYZZ \cite{liuEquivariantGromovWitten2022} for $X^{\mathrm{ad}}$, and establish a Hodge-twisted version for $X^{\mathrm{di}}$. Using Eynard--Orantin topological recursion and the DOSS correspondence, we bridge the spectral-curve descriptions of $\mathbb P^1$ and the resolved conifold by comparing their semisimple CohFT data and tracking the descendant linear functionals across the two curves, showing that the genus-zero LYZZ correspondences lift to all genera.

\medskip
\noindent\textbf{Main results.}
We prove the anti-diagonal correspondence and establish a
Hodge-twisted version of the diagonal correspondence. Let
$\mathbb E_g$ be the Hodge bundle and write
\[
\Lambda_g^\vee(1)=\sum_{j=0}^{g}(-1)^j c_j(\mathbb E_g).
\]
The twisted diagonal theory is defined by inserting
$\Lambda_g^\vee(1)^{-2}$ into the localized equivariant
Gromov--Witten integrals of $X^{\mathrm{di}}$, as specified in
\eqref{eq:corrected-diagonal-correlators}. Its genus-$g$ descendant
potential is denoted by $F_g^{X^{\mathrm{di}},\mathrm{tw}}$.
Throughout the paper, the equivariant parameter is specialized to $1$ after localization.

The following statement summarizes
Theorem~\ref{thm:full-correspondences}, proved in
Section~\ref{sec:proof}.

\begin{theorem}[LYZZ conjecture with Hodge twist]
\label{thm:intro-correspondences}
For every $g\geq0$, one has
\begin{equation}
\begin{aligned}
F_g^{\mathbb P^1}(s)
&=F_g^{X^{\mathrm{di}},\mathrm{tw}}
\bigl(t_{\mathrm{di}}(s)\bigr)
+\delta_{g0}A_{\mathrm{di}}(s),
\\
F_g^{\mathbb P^1}(s)
&=F_g^{X^{\mathrm{ad}}}
\bigl(t_{\mathrm{ad}}(s)\bigr)
+\delta_{g0}A_{\mathrm{ad}}(s).
\end{aligned}
\label{eq:intro-correspondences}
\end{equation}
For $k\geq0$, the diagonal descendant transformation is
\begin{equation}
\begin{aligned}
t_{\mathrm{di}}^{1,k}(s)
&=\sum_{m\geq0}\frac{m^k}{m!}\,s^{1,m}
-1+\delta_{k1},
\\
t_{\mathrm{di}}^{2,k}(s)
&=\sum_{m\geq0}
\frac{2k\,m^k s^{1,m}+(m+1)^k s^{2,m}}{m!}
-2k,
\end{aligned}
\label{eq:intro-di-descendant-transformation}
\end{equation}
and the anti-diagonal transformation is
\begin{equation}
\begin{aligned}
t_{\mathrm{ad}}^{1,k}(s)
&=\sum_{m\geq0}
\frac{\mathrm i^{\,m-1}m^k}{m!}\,s^{1,m}
-1+\delta_{k1},
\\
t_{\mathrm{ad}}^{2,k}(s)
&=\sum_{m\geq0}
\frac{
(2k+1)\mathrm i^{\,m-1}m^k s^{1,m}
+\mathrm i^{\,m}(m+1)^k s^{2,m}
}{m!}
-(2k+1).
\end{aligned}
\label{eq:intro-ad-descendant-transformation}
\end{equation}
The genus-zero correction polynomials are
\begin{equation}
\begin{aligned}
A_{\mathrm{di}}(s)
&=\sum_{k,\ell\geq0}
\frac{s^{1,k}s^{2,\ell}}{(k+\ell+1)k!\ell!}
-\sum_{\ell\geq0}\frac{s^{2,\ell}}{(\ell+2)\ell!},
\\
A_{\mathrm{ad}}(s)
&=\sum_{k,\ell\geq0}
\frac{\mathrm i^{\,k+\ell+1}s^{1,k}s^{2,\ell}}
{(k+\ell+1)k!\ell!}
+\sum_{\ell\geq0}
\frac{\mathrm i^{\,\ell}s^{2,\ell}}{(\ell+2)\ell!}.
\end{aligned}
\label{eq:intro-genus-zero-corrections}
\end{equation}
Here $\mathrm{i}$ denotes the imaginary unit.
\end{theorem}

\begin{remark}
The original  LYZZ conjecture for the diagonal case does not have a twist.
It might be also possible to encounter this twists in their framework by requiring some extra 
conditions on the integrable hierarchy in the higher genera case,
such as bihamiltonian structure.
Similar considerations have been used in the Hodge--GUE correspondence 
where on the Hodge hierarchy side some special choice has to been made to guarantee 
the existence of bihamiltonian structure.
\end{remark}

Since $\Lambda_0^\vee(1)=1$, the Hodge twist preserves the
genus-zero descendant theory and its $S$-calibration. Its effect
appears in the higher-genus reconstruction: the stationary-phase construction for the diagonal spectral curve produces an $R$-matrix differing from the geometric diagonal $R$-matrix by a scalar odd exponential. 
We identify this factor in \eqref{eq:di-geometric-R-correction} and show in Proposition~\ref{prop:Hodge-twisted-R} that it is precisely the transformation induced by the Hodge insertion above.

\begin{figure}

\begin{center}
\begin{tikzpicture}[
  >=Latex,
  font=\small,
  model/.style={
    align=center,
    text width=3.8cm
  },
  target/.style={
    align=center,
    text width=1.7cm,
    minimum height=0.6cm,
    inner sep=2pt
  },
  curve/.style={
    align=center,
    text width=2.5cm,
    minimum height=0.6cm,
    inner sep=2pt
  },
  flow/.style={
    ->,
    line width=0.7pt,
    shorten <=3pt,
    shorten >=3pt
  },
  flowlabel/.style={
    font=\scriptsize,
    align=center,
    inner sep=2pt,
    fill=white
  }
]
  \node[model] at (0,1.05)
    {{\bfseries A-model}\\[3pt]
     {\footnotesize Gromov--Witten theory}};

  \node[model] at (5.7,1.05)
    {{\bfseries B-model}\\[3pt]
     {\footnotesize Topological recursion}};

  \node[target] (P) at (0,-0.15)
    {$\bP^1$};

  \node[target] (X) at (0,-3.05)
    {$X^{\rm di},\;X^{\rm ad}$};

  \node[curve] (NS) at (5.7,-0.15)
    {Spectral curve\\of $\bP^1$};

  \node[curve] (Eff) at (5.7,-1.60)
    {Shifted curves};

  \node[curve] (Mirror) at (5.7,-3.05)
    {Spectral curves of $X^{\rm di}, \; X^{\rm ad}$};

  \draw[flow] (P.east) --
    node[flowlabel,above=3pt]
      {Norbury--Scott}
    (NS.west);

  \draw[flow] (X.east) --
    node[flowlabel,above=3pt]
      {Landau--Ginzburg}
    node[flowlabel,below=3pt] 
      {dimensional reduction}
    (Mirror.west);

  \draw[flow] (Mirror.north) --
    node[flowlabel,right=6pt,align=left]
      {Change $y$}
    (Eff.south);

  \draw[flow] (Eff.north) --
    node[flowlabel,right=6pt,align=left]
      {Compare $\omega_{g,n}$}
    (NS.south);
\end{tikzpicture}
\end{center}
\caption{Proof strategy of Theorem \ref{thm:intro-correspondences}. The diagonal spectral curve corresponds to the Hodge-twisted A-model theory.}\label{fig:proof-step}
\end{figure}
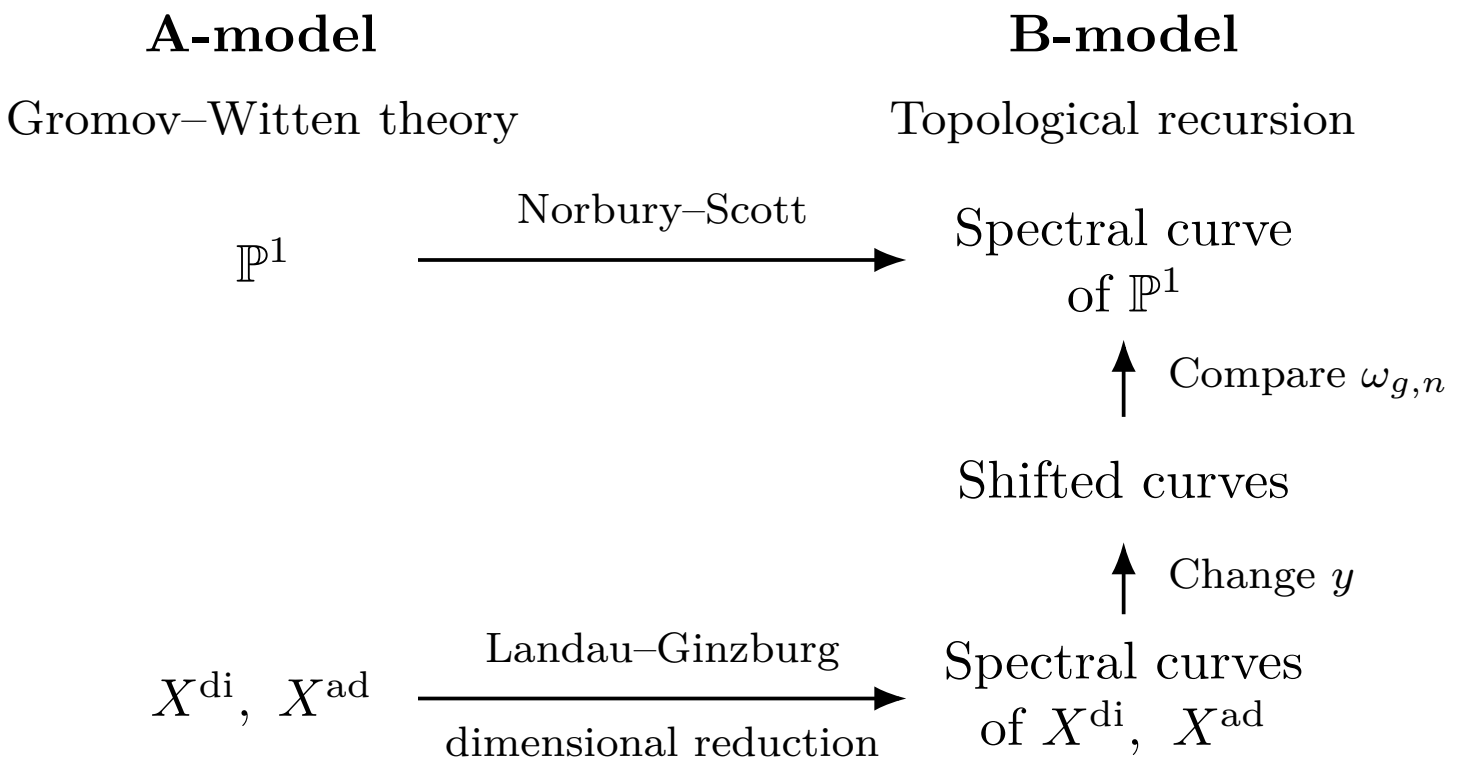

\medskip
\noindent\textbf{Outline of the proof.}
The proof compares the spectral curve realizations of the three theories and then recovers their descendant potentials (see Figure \ref{fig:proof-step}). 
It consists of three principal steps.

\medskip
\noindent\textit{1. The original spectral curves and their comparison
on a common surface.}
The one-dimensional spectral curve associated with the resolved conifold Frobenius structures are
\begin{equation}
\begin{aligned}
x^{X^{\mathrm{di}}}(u;p)
&=u^1+\log(1+p)
+\log\!\left(1+\frac{e^{u^2}}p\right),
&y^{X^{\mathrm{di}},\circ}(p)&=\log p,
\\
x^{X^{\mathrm{ad}}}(u;p)
&=u^1+\log(1+p)
-\log\!\left(1+\frac{e^{u^2}}p\right),
&y^{X^{\mathrm{ad}},\circ}(p)&=-\mathrm i\log p,
\end{aligned}
\label{eq:intro-conifold-spectral-curves}
\end{equation}
with Bergman kernel
\[
B(p_1,p_2)=\frac{dp_1\,dp_2}{(p_1-p_2)^2}.
\]
Here we use $u=(u^1,u^2)$ as the primary flat coordinates for the resolved conifold Frobenius manifolds.

The spectral curves at \eqref{eq:intro-conifold-spectral-curves} are obtained by dimensional reduction of the Landau-Ginzburg potentials.
Their critical values are the canonical coordinates of the corresponding Frobenius manifolds at the primary point $u=(u^1,u^2)$. 

We determine their $R$-matrices by normalized stationary-phase expansions in
Subsection~\ref{subsec:R-comparison}. The anti-diagonal spectral curve reconstructs
the geometric anti-diagonal theory, while the diagonal spectral curve
reconstructs the Hodge-twisted theory described above; see
Proposition~\ref{prop:three-reconstructions}.

Denote the primary flat coordinates of $\mathbb P^1$ by $v=(v^1,v^2)$. The Frobenius coordinate transformations are \cite{liuEquivariantGromovWitten2022}
\begin{equation}
\begin{aligned}
X^{\mathrm{di}}:\qquad
v^1&=e^{u^1}(1+e^{u^2}),
&v^2&=2u^1+u^2,
\\
X^{\mathrm{ad}}:\qquad
v^1&=\mathrm i\,e^{u^1}(1-2e^{u^2}),
&v^2&=2u^1+u^2+\log(1-e^{u^2}).
\end{aligned}
\label{eq:intro-Frobenius-coordinate-transformations}
\end{equation}
These are \eqref{eq:di-coordinate-transformation} and
\eqref{eq:ad-coordinate-transformation}. Indeed, write the $x$-function of $\mathbb P^1$ as 
\begin{equation}
x^{\mathbb{P}^1}(\ell)=\mathsf{x}(\ell)=v^1+\ell+\frac{e^{v^2}}{\ell}.
\label{eq:intro-common-projection}
\end{equation}
Then the changes of sphere coordinate
\begin{equation}
\begin{aligned}
X^{\mathrm{di}}:\qquad
\ell&=e^{u^1}p,
&\mathsf{x}(\ell)&=e^{x^{X^{\mathrm{di}}}(u;p)},
\\
X^{\mathrm{ad}}:\qquad
\ell&=\mathrm i\,e^{u^1}(p+e^{u^2}),
&\mathsf{x}(\ell)&=\mathrm i\,e^{x^{X^{\mathrm{ad}}}(u;p)}
\end{aligned}
\label{eq:intro-sphere-coordinate-transformations}
\end{equation}
give $dx^{X^{\mathrm{di}}}=d\log \mathsf{x}$ and
$dx^{X^{\mathrm{ad}}}=d\log \mathsf{x}$. Here $\ell$ is the common coordinate of the spectral curves. 

Both changes are affine in $p$ and preserve the Bergman kernel. 
Thus the three spectral curves are expressed on the common sphere $\Sigma=\mathbb P^1$ with coordinate $\ell$ and the same Bergman kernel $B$. Their $x$-functions are $\mathsf{x}$ for $\mathbb P^1$ and, up to additive constants, $\log \mathsf{x}$ for the two resolved conifold theories.

We then replace the original $y$-functions of conifolds:
\begin{equation}
  \begin{aligned}
    y^{X^{\mathrm{di}},\circ}(p)=\log p&\longrightarrow y^{X^{\mathrm{di}}}_{\mathrm{sh}}(\ell)=\mathsf{x}\log\ell,\\
    y^{X^{\mathrm{ad}},\circ}(p)=-\mathrm i\log p&\longrightarrow y^{X^{\mathrm{ad}}}_{\mathrm{sh}}(\ell)=-\mathsf{x}\log\ell.
  \end{aligned}
\end{equation}  
The subscript `$\mathrm{sh}$' indicates the corresponding shift
of the ancestor background.
The three spectral curve data to be compared become
\begin{equation}
\mathcal{C}^{\mathbb{P}^1}=(\mathsf{x},\log\ell,B),
\quad
\mathcal{C}^{\mathrm{di}}_{\mathrm{sh}}=(\log \mathsf{x},\mathsf{x}\log\ell,B),
\quad
\mathcal{C}^{\mathrm{ad}}_{\mathrm{sh}}=(\log \mathsf{x},-\mathsf{x}\log\ell,B),
\label{eq:intro-modified-spectral-curves}
\end{equation}
where
\[
B(\ell_1,\ell_2)
=\frac{d\ell_1\,d\ell_2}{(\ell_1-\ell_2)^2}.
\]
On the locus $(v^1)^2\ne4e^{v^2}$, the functions $\mathsf{x}$ and $\log \mathsf{x}$ have the same simple critical points and the same local involution
\[
\ell\longmapsto\frac{e^{v^2}}{\ell}.
\]
The factor $\mathsf{x}$ introduced in the shifted $y$-function cancels the
factor $\mathsf{x}^{-1}$ in $d\log \mathsf{x}$. Consequently, the first two curves in
\eqref{eq:intro-modified-spectral-curves} have identical Eynard--Orantin recursion
kernels. The additional minus sign in the third curve gives
\begin{equation}
\omega^+_{g,n}=\omega^{\mathbb P^1}_{g,n},
\qquad
\omega^-_{g,n}=(-1)^n\omega^{\mathbb P^1}_{g,n},
\qquad n\geq1,\quad 2g-2+n>0.
\label{eq:intro-recursion-comparison}
\end{equation}
Here $\omega^+_{g,n}$ and $\omega^-_{g,n}$ denote the multidifferentials
of the shifted diagonal and anti-diagonal curves.
These are identities of multidifferentials on the common surface.


The relations between these spectral curves are summarized in the following diagram:
\par\medskip

\begin{center}

\begin{tikzpicture}[
  >=Latex,
  font=\small,
  curve/.style={align=center,inner sep=3pt,text width=4.05cm},
  flow/.style={->,line width=.65pt,shorten <=4pt,shorten >=4pt},
  compare/.style={flow,dashed},
  shiftlabel/.style={font=\scriptsize,align=center,inner sep=2pt},
  maplabel/.style={font=\small,align=center,inner sep=3pt,fill=white}
]
\node[curve] (D0) at (-4.65,0)
  {Diagonal\\original curve\\[5pt]
   $\bigl(\log \mathsf{x},\,y^{X^{\mathrm{di}},\circ},\,B\bigr)$};
\node[curve] (A0) at (4.65,0)
  {Anti-diagonal\\original curve\\[5pt]
   $\bigl(\log \mathsf{x},\,y^{X^{\mathrm{ad}},\circ},\,B\bigr)$};
\node[curve] (P) at (0,0)
  {$\mathcal C^{\mathbb P^1}$\\[5pt]
   $\bigl(\mathsf{x},\,y^{\mathbb P^1,\circ},\,B\bigr)$\\[5pt]
   $y^{\mathbb P^1,\circ}=\log\ell$};
\node[curve] (D) at (-4.65,-3.6)
  {$\mathcal C_{\mathrm{sh}}^{\mathrm{di}}$\\[5pt]
   $\bigl(\log \mathsf{x},\,y_{\mathrm{sh}}^{\mathrm{di}},\,B\bigr)$\\[5pt]
   $y_{\mathrm{sh}}^{\mathrm{di}}=\mathsf{x}\log\ell$};
\node[curve] (A) at (4.65,-3.6)
  {$\mathcal C_{\mathrm{sh}}^{\mathrm{ad}}$\\[5pt]
   $\bigl(\log \mathsf{x},\,y_{\mathrm{sh}}^{\mathrm{ad}},\,B\bigr)$\\[5pt]
   $y_{\mathrm{sh}}^{\mathrm{ad}}=-\mathsf{x}\log\ell$};
\draw[flow] (D0.south) --
  node[shiftlabel,left=4pt] {$y$ shift\\[3pt]$x,B$ fixed}
  (D.north);
\draw[flow] (A0.south) --
  node[shiftlabel,right=4pt] {$y$ shift\\[3pt]$x,B$ fixed}
  (A.north);
\draw[compare] (P.south) -- (D.north);
\draw[compare] (P.south) -- (A.north);
\node[maplabel] at (-1.35,-2.35)
  {$x\mapsto\log x$\\[4pt]$y\mapsto xy$};
\node[maplabel] at (1.35,-2.35)
  {$x\mapsto\log x$\\[4pt]$y\mapsto-xy$};
\end{tikzpicture}


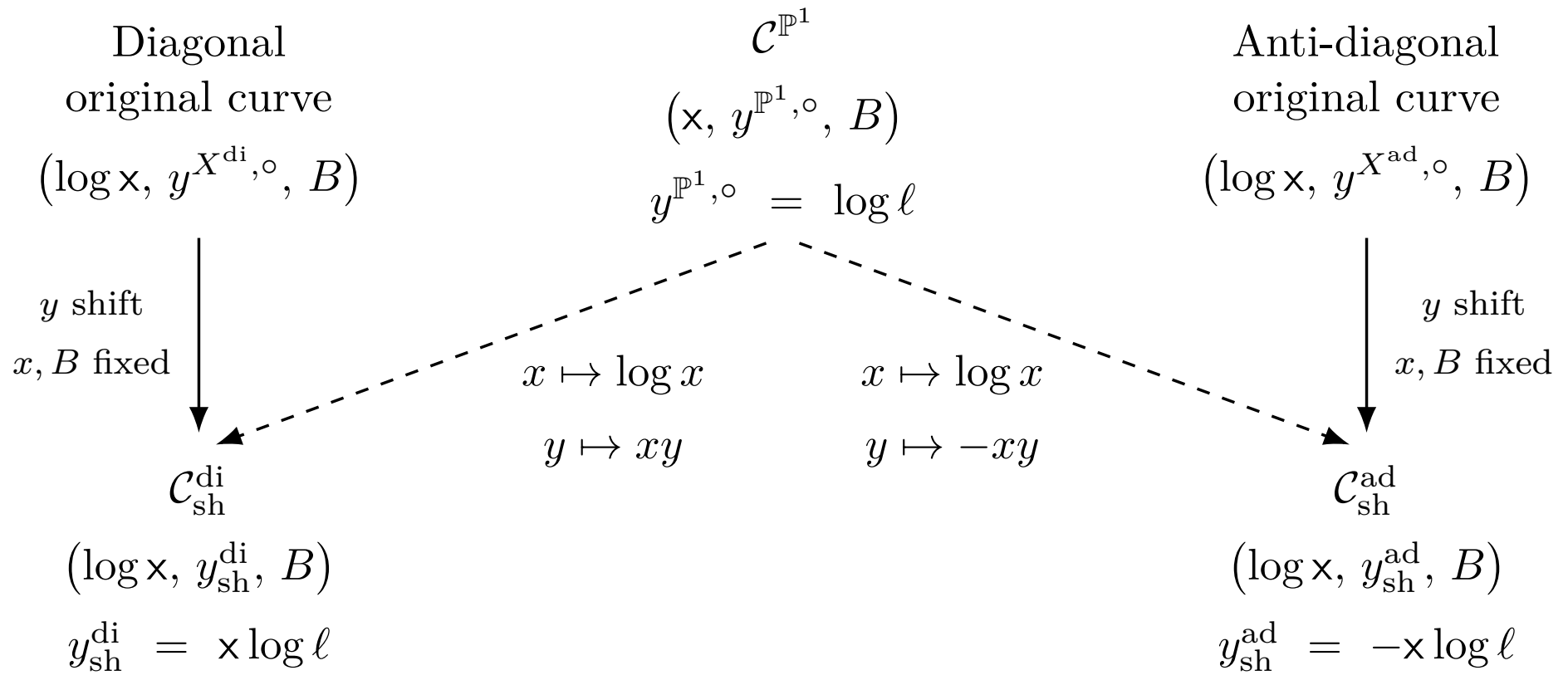
\captionof{figure}{Original curves and shifted curves.
The vertical arrows vary $y$ at fixed $x$ and $B$.
The slanted arrows describe the transformations used to compare
the multidifferentials; $B$ is unchanged throughout.}
\label{fig:ancestor-shifts-and-comparison}
\end{center}

\medskip
\noindent\textit{2. Interpretation of the shifted $y$-functions as
ancestor background shifts.}
The modification of $y$ has a precise meaning in CohFTs. The DOSS correspondence \cite{dunin-barkowskiIdentificationGiventalFormula2014} needs a $B$-compatible $y$-function $y^\circ$ with semisimple reconstruction. The extension developed in \cite{xiongEYNARDORANTINTOPOLOGICALRECURSION} describes how replacing the $B$-compatible $y^\circ$ by a general $y$, while keeping $x$, $B$, and the corresponding $R$-matrix fixed, changes the point at which the ancestor potential is evaluated.
Thus a modified $y$-function encodes a background shift of ancestor theory. 
In brief,
\begin{equation}
\omega_{g,n}[y^\circ]
\longleftrightarrow
\left.\partial^n\overline F_g(\bar t;v)\right|_{\bar t=0},
\qquad
\omega_{g,n}[y]
\longleftrightarrow
\left.\partial^n\overline F_g(\bar t;v)\right|_{\bar t=a}.
\end{equation}

At a fixed primary point $v$, set $s_v^{\alpha,k}=v^\alpha\delta_{k0}$.
The descendant transformations
\eqref{eq:intro-di-descendant-transformation} and
\eqref{eq:intro-ad-descendant-transformation} send $s_v$ to
$t_{\mathrm{di}}(s_v)$ and $t_{\mathrm{ad}}(s_v)$.
Proposition~\ref{prop:effective-background} identifies their ancestor
backgrounds at the corresponding primary points $u$ with the shifts
given by $\mathsf{x}\log\ell$ and $-\mathsf{x}\log\ell$, respectively.

\medskip
\noindent\textit{3. Extraction of descendant
coefficients.}
To pass from the equality of multidifferentials to descendant invariants, we define linear functionals $\lambda_{\alpha,k}$ on the space of one-forms on the common sphere. This approach is related to
the extraction of descendants by residues and regularized integrals
in Borot--Norbury \cite{borotLoopEquationsGromovWitten2019}.
In the present proof, the functionals are defined algebraically
from the $S$-calibration of $\mathbb P^1$ theory.

More precisely, let $W_j^{\mathbb P^1,i}$ be the ancestor one-forms 
of the $\mathbb{P}^1$ curve defined in Section \ref{sec:TR} and let $\Phi^{\mathbb P^1}$ be its transition matrix, both evaluated at the fixed primary point $v$. We define the functionals $\lambda_{\alpha,k}$ by
\begin{equation}
\lambda_{\alpha,k}\bigl(W_j^{\mathbb P^1,i}\bigr)
=\sum_\beta\bigl((\Phi^{\mathbb P^1}(v))^{-1}\bigr)^i{}_\beta
\bigl(S_{k-j}^{\mathbb P^1}(v)\bigr)^\beta{}_\alpha,
\qquad S_m^{\mathbb P^1}=0\quad(m<0).
\label{eq:intro-descendant-functionals}
\end{equation}
The linear independence of the ancestor one-forms
makes this prescription well defined; see
\eqref{eq:source-functional}. For $2g-2+n>0$, applying these
functionals to the recursion differentials gives
\begin{equation}
\left(\bigotimes_{r=1}^{n}\lambda_{\alpha_r,k_r}\right)
\omega^{\mathbb P^1}_{g,n}
=\left.
\frac{\partial^n F_g^{\mathbb P^1}(s)}
{\partial s^{\alpha_1,k_1}\cdots\partial s^{\alpha_n,k_n}}
\right|_{s=s_v}.
\label{eq:intro-descendant-extraction}
\end{equation}

\begin{figure}
\centering
\begin{tikzpicture}[
  >=Latex,
  font=\small,
  heading/.style={font=\small\bfseries,text=black},
  source/.style={align=center,text width=1.8cm,inner sep=3pt},
  result/.style={align=center,text width=4.6cm,inner sep=3pt},
  flow/.style={
    ->,draw=black,line width=0.8pt,
    shorten <=3pt,shorten >=3pt
  },
  tag/.style={font=\scriptsize,align=center,inner sep=3pt}
]
  \node[heading] at (0,2.15) {Multidifferentials};
  \node[heading] at (6.2,2.15) {Descendant invariants};

  \node[source] (AP) at (0,1.15)
    {$\displaystyle\omega_{g,n}^{\bP^1}$};
  \node[source] (AX) at (0,-1.65)
    {$\displaystyle\omega_{g,n}^{+},\;\omega_{g,n}^{-}$};

  \node[result] (DP) at (6.2,1.15)
    {$\displaystyle
      \left.
      \frac{\partial^n F_g^{\bP^1}}
      {\partial s^{\alpha_1,k_1}\cdots
       \partial s^{\alpha_n,k_n}}
      \right|_{s=s_v}$};

  \node[result] (DX) at (6.2,-1.65)
    {$\displaystyle
      \left.
      (\pm)^n\frac{\partial^n\!\left[F_g^{X}(t_\nu(s))\right]}
      {\partial s^{\alpha_1,k_1}\cdots
       \partial s^{\alpha_n,k_n}}
      \right|_{s=s_v}$};

  \draw[flow] (AP.east) --
    node[tag,above=4pt]
      {$\bigl(\Phi^{\bP^1}\bigr)^{-1}S^{\bP^1}$}
    node[tag,below=4pt]
      {$\bigotimes_{r=1}^{n}\lambda_{\alpha_r,k_r}$}
    (DP.west);

  \draw[flow] (AX.east) --
    node[tag,above=4pt]
      {$\bigotimes_{r=1}^{n}\lambda_{\alpha_r,k_r}$}
    (DX.west);

  \draw[flow,<->] (AP.south) --
    (AX.north);

  \draw[flow,<->] (DP.south) --
    node[tag,right=5pt,align=left]
      {LYZZ coordinate\\transformation}
    (DX.north);
\end{tikzpicture}
\caption{From multidifferentials to descendant Gromov--Witten invariants.
In the lower row $F_g^X$ denotes the Hodge-twisted diagonal potential
for the upper sign and the anti-diagonal potential for the lower sign.}\label{fig:form-descendant}
\end{figure}
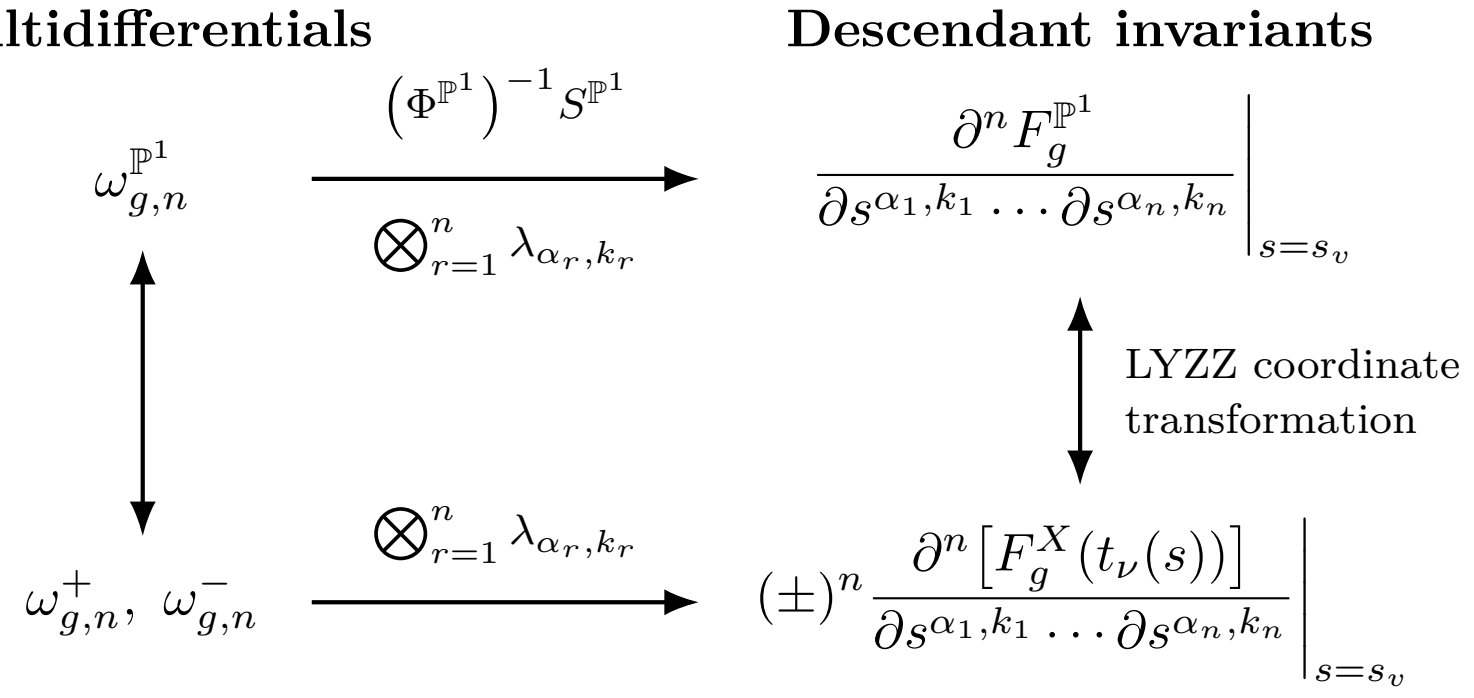

We compute the action of the same functionals on the resolved conifold
ancestor forms and show that the resulting coefficients reproduce
\eqref{eq:intro-di-descendant-transformation} and
\eqref{eq:intro-ad-descendant-transformation}. 
This proves the correspondence for positive-genus potentials up to constants. 
We fix these constants using the dilaton equation for $g\geq2$ and a direct computation of the genus-one ancestor vacuum contribution. 
Combining the result with the genus-zero identities already proved by LYZZ establishes Theorem~\ref{thm:intro-correspondences} (see Figure \ref{fig:form-descendant}).

\subsection*{Overview of the paper}
The remainder of this paper is organized as follows. 
Section~\ref{sec:geometry} reviews the Gromov--Witten theories of $\mathbb P^1$ and the resolved conifold.
Section~\ref{sec: B-theory} derives the spectral curves from
Landau--Ginzburg mirrors, determines the geometric $R$-matrices,
and introduces the Hodge twist in the diagonal theory.
Section~\ref{sec:TR} recalls topological recursion and the DOSS
correspondence, explains the effect of varying $y$, and establishes
the comparison of recursion differentials.
Section~\ref{sec:proof} identifies the ancestor backgrounds,
compares descendant invariants through the functionals
$\lambda_{\alpha,k}$, determines the vacuum terms, and proves
the main theorem.
Section~\ref{sec:outlook} summarizes the method and discusses possible extensions.

\subsection*{Acknowledgements}
The second-named author would like to thank Professor Zhengyu Zong for useful discussions.
The third-named author is partially supported by the Natural Science Foundation, China
(grant No.~12371254, No.~11890662, No.~12061131014).

\section{Gromov--Witten theories of \texorpdfstring{$\mathbb P^1$}{P1} and the resolved conifold}
\label{sec:geometry}

In this section, we describe the Gromov--Witten theories of $\mathbb P^1$
and the resolved conifold with anti-diagonal and diagonal torus actions.
We first introduce their descendant invariants and generating functions,
and then recall the associated Frobenius structures and $S$-calibrations.
These calibrations determine the transformations of descendant variables appearing
in the LYZZ correspondences. 
\subsection{The geometric theories and their potentials}
\label{subsec:geometric-theories}

Let
\[
X=\operatorname{Tot}\!\left(
\mathcal O_{\mathbb P^1}(-1)\oplus
\mathcal O_{\mathbb P^1}(-1)
\right)
\]
be the resolved conifold. We consider two actions of $T=\mathbb C^*$ on
$X$, both acting trivially on the base $\mathbb P^1$. They act on the
two line-bundle summands by
\begin{equation}
\begin{aligned}
X^{\mathrm{ad}}:\qquad
a\cdot(p,v_1,v_2)&=(p,av_1,a^{-1}v_2),\\
X^{\mathrm{di}}:\qquad
a\cdot(p,v_1,v_2)&=(p,av_1,av_2),
\end{aligned}
\qquad a\in\mathbb C^*,
\label{eq:conifold-actions}
\end{equation}
in local coordinates. 
Thus $X^{\mathrm{ad}}$ and $X^{\mathrm{di}}$ denote the same underlying
threefold equipped with the anti-diagonal and diagonal actions,
respectively. 

Let $H\in H^2(\mathbb P^1;\mathbb C)$ be the hyperplane class,
normalized by $\int_{\mathbb P^1}H=1$, and set
\[
\phi_1=1,\qquad \phi_2=H.
\]
We use the same notation for their pullbacks to $X$. Write $\lambda$
for the equivariant parameter corresponding to the standard character
of $T$. Since both actions fix the base, we have
\[
H_T^*(X^{\mathrm{ad}};\mathbb C)
\cong
H_T^*(X^{\mathrm{di}};\mathbb C)
\cong
\mathbb C[H,\lambda]/(H^2).
\]
The notation $\mathcal O_{\mathbb P^1}(-1)_w$ indicates that the line
bundle carries the fiber character $a\mapsto a^w$. In particular,
\begin{equation}
e_T\!\left(\mathcal O_{\mathbb P^1}(-1)_1\right)=\lambda-H,
\qquad
e_T\!\left(\mathcal O_{\mathbb P^1}(-1)_{-1}\right)=-\lambda-H.
\label{eq:equivariant-line-euler-classes}
\end{equation}
Throughout the paper, the equivariant resolved conifold invariants are specialized
to $\lambda=1$ after localization.

For nonnegative integers $g,n,d$, let
$\overline{\mathcal M}_{g,n}(\mathbb P^1,d)$ be the moduli space of
connected stable maps of genus $g$ and degree $d$, with $n$ marked points.
Denote the universal curve and universal map by
\[
\pi:\mathcal C\longrightarrow
\overline{\mathcal M}_{g,n}(\mathbb P^1,d),
\qquad
f:\mathcal C\longrightarrow\mathbb P^1.
\]
For each marking, let $\operatorname{ev}_i$ be the evaluation map and
let $\psi_i$ be the first Chern class of the cotangent line bundle.

The $T$-action on the resolved conifold $X$ induces a $T$-action on $\overline{\mathcal M}_{g,n}(X,d)$.
The Gromov-Witten invariants of $X$ are given by the twisted integrals over the virtual fundamental class $[\overline{\mathcal M}_{g,n}(\mathbb P^1,d)]^\vir$
\cite{graberLocalizationVirtualClasses1999}. 
Explicitly, for $k_i\geq0$
and $\alpha_i\in\{1,2\}$, the three sets of descendant invariants are
\begin{equation}
\begin{aligned}
\left\langle
\prod_{i=1}^{n}\tau_{k_i}(\phi_{\alpha_i})
\right\rangle^{\mathbb P^1}_{g,d}
&=
\int_{[\overline{\mathcal M}_{g,n}(\mathbb P^1,d)]^{\mathrm{vir}}}
\prod_{i=1}^{n}\psi_i^{k_i}
\operatorname{ev}_i^*(\phi_{\alpha_i}),
\\[6pt]
\left\langle
\prod_{i=1}^{n}\tau_{k_i}(\phi_{\alpha_i})
\right\rangle^{X^{\mathrm{ad}}}_{g,d}
&=
\Biggl[
\int_{[\overline{\mathcal M}_{g,n}(\mathbb P^1,d)]^{\mathrm{vir}}}
e_T\left(
-R^\bullet\pi_*f^*
\bigl(
\mathcal O_{\mathbb P^1}(-1)_1
\oplus
\mathcal O_{\mathbb P^1}(-1)_{-1}
\bigr)
\right)
\\
&\qquad\times
\prod_{i=1}^{n}\psi_i^{k_i}
\operatorname{ev}_i^*(\phi_{\alpha_i})
\Biggr]_{\lambda=1},
\\[6pt]
\left\langle
\prod_{i=1}^{n}\tau_{k_i}(\phi_{\alpha_i})
\right\rangle^{X^{\mathrm{di}}}_{g,d}
&=
\Biggl[
\int_{[\overline{\mathcal M}_{g,n}(\mathbb P^1,d)]^{\mathrm{vir}}}
e_T\left(
-R^\bullet\pi_*f^*(
\mathcal O_{\mathbb P^1}(-1)_1
\oplus 
\mathcal O_{\mathbb P^1}(-1)_1)
\right)
\\
&\qquad \times 
\prod_{i=1}^{n}\psi_i^{k_i}
\operatorname{ev}_i^*(\phi_{\alpha_i})
\Biggr]_{\lambda=1}.
\end{aligned}
\label{eq:geometric-correlators}
\end{equation}
Here $R^\bullet\pi_*=R^0\pi_*-R^1\pi_*$, and the equivariant Euler
class is extended multiplicatively to virtual bundles. The integrals
for $X^{\mathrm{ad}}$ and $X^{\mathrm{di}}$ are taken in localized
equivariant cohomology before the specialization $\lambda=1$.
For $d=0$, the moduli space is defined only when $2g-2+n>0$;
the remaining degree-zero terms are assigned the value zero in the
generating functions below.

Introduce descendant variables
\[
s=(s^{\alpha,k})_{\alpha=1,2;\,k\geq0},
\qquad
t=(t^{\alpha,k})_{\alpha=1,2;\,k\geq0},
\]
for $\mathbb P^1$ and the resolved conifold, respectively. 
The genus-$g$ descendant potentials are
\begin{equation}
\begin{aligned}
F_g^{\mathbb P^1}(s)
&=
\sum_{d,n\geq0}\frac{1}{n!}
\sum_{\substack{k_1,\ldots,k_n\geq0\\
\alpha_1,\ldots,\alpha_n\in\{1,2\}}}
s^{\alpha_1,k_1}\cdots s^{\alpha_n,k_n}\times
\left\langle
\prod_{i=1}^{n}\tau_{k_i}(\phi_{\alpha_i})
\right\rangle^{\mathbb P^1}_{g,d},
\\[4pt]
F_g^{X^{\mathrm{ad}}}(t)
&=
\sum_{d,n\geq0}\frac{1}{n!}
\sum_{\substack{k_1,\ldots,k_n\geq0\\
\alpha_1,\ldots,\alpha_n\in\{1,2\}}}
t^{\alpha_1,k_1}\cdots t^{\alpha_n,k_n}\times
\left\langle
\prod_{i=1}^{n}\tau_{k_i}(\phi_{\alpha_i})
\right\rangle^{X^{\mathrm{ad}}}_{g,d},
\\[4pt]
F_g^{X^{\mathrm{di}}}(t)
&=
\sum_{d,n\geq0}\frac{1}{n!}
\sum_{\substack{k_1,\ldots,k_n\geq0\\
\alpha_1,\ldots,\alpha_n\in\{1,2\}}}
t^{\alpha_1,k_1}\cdots t^{\alpha_n,k_n}\times
\left\langle
\prod_{i=1}^{n}\tau_{k_i}(\phi_{\alpha_i})
\right\rangle^{X^{\mathrm{di}}}_{g,d}.
\end{aligned}
\label{eq:genus-descendant-potentials}
\end{equation}

The total descendant potentials, or partition functions, are
\begin{equation}
\begin{aligned}
\mathcal D_{\mathbb P^1}(s;\hbar)
&=
\exp\left(
\sum_{g\geq0}\hbar^{g-1}F_g^{\mathbb P^1}(s)
\right),\\
\mathcal D_{X^{\mathrm{ad}}}(t;\hbar)
&=
\exp\left(
\sum_{g\geq0}\hbar^{g-1}F_g^{X^{\mathrm{ad}}}(t)
\right),\\
\mathcal D_{X^{\mathrm{di}}}(t;\hbar)
&=
\exp\left(
\sum_{g\geq0}\hbar^{g-1}F_g^{X^{\mathrm{di}}}(t)
\right),
\end{aligned}
\label{eq:descendant-potentials}
\end{equation}
where $\hbar$ is a formal parameter keeping track of the genus. The
genus-zero restrictions of these potentials to primary parameters determine the Frobenius structures described in the next subsection.


\subsection{Frobenius manifolds}
\label{subsec:Frobenius}

The genus-zero primary potentials define Frobenius structures on the
corresponding parameter spaces \cite{dijkgraafMeanFieldTheory1990,dubrovinGeometry2DTopological1996,kontsevichGromovWittenClassesQuantum1994}.
It is the Frobenius structure of the equivariant quantum cohomology ring of $\bP^1$ and the resolved conifold.
We use flat coordinates $v^1,v^2$ for $\mathbb P^1$ and $u^1,u^2$
for the resolved conifold, corresponding to the basis $\phi_1=1$,
$\phi_2=H$. Thus the primary restrictions are obtained by setting
\[
s^{\alpha,0}=v^\alpha,\qquad
t^{\alpha,0}=u^\alpha,\qquad
s^{\alpha,k}=t^{\alpha,k}=0
\quad(k>0).
\]
We write $\partial_\alpha$ for the coordinate vector field corresponding
to $\phi_\alpha$. In each theory, $\partial_1$ is the unit vector field.

The flat metrics, in the ordered basis $(\partial_1,\partial_2)$, are
\begin{equation}
\eta^{\mathbb P^1}
=
\begin{pmatrix}
0&1\\
1&0
\end{pmatrix},
\qquad
\eta^{X^{\mathrm{ad}}}
=
\begin{pmatrix}
0&-1\\
-1&0
\end{pmatrix},
\qquad
\eta^{X^{\mathrm{di}}}
=
\begin{pmatrix}
2&1\\
1&0
\end{pmatrix}.
\label{eq:flat-metrics}
\end{equation}

The genus-zero primary potentials are 
\begin{equation}
\begin{aligned}
F_0^{\mathbb P^1}(v^1,v^2)
&=\frac12(v^1)^2v^2+e^{v^2},\\
F_0^{X^{\mathrm{ad}}}(u^1,u^2)
&=-\frac12(u^1)^2u^2+\operatorname{Li}_3(e^{u^2}),\\
F_0^{X^{\mathrm{di}}}(u^1,u^2)
&=\frac13(u^1)^3+\frac12(u^1)^2u^2
+\operatorname{Li}_3(e^{u^2}).
\end{aligned}
\label{eq:primary-potentials}
\end{equation}
Here
\[
\operatorname{Li}_3(e^{u^2})
=
\sum_{d=1}^{\infty}\frac{e^{du^2}}{d^3}.
\]
See e.g. \cite{dijkgraafMeanFieldTheory1990,briniLocalGromovWitten2012}.
For the resolved conifold, the degree sums are initially understood in the region
$\operatorname{Re}u^2<0$. We subsequently use their local analytic
continuations on simply connected domains avoiding $e^{u^2}=1$.

The potentials are also
described in \cite{liuEquivariantGromovWitten2022}, with a different
overall sign in the anti-diagonal case. Our anti-diagonal potential
and metric are the negatives of those used there; the convention here
follows directly from the localization definition in
Section~\ref{subsec:geometric-theories}.
It matches with the convention used in \cite{briniLocalGromovWitten2012}.

In each theory, the quantum product is characterized by
\begin{equation}
\eta\!\left(
\partial_\alpha\star\partial_\beta,\partial_\gamma
\right)
=
\partial_\alpha\partial_\beta\partial_\gamma F_0,
\qquad
\alpha,\beta,\gamma\in\{1,2\},
\label{eq:quantum-product-definition}
\end{equation}
where $F_0$ and $\eta$ are the corresponding potential and metric.
Note $\partial_1$ is the unit, and
\begin{equation}
\begin{aligned}
\partial_2\star_{\mathbb P^1}\partial_2
&=e^{v^2}\partial_1,\\
\partial_2\star_{X^{\mathrm{ad}}}\partial_2
&=-\frac{e^{u^2}}{1-e^{u^2}}\partial_1,\\
\partial_2\star_{X^{\mathrm{di}}}\partial_2
&=\frac{e^{u^2}}{1-e^{u^2}}
\bigl(\partial_1-2\partial_2\bigr).
\end{aligned}
\label{eq:quantum-products}
\end{equation}
These multiplication tables show that the three structures are
semisimple at finite flat coordinates on the domains under consideration.

We now describe the canonical coordinates and normalized idempotents
needed for the $R$-matrix construction. Label the two canonical
directions by $+$ and $-$. In each theory, let
\begin{equation}
e_\pm=\frac{\partial}{\partial U_\pm},
\qquad
e_i\star e_j=\delta_{ij}e_i,
\qquad
\eta(e_i,e_j)=\frac{\delta_{ij}}{\Delta_i},
\qquad i,j\in\{+,-\}.
\label{eq:canonical-idempotents}
\end{equation}
We use branches for which $\sqrt{1-e^{u^2}}\to1$,
$\arcsin(e^{u^2/2})\to0$, and $\log(1\pm e^{u^2/2})\to0$
as $\operatorname{Re}u^2\to-\infty$, and keep all branch and sign
choices fixed throughout the paper.

For $\mathbb P^1$, one obtains
\begin{equation}
\begin{aligned}
U_\pm^{\mathbb P^1}
&=v^1\pm2e^{v^2/2},\\
e_\pm^{\mathbb P^1}
&=\frac12\partial_1
\pm\frac12e^{-v^2/2}\partial_2,\\
\Delta_\pm^{\mathbb P^1}
&=\pm2e^{v^2/2}.
\end{aligned}
\label{eq:P1-canonical-data}
\end{equation}

For the anti-diagonal resolved conifold, the corresponding expressions are
\begin{equation}
\begin{aligned}
U_\pm^{X^{\mathrm{ad}}}
&=u^1\pm2\mathrm i\,\arcsin(e^{u^2/2}),\\
e_\pm^{X^{\mathrm{ad}}}
&=\frac12\partial_1
\mp
\frac{\mathrm i\sqrt{1-e^{u^2}}}{2e^{u^2/2}}
\partial_2,\\
\Delta_\pm^{X^{\mathrm{ad}}}
&=\mp
\frac{2\mathrm i\,e^{u^2/2}}{\sqrt{1-e^{u^2}}}.
\end{aligned}
\label{eq:ad-canonical-data}
\end{equation}

For the diagonal resolved conifold, one has
\begin{equation}
\begin{aligned}
U_\pm^{X^{\mathrm{di}}}
&=u^1+2\log(1\pm e^{u^2/2}),\\
e_\pm^{X^{\mathrm{di}}}
&=\frac{1\pm e^{u^2/2}}{2}\partial_1
\pm
\frac{1-e^{u^2}}{2e^{u^2/2}}\partial_2,\\
\Delta_\pm^{X^{\mathrm{di}}}
&=\frac{\pm2e^{u^2/2}}
{(1\pm e^{u^2/2})^2}.
\end{aligned}
\label{eq:di-canonical-data}
\end{equation}
The idempotent formulas in
\eqref{eq:P1-canonical-data}--\eqref{eq:di-canonical-data} follow
directly from the multiplication table \eqref{eq:quantum-products};
applying the metrics \eqref{eq:flat-metrics} then gives the values of $\Delta_\pm$.

Finally, in each theory let $\Phi$ be the matrix whose columns are
the normalized idempotents in the flat basis. From now on, we use
the Einstein summation convention for repeated upper and lower indices.
Thus
\begin{equation}
\sqrt{\Delta_i}\,e_i
=
\Phi^\alpha{}_i\,\partial_\alpha,
\qquad i\in\{+,-\}.
\label{eq:normalized-idempotent-matrix}
\end{equation}
Explicitly, in the ordered flat basis $(\partial_1,\partial_2)$
and with columns ordered by $+,-$, these matrices are
\begin{equation}
\begin{aligned}
\Phi^{\mathbb P^1}
&=
\frac{1}{2}
\begin{pmatrix}
\sqrt{\Delta_+^{\mathbb P^1}}
&
\sqrt{\Delta_-^{\mathbb P^1}}
\\[2pt]
e^{-v^2/2}\sqrt{\Delta_+^{\mathbb P^1}}
&
-e^{-v^2/2}\sqrt{\Delta_-^{\mathbb P^1}}
\end{pmatrix},
\\[6pt]
\Phi^{X^{\mathrm{ad}}}
&=
\frac{1}{2}
\begin{pmatrix}
\sqrt{\Delta_+^{X^{\mathrm{ad}}}}
&
\sqrt{\Delta_-^{X^{\mathrm{ad}}}}
\\[2pt]
-\dfrac{i\sqrt{1-e^{u^2}}}{e^{u^2/2}}
 \sqrt{\Delta_+^{X^{\mathrm{ad}}}}
&
\dfrac{i\sqrt{1-e^{u^2}}}{e^{u^2/2}}
 \sqrt{\Delta_-^{X^{\mathrm{ad}}}}
\end{pmatrix},
\\[6pt]
\Phi^{X^{\mathrm{di}}}
&=
\frac{1}{2}
\begin{pmatrix}
(1+e^{u^2/2})\sqrt{\Delta_+^{X^{\mathrm{di}}}}
&
(1-e^{u^2/2})\sqrt{\Delta_-^{X^{\mathrm{di}}}}
\\[2pt]
\dfrac{1-e^{u^2}}{e^{u^2/2}}
 \sqrt{\Delta_+^{X^{\mathrm{di}}}}
&
-\dfrac{1-e^{u^2}}{e^{u^2/2}}
 \sqrt{\Delta_-^{X^{\mathrm{di}}}}
\end{pmatrix}.
\end{aligned}
\label{eq:Phi-matrices}
\end{equation}
Here the factors $\Delta_\pm$ are given in
(2.11)--(2.13), and their square roots are chosen locally
and kept fixed throughout.
This matrix satisfies
\begin{equation}
\Phi^{T}\eta\Phi=I,
\qquad
\Phi^{-1}=\Phi^{T}\eta,
\qquad
\Phi^{-1}
\begin{pmatrix}1\\0\end{pmatrix}
=
\begin{pmatrix}
\Delta_+^{-1/2}\\
\Delta_-^{-1/2}
\end{pmatrix}.
\label{eq:canonical-frame}
\end{equation}
$\Phi$ converts components in the normalized canonical basis
to components in the flat basis. 


\subsection{Cohomological field theory}\label{subsec:cohft}
In this subsection, we review the cohomological field theory (CohFT) introduced by
Kontsevich and Manin \cite{kontsevichGromovWittenClassesQuantum1994}.
Let $\cH$ be a complex vector space of dimension $N$ with a non-degenerate symmetric bilinear form $\eta$.
Let $\mathbb{A}$ be a $\bC$-algebra. A CohFT consists of $\{\Omega_{g,n}\}_{2g-2+n>0}$ on
$(\cH,\eta)$ where 
$$
  \Omega_{g,n}: \cH^{\otimes n}\rightarrow \mathbb{A}\otimes H^*(\overline{\cM}_{g,n},\bQ)
$$
satisfying 
\begin{itemize}
  \item[(1)]$\Omega_{g,n}(\gamma_1,\ldots,\gamma_n)$ is invariant under permutation of $\gamma_1,\ldots,\gamma_n\in \cH$.
  \item[(2)]The pull-backs $q^*\Omega_{g,n}$ and $r^*\Omega_{g,n}$ of the gluing maps 
  \begin{align*}
      &q: \overline{\cM}_{g-1,n+2}\rightarrow\overline{\cM}_{g,n}
      \\
      &r:\overline{\cM}_{g_1,n_1+1}\times\overline{\cM}_{g_2,n_2+1}\rightarrow \overline{\cM}_{g_1+g_2,n_1+n_2}
  \end{align*}
  are equal to the contraction with the bivector $\sum_{\alpha}\partial_\alpha\otimes \partial^\alpha$. Here $\{\partial_\alpha\}_{\alpha=1}^N$ is a flat basis of $\cH$, and
  $\{\partial^\alpha\}_{\alpha=1}^N$ is its dual basis with respect to $\eta$. 
  \item[(3)]There exists a flat unit ${\bf 1}\in\cH$ such that 
      $$
        \pi^*\Omega_{g,n}(\gamma_1,\ldots,\gamma_n) = \Omega_{g,n+1}(\gamma_1,\ldots,\gamma_n,{\bf 1}),
      $$
      where $\pi$ is the forgetful map $\pi:\overline{\mathcal M}_{g,n+1}\rightarrow\overline{\mathcal M}_{g,n}$ forgetting the last marked point.
\end{itemize}

Let $v\in \cH$, we define the shifted CohFT $\Omega_{g,n;v}$
$$
    \Omega_{g,n;v}(\gamma_1,\ldots,\gamma_n) := \sum_{m\ge0}\frac{1}{m!}(\pi_{m})_*\Omega_{g,n+m}(\gamma_1,\ldots,\gamma_n,v,\dots,v),
$$ 
where $(\pi_{m})_*$ is the push-forward of the forgetful map $\pi_{m}:\overline{\mathcal M}_{g,n+m}\rightarrow\overline{\mathcal M}_{g,n}$
forgetting the last $m$ points.
The ancestor correlators of the shifted CohFT $\Omega_{g,n;v}$ are
\begin{equation}
\left\langle
\prod_{r=1}^n\bar\tau_{k_r}(\gamma_r)
\right\rangle_{g,n;v}
=
\int_{\overline{\mathcal M}_{g,n}}
\Omega_{g,n;v}(\gamma_1,\ldots,\gamma_n)
\prod_{r=1}^n\psi_r^{k_r},
\qquad 2g-2+n>0.
\label{eq:ancestor-correlators}
\end{equation}

Let $\bar{t}(z) = \sum_{\alpha,k}\bar{t}^{\alpha,k}z^k\partial_\alpha $, we define the ancestor potential $\overline F_g(\bar t;v)$ as 
\begin{align}\label{eq:ancestor-invariants}
    \overline F_g(\bar t;v) = \sum_{n\geq0}\frac{1}{n!}\sum_{\substack{k_1,\ldots,k_n\ge0\\\alpha_1,\ldots,\alpha_n}}\bar{t}^{\alpha_1,k_1}\ldots\bar{t}^{\alpha_n,k_n}
    \times\left\langle
\prod_{r=1}^n\bar\tau_{k_r}(\partial_{\alpha_r})
\right\rangle_{g,n;v}.
\end{align}

In our case, we take $\cH$ as the Frobenius structure, i.e. the quantum cohomology ring, in Section \ref{subsec:Frobenius}.
Let $\ev:\overline{\cM}_{g,n}(\bP^1,d)\rightarrow (\bP^1)^n$ and $p:\overline{\cM}_{g,n}(\bP^1,d)\rightarrow \overline\cM_{g,n}$. 
The Gromov-Witten theory of $\bP^1$ and the resolved conifold in Section \ref{subsec:geometric-theories} is captured by the CohFT:
\begin{align*}
    \Omega^{\bP^1}_{g,n}(\gamma_1,\ldots,\gamma_n) := &\sum_{d\geq0} \mathfrak{Q}^d\cdot p_*\Big([\overline{\cM}_{g,n}(\bP^1,d)]^\vir\cap\ev^*(\gamma_1\otimes\ldots\otimes\gamma_n)\Big),
    \\
    \Omega^{X^{\rm ad}}_{g,n}(\gamma_1,\ldots,\gamma_n) :=& \sum_{d\geq0} \mathfrak{Q}^d\cdot p_*\Big(e_T\left(
    -R^\bullet\pi_*f^*
    \bigl(
    \mathcal O_{\mathbb P^1}(-1)_1
    \oplus
    \mathcal O_{\mathbb P^1}(-1)_{-1}
    \bigr)
    \right)\\
    &\qquad\qquad\qquad\cap[\overline{\cM}_{g,n}(\bP^1,d)]^\vir\cap\ev^*(\gamma_1\otimes\ldots\otimes\gamma_n)\Big),
    \\
    \Omega^{X^{\rm di}}_{g,n}(\gamma_1,\ldots,\gamma_n) :=& \sum_{d\geq0} \mathfrak{Q}^d\cdot p_*\Big(e_T\left(
    -R^\bullet\pi_*f^*
    \bigl(
    \mathcal O_{\mathbb P^1}(-1)_1
    \oplus
    \mathcal O_{\mathbb P^1}(-1)_{1}
    \bigr)
    \right)
    \\
    &\qquad\qquad\qquad\cap[\overline{\cM}_{g,n}(\bP^1,d)]^\vir\cap\ev^*(\gamma_1\otimes\ldots\otimes\gamma_n)\Big).
\end{align*}
Here $\fQ$ is the Novikov variable, we will use $\fQ=1$ in the following paragraph. We remind the reader that the ancestor potential \eqref{eq:ancestor-invariants} is 
different from the descendant potentials in \eqref{eq:genus-descendant-potentials}. The relation between ancestor and descendant potential is given by the $S$-calibration (see Section \ref{subsec:backgrounds}).


\subsection{\texorpdfstring{$S$}{S}-calibrations and the descendant transformations}
\label{subsec:S-calibration}

The genus-zero descendant theory determines an $S$-calibration consisting
of tangent-vector-valued formal series in $z^{-1}$ \cite{dubrovinGeometry2DTopological1996,giventalGromovWittenInvariantsQuantization2001}:
\[
S_\alpha(z)
=\partial_\alpha
+\sum_{k\geq1}z^{-k}\,(S_k)^\beta{}_\alpha\,\partial_\beta,
\qquad \alpha=1,2.
\]
Writing these series as columns in the flat basis gives the matrix
\[
S(z)=I+\sum_{k\geq1}S_kz^{-k}.
\]
It satisfies the quantum differential equations (QDE) and the symplectic
condition \cite{dubrovinGeometry2DTopological1996,giventalGromovWittenInvariantsQuantization2001}
\begin{equation}
z\partial_\alpha S(z)
=
(\partial_\alpha\star)S(z),
\qquad
S(-z)^*S(z)=I.
\label{eq:S-equations}
\end{equation}
Here the adjoint $^*$ is taken with respect to the metric of the corresponding
theory.

It is convenient to express these matrices through $\theta$-calibration functions:
\begin{equation}
\theta_\alpha(z)
=
\sum_{k\geq0}\theta_{\alpha,k}z^k,
\qquad
S^\alpha{}_\beta(z)
=
\eta^{\alpha\gamma}
\partial_\gamma\theta_\beta(z^{-1}).
\label{eq:scalar-matrix-calibration}
\end{equation}
The dependence on the primary variables is written explicitly below.
The scalar calibrations $\theta$ used here can be found in
\cite[Section 2]{liuEquivariantGromovWitten2022}, with the anti-diagonal
sign adjusted to the convention of Section~\ref{subsec:Frobenius}.

For $\mathbb P^1$, they are
\begin{equation}
\begin{aligned}
\theta_1^{\mathbb P^1}(v;z)
&=
e^{zv^1}
\sum_{d\geq0}
\left(
v^2-2\sum_{j=1}^{d}\frac1j
\right)
\frac{e^{dv^2}z^{2d}}{(d!)^2},
\\
\theta_2^{\mathbb P^1}(v;z)
&=
\frac1z
\left[
e^{zv^1}
\sum_{d\geq0}
\frac{e^{dv^2}z^{2d}}{(d!)^2}
-1
\right].
\end{aligned}
\label{eq:P1-calibration}
\end{equation}

For the anti-diagonal resolved conifold, one has
\begin{equation}
\begin{aligned}
\theta_1^{X^{\mathrm{ad}}}(u;z)
&=
-e^{zu^1}
\sum_{d\geq0}
\frac{(z)_d(-z)_d}{(d!)^2}e^{du^2}
\left(
u^2+
2\sum_{j=0}^{d-1}
\frac{j+z^2}{(j+1)(j^2-z^2)}
\right),
\\
\theta_2^{X^{\mathrm{ad}}}(u;z)
&=
-\frac1z
\left[
e^{zu^1}
\sum_{d\geq0}
\frac{(z)_d(-z)_d}{(d!)^2}e^{du^2}
-1
\right].
\end{aligned}
\label{eq:ad-calibration}
\end{equation}
For the diagonal resolved conifold, the corresponding formulas are
\begin{equation}
\begin{aligned}
\theta_1^{X^{\mathrm{di}}}(u;z)
&=
e^{zu^1}
\sum_{d\geq0}
\frac{(-z)_d^2}{(d!)^2}e^{du^2}
\left(
u^2+\frac2z
-
2\sum_{j=1}^{d}
\frac{z+1}{j(z-j+1)}
\right)
-\frac2z,
\\
\theta_2^{X^{\mathrm{di}}}(u;z)
&=
\frac1z
\left[
e^{zu^1}
\sum_{d\geq0}
\frac{(-z)_d^2}{(d!)^2}e^{du^2}
-1
\right].
\end{aligned}
\label{eq:di-calibration}
\end{equation}
Here $(a)_d=a(a+1)\cdots(a+d-1)$, with $(a)_0=1$. The apparent poles in these expressions cancel; in particular, the functions are regular at $z=0$.

The leading terms are
\[
\begin{aligned}
\bigl(\theta_{1,0}^{\mathbb P^1},\theta_{2,0}^{\mathbb P^1}\bigr)
&=(v^2,v^1),\\
\bigl(\theta_{1,0}^{X^{\mathrm{ad}}},\theta_{2,0}^{X^{\mathrm{ad}}}\bigr)
&=(-u^2,-u^1),\\
\bigl(\theta_{1,0}^{X^{\mathrm{di}}},\theta_{2,0}^{X^{\mathrm{di}}}\bigr)
&=(2u^1+u^2,u^1),
\end{aligned}
\]
so \eqref{eq:scalar-matrix-calibration} gives $S(z)=I+O(z^{-1})$ in
all three cases. 
In the anti-diagonal theory, these signs cancel in \eqref{eq:scalar-matrix-calibration}, leaving the $S$-matrix unchanged.

As explained in \cite{liuEquivariantGromovWitten2022},
the above three Frobenius manifolds are related to each other by some transformations of 
Frobenius manifolds introduced by Dubrovin.
The $\bP^1$-Frobenius manifold is related to the Frobenius manifold with potential $F_0^{X^{\rm di}}$
by almost duality \cite{dubrovinAlmostDualityFrobenius2004},
the latter is related to the Frobenius manifold with potential $F_0^{X^{\rm ad}}$ by 
a generalized Legendre transformation \cite{strachanGeneralizedLegendre2017}. 
Originally, Dubrovin \cite{dubrovinGeometry2DTopological1996} introduced the Legendre transformation.
This led to the generalized Legendre transformation \cite{strachanGeneralizedLegendre2017}
and was further developed in \cite{liuLegendreTransformationsClass2025}. 
Alternatively, 
the $\bP^1$-Frobenius manifold 
is related to the Frobenius manifold related to the Ablowitz-Ladik hierarchy via generalized Frobenius manifold \cite{liuEquivariantGromovWitten2022},
and the latter is almost dual to the  Frobenius manifold with potential $F_0^{X^{\rm ad}}$ 
as noted by \cite{briniLocalGromovWitten2012}.
Using such transformations,
one can write down the following changes of coordinates \cite[Section 3]{liuEquivariantGromovWitten2022}.
For the diagonal case, 
the flat coordinates and the calibration are changed to the $\bP^1$-case via
\begin{equation}
v^1=e^{u^1}(1+e^{u^2}),
\qquad
v^2=2u^1+u^2,
\label{eq:di-coordinate-transformation}
\end{equation}
and
\begin{equation}
\theta_{\alpha,k}^{\mathbb P^1}(v)
=
b_{\alpha,k}^{\beta,m,\mathrm{di}}
\theta_{\beta,m}^{X^{\mathrm{di}}}(u)
+
\frac{\delta_{\alpha2}}{(k+1)!}.
\label{eq:di-calibration-transformation}
\end{equation}
For the anti-diagonal theory, the coordinate transformation is
\begin{equation}
v^1=\mathrm i\,e^{u^1}(1-2e^{u^2}),
\qquad
v^2=2u^1+u^2+\log(1-e^{u^2}).
\label{eq:ad-coordinate-transformation}
\end{equation}
The corresponding identity is expressed in terms of gradients:
\begin{equation}
\begin{pmatrix}
\displaystyle
\frac{\partial\theta_{\alpha,k}^{\mathbb P^1}}{\partial v^1}(v)
\\[5pt]
\displaystyle
\frac{\partial\theta_{\alpha,k}^{\mathbb P^1}}{\partial v^2}(v)
\end{pmatrix}
=
-
b_{\alpha,k}^{\beta,m,\mathrm{ad}}
\begin{pmatrix}
e^{-u^1}&-e^{-u^1}\\
\mathrm i\,e^{u^2}&\mathrm i(1-e^{u^2})
\end{pmatrix}
\begin{pmatrix}
\displaystyle
\frac{\partial\theta_{\beta,m}^{X^{\mathrm{ad}}}}{\partial u^1}(u)
\\[5pt]
\displaystyle
\frac{\partial\theta_{\beta,m}^{X^{\mathrm{ad}}}}{\partial u^2}(u)
\end{pmatrix}.
\label{eq:ad-calibration-transformation}
\end{equation}
In the above we use the coefficients in \cite{liuEquivariantGromovWitten2022}:
\begin{equation}
\begin{aligned}
b_{1,k}^{\alpha,m,\mathrm{di}}
&=
\frac{k^m}{k!}
\left(
\delta^\alpha_1+2m\delta^\alpha_2
\right),
&
b_{2,k}^{\alpha,m,\mathrm{di}}
&=
\frac{(k+1)^m}{k!}\delta^\alpha_2,
\end{aligned}
\label{eq:di-b-coefficients}
\end{equation}
and
\begin{equation}
\begin{aligned}
b_{1,k}^{\alpha,m,\mathrm{ad}}
&=
\frac{\mathrm i^{\,k-1}k^m}{k!}
\left(
\delta^\alpha_1+(2m+1)\delta^\alpha_2
\right),
&
b_{2,k}^{\alpha,m,\mathrm{ad}}
&=
\frac{\mathrm i^{\,k}(k+1)^m}{k!}\delta^\alpha_2.
\end{aligned}
\label{eq:ad-b-coefficients}
\end{equation}
Here $k,m\geq0$.

The contractions in \eqref{eq:di-calibration-transformation} and \eqref{eq:ad-calibration-transformation} can be evaluated through the generating functions \eqref{eq:ad-calibration}--\eqref{eq:di-calibration}.
For example,
\[
\begin{aligned}
b_{1,k}^{\beta,m,\mathrm{di}}
\theta_{\beta,m}^{X^{\mathrm{di}}}(u)
&=
\frac1{k!}
\left[
\theta_1^{X^{\mathrm{di}}}(u;z)
+
2z\frac{\partial\theta_2^{X^{\mathrm{di}}}}{\partial z}(u;z)
\right]_{z=k},
\\
b_{1,k}^{\beta,m,\mathrm{ad}}
\theta_{\beta,m}^{X^{\mathrm{ad}}}(u)
&=
\frac{\mathrm i^{\,k-1}}{k!}
\left[
\theta_1^{X^{\mathrm{ad}}}(u;z)
+
\left(2z\frac{\partial}{\partial z}+1\right)
\theta_2^{X^{\mathrm{ad}}}(u;z)
\right]_{z=k}.
\end{aligned}
\]
At the integral values of $z$ appearing here, the expressions are
evaluated after cancellation of their apparent poles.

The same coefficients define the linear transformations of dilaton-shifted  descendant
variables:
\begin{equation}
\begin{aligned}
t_{\mathrm{di}}^{\alpha,k}(s) - \delta^\alpha_1\delta^k_{1}
&=
b_{\beta,\ell}^{\alpha,k,\mathrm{di}}[s^{\beta,\ell}
- \delta^\beta_1\delta^\ell_{1}],
\\
t_{\mathrm{ad}}^{\alpha,k}(s) -\delta^\alpha_1\delta^k_{1}
&=
b_{\beta,\ell}^{\alpha,k,\mathrm{ad}}[s^{\beta,\ell}
- \delta^\beta_1\delta^\ell_{1}].
\end{aligned}
\label{eq:descendant-transformations}
\end{equation} 
These are the descendant transformations entering the LYZZ
correspondences. 
LYZZ proved the genus-zero identities associated with these transformations \cite{liuEquivariantGromovWitten2022}. Section~\ref{sec:proof} establishes their positive-genus extensions, with the Hodge twist in the diagonal case.

\section{Landau-Ginzburg mirror, spectral curve and R-matrix}
\label{sec: B-theory} 
	 
 In this section, we recall the Landau-Ginzburg mirror of the resolved conifold and derive the
spectral curve by dimensional reduction. We subsequently determine the geometric
$R$-matrices and the Hodge twist in the diagonal case.

\subsection{Landau-Ginzburg potentials}

In the above we have recalled the A-theory side of the relevant theories.
Now we move on to the B-theory side.

Using the full equivariant mirror theorem,
the genus zero theory of a toric manifold
is encoded in the Jacobian ring of a Landau-Ginzburg mirror. 
For the resolved conifold $X={\rm Tot}(\cO_{\bP^1}(-1)\oplus\cO_{\bP^1}(-1))$, we introduce a torus action $\mathbb{T} = (\mathbb{C}^*)^3$ on $X$. Let $\fp_i$ be the $\bT$-fixed point and $\fl_j$ $(j=1,2,3)$ be the 1-dimensional $\bT$-orbit connected to $\fp_i$.
The weights at the two fixed points $\fp_+=0$ and $\fp_-=\infty$ are given by
\begin{equation}
(w_{+,1},w_{+,2},w_{+,3})=(\mu,a-\mu,b),
\qquad
(w_{-,1},w_{-,2},w_{-,3})=(-\mu,a,b+\mu),
\label{eq:full-torus-tangent-weights}
\end{equation}
where $w_{i,j}=e_\bT(T_{\fp_i}\fl_j)$, $\mu$ is the $\bT$-character on the base $\bP^1$ of $X$, and $a, b$ are the $\bT$-characters on $\cO(-1)$ fibers of $X$ respectively.
We will later take the limits $\mu \to 0$ and  $(a,b) \to (1,1)$ or $(1,-1)$
to obtain the B-theory mirror for $X^{\rm di}$ and $X^{\rm ad}$ respectively.

According to the recipe in \cite[Section 4.5]{fangRemodelingConjectureToric2019},
 the equivariant Landau-Ginzburg potential for the above three-torus action is 
\begin{equation}
W^{\mathbb {T}}(Y_1,Y_2,p)
=u^1+Y_1(1+p)
+Y_2\left(1+\frac{q}p\right)
-a\log Y_1-b\log Y_2-\mu\log p,
\label{eq:full-torus-LG-phase}
\end{equation}
with volume form
\[
\Omega=\frac{dY_1}{Y_1}\wedge\frac{dY_2}{Y_2}\wedge\frac{dp}{p}.
\]
Here $u^1$ is a constant, which is identified with the coordinate of the basis $\phi_1=1$ in $H^*(\bP^1;\bC)$; $q=e^{u^2}$ and 
$u^2$ is the coordinate of the hyperplane class $\phi_2=H$; $\mu, a, b$ are the weights of $\bT$.

Let $\text{Jac}(W^\mathbb{T})$ be the Jacobian ring of $W^\mathbb{T}$:
$$
    \text{Jac}(W^\mathbb{T}):=\frac{A_\mathbb{T}[u^1,q^{\pm 1},Y_1,Y_1^{-1},Y_2,Y_2^{-1},p,p^{-1}]}
    {\left\langle\frac{\partial W^\mathbb{T}}{\partial Y_1},\frac{\partial W^\mathbb{T}}{\partial Y_2},\frac{\partial W^\mathbb{T}}{\partial p}\right\rangle},
$$
where $A_\mathbb{T}$ is the minimal field extension of $\bC(a,b,\mu)$ that contains $\sqrt{e_\bT(T_{\fp_i}X)}$. 
For $i\in \{+,-\}$, let ${\bf p}_i$ be the critical point of $W^\bT$, $\text{Jac}(W^\mathbb{T})$ is a Frobenius algebra under multiplication and residue pairing:
$$
  (f,g) := \sum_{i=\pm} \frac{f({\bf p}_i)g({\bf p}_i)}{\det(\text{Hess}(W^\bT))({\bf p}_i)},
$$
where $\text{Hess}(W^\bT)$ is the Hessian of $W^\bT$ with respect to $(\log Y_1,\log Y_2, \log p)$.
Let $\text{QH}^*_\mathbb{T}(X)$ be the $\bT$-equivariant quantum cohomology ring of $X$.
By the argument in \cite{tsengSeidelRepresentationsQuantum2012}, we have the isomorphism of Frobenius algebras:
$$
    \text{QH}^*_\mathbb{T}(X)\rightarrow \text{Jac}(W^\mathbb{T}): 1\rightarrow [1],\quad H\rightarrow \left[\frac{\partial W^\mathbb{T}}{\partial u^2}\right] = \left[\frac{qY_2}{p}\right] .
$$
In the Jacobian ring $\text{Jac}(W^{\bT})$, 
the flat metric in the ordered basis $(1,H)$ is given by 
\[
\eta^{\bT}=
\begin{pmatrix}
\dfrac{a+b}{ab(a-\mu)(b+\mu)}
&
\dfrac1{a(b+\mu)}
\\[6pt]
\dfrac1{a(b+\mu)}
&
-\dfrac{\mu}{a(b+\mu)}
\end{pmatrix},
\]
and the quantum product is given by
$$
    H\star(H+\mu) = q(H-a+\mu)\star(H-b).
$$
For $i\in \{+,-\}$, let 
$V_i(u)$ be a Laurent polynomial in $Y_1,Y_2,p$ such that 
$V_i(u)({\bf p}_{j})=\delta_{i,j}$. The class of $V_i(u)$ in the Jacobian ring, denoted by $[V_i(u)]$,
is a canonical basis of the semisimple Frobenius algebra $\Jac(W^{\bT})$. 
The scalar factor $\Delta_i$ in \eqref{eq:canonical-idempotents} is identified with $\det(\text{Hess}(W^\bT))({\bf p}_i)$.
Let $\Phi$ be the transition matrix between the flat basis and normalized canonical basis as \eqref{eq:normalized-idempotent-matrix}.

We can write $[V_i(u)]$ in the following decomposition
$$
    [V_i(u)] = B_i(u)\left[\frac{\partial W^\mathbb{T}}{\partial u^2}\right] + C_i(u),
$$
where $$B_i(u) = \frac{\Phi^2_{\ i}}{\sqrt{\Delta_i}}, \quad C_i(u) = \frac{\Phi^1_{\ i}}{\sqrt{\Delta_i}}.$$
Define the standard form of $[V_i(u)]$ to be
$$
    \overline{V_i(u)} = B_i(u)\frac{\partial W^\mathbb{T}}{\partial u^2} + C_i(u).
$$
The standard form $\overline{V_i(u)}$ will be used to define the $R$-matrix given by the period integrals.

Let $\Gamma$ be a 3-cycle on $H_3((\mathbb{C}^*)^3,\{\mathfrak{R}(\frac{W^\mathbb{T}}{z}\ll0)\};\mathbb{Z})$ and consider the period integral
$$
    f(u,z) = \int_\Gamma e^{\frac{W^{\mathbb{T}}}{z}}\Omega,
$$
which satisfies
\begin{equation}
\left[
(z\partial_{u^2})(z\partial_{u^2}+\mu)
-e^{u^2}(z\partial_{u^2}-a+\mu)(z\partial_{u^2}-b)
\right]f=0.
\label{eq:full-torus-scalar-equation}
\end{equation}
The full equivariant mirror theorem identifies this equation with
the quantum differential equation on the A-theory side. 

\subsection{Spectral curve from the Landau-Ginzburg potential}

We now reduce the 3-dimensional period integrals to the 1-dimensional oscillating integrals
to see how the spectral curves appear on the B-theory side.
We integrate out the two fiber variables. The substitutions
\[
A_1=(1+p)Y_1,
\qquad
A_2=\left(1+\frac{e^{u^2}}p\right)Y_2
\]
preserve the logarithmic volume form and give the exact separation
\begin{equation}
\begin{aligned}
W^\mathbb{T}
={}&(A_1-a\log A_1)+(A_2-b\log A_2)
\\
&+u^1+a\log(1+p)
+b\log\!\left(1+\frac{e^{u^2}}p\right)-\mu\log p.
\end{aligned}
\label{eq:full-torus-phase-separation}
\end{equation}
Let $q\neq 0,1$ and $0<\mu <\min\{|a|, |b|\}$, 
\begin{align}\label{eq:a-b-spectral-curve}
x(p)
=
u^1+a\log(1+p)
+b\log\left(1+\frac qp\right)-\mu\log p,
\qquad
y^\circ(p)=\frac{\log p}{\sqrt{ab}}.
\end{align}
Choose $q$ such that all critical points of $x(p)$ are simple. Let $p_i$ be a critical point of $x$, $\gamma_i$ be the Lefschetz thimble passing through $p_i$ such that
$x(p)-x(p_i)\in\bR_{\geq0}$.

We consider the 3-cycle $\Gamma_{i} \in H_3((\mathbb{C}^*)^3,\{\mathfrak{R}(\frac{W^\mathbb{T}}{z}\ll0)\};\mathbb{Z})$ as
\[
\Gamma_{i}
=
\left\{
(Y_1,Y_2,p):
\ p\in\gamma_{i},\
(1+p)Y_1\in\mathbb R_+,\
\left(1+\frac qp\right)Y_2\in\mathbb R_+
\right\}.
\]
This cycle is parameterized by
\[
(A_1,A_2,p)\longmapsto
\left(
\frac{A_1}{1+p},
\frac{A_2}{1+q/p},
p
\right),
\qquad
(A_1,A_2,p)\in\mathbb R_+^2\times\gamma_{i}.
\]
We have the following dimensional reduction of the period integrals:
\begin{lemma}\label{lemma:dimensional-reduction}
    Let 
    $$
    \widetilde{I}_i(u,z)=\int_{\Gamma_{i}}e^{W^\mathbb{T}/z}\Omega,
    \qquad
    I_i(u,z)=\int_{\gamma_{i}}e^{x/z}y^\circ\,dx.
    $$
    Assume $\Re z<0$ and $a, b >0$, then 
    $$
        \widetilde{I}_i(u,z) = \sqrt{ab}(-z)^{-(a+b)/z-1}\Gamma(-\frac{a}{z})\Gamma(-\frac{b}{z})I_i(u,z).
    $$
    Formally extend $\widetilde{I}_i(u,z)$ to general $a,b$. When $z\rightarrow 0$, we have the asymptotic expansion
    $$
\frac{e^{-c_{a,b}/z}}{2\pi}\,
\widetilde{I}_i(u,z)
\sim
\exp\!\left[-
\sum_{m\ge1}
\frac{B_{2m}}
{2m(2m-1)}
\bigl(a^{1-2m}+b^{1-2m}\bigr)z^{2m-1}
\right]I_i(u,z).
    $$
    where $c_{a,b}=a+b-a\log a-b\log b.$
\end{lemma}

The one-dimensional oscillating integrals have the following phase
functions:
\begin{subequations}\label{eq:oscillating-integral-phases}
\renewcommand{\theequation}{\theparentequation-\roman{equation}}
\begin{align}
x^{\mathbb P^1}(v;p)
&=v^1+p+\frac{e^{v^2}}p,
&
y^{\mathbb P^1,\circ}(p)
&=\log p,
\label{eq:oscillating-integral-phases-p1}
\\
x^{X^{\mathrm{ad}}}(u;p)
&=u^1+\log(1+p)-\log\!\left(1+\frac{e^{u^2}}p\right),
&
y^{X^{\mathrm{ad}},\circ}(p)
&=-\mathrm i\log p,
\label{eq:oscillating-integral-phases-ad}
\\
x^{X^{\mathrm{di}}}(u;p)
&=u^1+\log(1+p)+\log\!\left(1+\frac{e^{u^2}}p\right),
&
y^{X^{\mathrm{di}},\circ}(p)
&=\log p.
\label{eq:oscillating-integral-phases-di}
\end{align}
\end{subequations}
The Equation~\eqref{eq:oscillating-integral-phases-p1} is the spectral curve of $\bP^1$ \cite{norburyGromovWittenInvariants2014}. 
The \eqref{eq:oscillating-integral-phases-ad} and \eqref{eq:oscillating-integral-phases-di} are the limits of \eqref{eq:a-b-spectral-curve} when $\mu\rightarrow 0$ and $a=-b=1$ (\emph{resp.} $a=b=1$).
In each case, the one-dimensional oscillating integral is of the form $\int e^{x/z}\,dy^\circ$.
The critical points of the three phase functions $x$ are
\begin{equation}
\begin{aligned}
p_\pm^{\mathbb P^1}
&=\pm e^{v^2/2},
\\
p_\pm^{X^{\mathrm{ad}}}
&=-e^{u^2}\pm\mathrm i\,e^{u^2/2}\sqrt{1-e^{u^2}},
\\
p_\pm^{X^{\mathrm{di}}}
&=\pm e^{u^2/2}.
\end{aligned}
\label{eq:phase-critical-points}
\end{equation}
Their critical values are the canonical coordinates in
Section~\ref{subsec:Frobenius} and the critical points are nondegenerate considered there.

\subsection{The \texorpdfstring{$R$}{R}-matrices}
\label{subsec:R-matrices}

At a semisimple point, the quantum differential equations \eqref{eq:S-equations} admit an
asymptotic fundamental solution of the form $\Phi R(z)e^{U/z}$ for $z\to 0$ with
\begin{equation}
U=\operatorname{diag}(U_+,U_-),
\qquad
R(z)=I+\sum_{k\geq1}R_kz^k,
\qquad
R(-z)^T R(z)=I.
\label{eq:R-fundamental-solution}
\end{equation}
Here $R(z)$ is an endomorphism-valued formal power series in the
normalized canonical basis, and $\Phi$ is the matrix defined in
\eqref{eq:normalized-idempotent-matrix}. The existence of this solution
and its ambiguity are described in
\cite{dubrovinGeometry2DTopological1996,dubrovinPainleveTranscendentsTwoDimensional1999,giventalGromovWittenInvariantsQuantization2001}.
Substitution into the quantum differential equations gives
\begin{equation}
[dU,R_{k+1}]
=dR_k+\Phi^{-1}d\Phi\,R_k,
\qquad k\geq0,
\qquad R_0=I.
\label{eq:R-recursion}
\end{equation}
The symplectic condition in \eqref{eq:R-fundamental-solution} leaves the freedom of right multiplication by
\begin{equation}
\exp\!\left(\sum_{m\geq1}A_{2m-1}z^{2m-1}\right),
\label{eq:scalar-factor}
\end{equation}
where the matrices $A_{2m-1}$ are constant and diagonal. For
$\mathbb P^1$, homogeneity fixes this ambiguity. For the
resolved conifold, the ambiguity is fixed by its equivariant Gromov--Witten invariants in the large radius limit.
We write $R^{\mathbb P^1}$, $R^{X^{\mathrm{ad}}}$, and
$R^{X^{\mathrm{di}}}$ for the resulting A-model $R$-matrices.
 We also write $R^{X,\bT}$ for the A-model $R$-matrix of full equivariant resolved conifold in \eqref{eq:full-torus-tangent-weights}.

We will compute these matrices in the B-model using stationary phase asymptotics of oscillating integrals.
This is Givental's construction of asymptotic
fundamental solutions from a Landau--Ginzburg mirror; see
\cite[Section 10]{giventalGromovWittenInvariantsQuantization2001}. 
For our purpose we need to apply the same strategy above to reduce to an integral of 
a three-form on a three-cycle to an integral of a one-form on a one-cycle,
which is a generalized Laplace transform one sees very often in the literature of topological recursions
(see e.g. \cite{eynardInvariantsSpectralCurves2014}).
In this subsection we will identify the A-side $R$-matrix with the B-side $R$-matrix.
In the next subsection we will reduce the computations of B-side $R$-matrix to integrals over one-cycles.
  
Let ${\bf U}_i := W^\bT({\bf p}_i)$, we define 
\begin{equation}\label{eq:3-cycle-integrals}
      \check{R}_{\ j}^{i}(u,z)\sim \frac{\sqrt{\Delta_i}}{(-2\pi z)^{3/2}}e^{-{\bf U}_j/z}\int_{\Ga_j}e^{\frac{W^\bT}{z}}\overline{V_i(u)}\Omega.
\end{equation}
We refer to $\check{R}(u,z)$ as the $R$-matrix defined by the three-dimensional period integrals.
Let's consider the column vector
\begin{equation}
{\bf s}_i
:= \frac{1}{(-2\pi z)^{3/2}}(\eta^{\bT})^{-1}\begin{pmatrix}\widetilde{I}_i\\z\partial_{u^2}\widetilde{I}_i\end{pmatrix}.
\label{eq:R-matrix-vector}
\end{equation}
The vector ${\bf s}_i$ admits the asymptotic expansion
\begin{equation}
\begin{aligned}
(\Phi^{-1}\mathbf s_j)^i
&=
\frac{\Phi^1{}_i\widetilde I_j+
      \Phi^2{}_i(z\partial_{u^2}\widetilde I_j)}{(-2\pi z)^{3/2}}
\\
&=
\frac{\sqrt{\Delta_i}}{(-2\pi z)^{3/2}}
\int_{\Gamma_j}
e^{W^{\mathbb T}/z}\overline{V_i(u)}\,\Omega
\sim
\check R^i{}_j(u,z)e^{\mathbf U_j/z}.
\end{aligned}    
\end{equation}
The Equation~\eqref{eq:full-torus-scalar-equation} shows ${\bf s}_i$ satisfies the same recursion relation as the quantum differential equation of the $\bT$-equivariant
resolved conifold. The $R$-matrices $\check{R}(u,z)$ and $R^{X,\bT}(u,z)$ differ up to a constant matrix as \eqref{eq:scalar-factor}. The ambiguity is fixed by the initial value in the large radius limit.
\begin{lemma}\label{lemma:period-integral-initial-value}
In the large radius limit $u^1,q\rightarrow0$, the $R$-matrix given by the period integrals in \eqref{eq:3-cycle-integrals} is 
\begin{align}\label{eq:period-integral-initial-value}
\left.\check R^{i}_{\ j}(u,z)\right|_{u^1=q=0}
=
\delta^{i}_{\ j}
\exp\left[-
\sum_{m\ge1}
\frac{B_{2m}}{2m(2m-1)}
\sum_{k=1}^3w_{j,k}^{1-2m}z^{2m-1}
\right].
\end{align}
\end{lemma}
\begin{proof}
We verify this initial value separately for the two branches.
After integrating out the fiber variables, Lemma
\ref{lemma:dimensional-reduction} gives
\[
\widetilde I_i(u,z)
=
(-z)^{-(a+b)/z}
\Gamma(-a/z)\Gamma(-b/z)
\int_{\gamma_i}e^{x(p)/z}\frac{dp}{p}.
\]

For the branch corresponding to $\fp_+$, the critical point
approaches $\mu/(a-\mu)$ as $q\to0$. The limiting phase is
\[
u^1+a\log(1+p)-\mu\log p.
\]
After extracting $e^{u^1/z}$, the remaining integral is
\[
\int_0^\infty p^{-\mu/z-1}(1+p)^{a/z}\,dp
=
\frac{\Gamma(-\mu/z)\Gamma(-(a-\mu)/z)}
     {\Gamma(-a/z)}.
\]
Multiplying by the two fiber Gamma factors cancels
$\Gamma(-a/z)$ and leaves
\[
\Gamma(-\mu/z)\Gamma(-(a-\mu)/z)\Gamma(-b/z),
\]
which corresponds to the three weights $(\mu,a-\mu,b)$.

For the branch corresponding to $\fp_-$, the critical point
satisfies $p\sim-(b+\mu)q/\mu$. To retain this branch in the
limit, set $t=q/p$. In this coordinate,
\[
x
=
u^1-\mu\log q
+a\log(1+q/t)+b\log(1+t)+\mu\log t.
\]
Thus, after extracting $e^{u^1/z}q^{-\mu/z}$ and taking
$q\to0$, the remaining integral is
\[
\int_0^\infty t^{\mu/z-1}(1+t)^{b/z}\,dt
=
\frac{\Gamma(\mu/z)\Gamma(-(b+\mu)/z)}
     {\Gamma(-b/z)}.
\]
Here the orientation is chosen to match the Gaussian
normalization. Multiplication by the fiber Gamma factors
now cancels $\Gamma(-b/z)$ and leaves
\[
\Gamma(\mu/z)\Gamma(-a/z)\Gamma(-(b+\mu)/z),
\]
corresponding to the weights $(-\mu,a,b+\mu)$.

For $a,b>0$ and $z<0$, the first and second beta integrals
converge in the regions $0<\mu<a$ and $-b<\mu<0$,
respectively. We first perform the calculations in these
regions and then continue the coefficients of the normalized
stationary-phase expansions in the equivariant parameters.
In each branch, Stirling's formula, together with
$\Delta_i\to\prod_{k=1}^3w_{i,k}$, shows that the normalization
in \eqref{eq:3-cycle-integrals} removes the critical-value
exponential and the Gaussian leading factor. Thus, for $i\in\{+,-\}$, we have
$$
    \lim_{u^1,q\rightarrow0}\frac{\sqrt{\Delta_i}}{(-2\pi z)^{3/2}}e^{-{\bf U}_i/z}\widetilde{I}_i \sim \exp\left[-
\sum_{m\ge1}
\frac{B_{2m}}{2m(2m-1)}
\sum_{k=1}^3w_{i,k}^{1-2m}z^{2m-1}
\right].
$$

At $q=0$, the two idempotents are
$V_+=1+H/\mu$ and $V_-=-H/\mu$.
Their standard insertions act on the periods as
$1+(z/\mu)\partial_{u^2}$ and
$-(z/\mu)\partial_{u^2}$, respectively.
The leading $q$-dependence of the two periods is
$q^0$ and $q^{-\mu/z}$.
Thus each insertion preserves the corresponding leading
period and annihilates the other one, so the off-diagonal
entries vanish in the large-radius limit.
The diagonal entries are given by the normalized Gamma
products computed above, yielding the stated
Bernoulli-series initial value.
\end{proof}

The results in \cite{giventalGromovWittenInvariantsQuantization2001,giventalSemisimpleFrobeniusStructures2001} imply that
the A-model  $R$-matrix of the equivariant Gromov-Witten theory of $X$ with torus action $\mathbb{T}$ is uniquely determined by the initial value:
\begin{align}\label{eq:A-model-initial-value}
R^{X,\mathbb{T}}(u,z)\big|_{u^1=e^{u^2}=0} = \text{diag}_{i=\pm}\exp\left[-\sum_{m\ge1}\frac{B_{2m}}{2m(2m-1)}
\sum_{k=1}^3 w_{i,k}^{1-2m}z^{2m-1}\right].
\end{align} 
The initial value in the large radius limit (see \eqref{eq:A-model-initial-value} and \eqref{eq:period-integral-initial-value}) fixes the ambiguity of right multiplication in \eqref{eq:scalar-factor}.
We get the identification of A-model $R$-matrix and $R$-matrix defined by period integrals
\begin{align}\label{eq:A-period-R-identification}
    \check{R}(u,z) = R^{X,\bT}(u,z).
\end{align}

\medskip
\subsection{A-model $R$-matrices vs B-model $R$-matrices}\label{subsec:R-comparison}
We fix the following convention for local Laplace transforms below; see also \cite{eynardInvariantsSpectralCurves2014}.
Let $p_i$ be a nondegenerate critical point of $x$, set $U_i=x(p_i)$,
and choose a local Airy coordinate such that
\[
x=U_i+\frac{\zeta_i^2}{2}.
\]
For the three functions $y^\circ$ in \eqref{eq:oscillating-integral-phases}, we choose the square roots consistently so that $$(dy^\circ/d\zeta_i)(0)=\Delta_i^{-1/2}.$$

\begin{definition}
For a one-form $\alpha$ holomorphic near $p_i$, its normalized local
Laplace transform is the formal series
\begin{equation}
\mathcal G_{i,z}(\alpha)
:=\sum_{m\geq0}(-1)^m(2m-1)!!
\left[\zeta_i^{2m}\right]\frac{\alpha}{d\zeta_i}\,z^m.
\label{eq:normalized-Gaussian-expansion}
\end{equation}
We also use the column notation
\begin{equation}
\bm{\mathcal G}_z(\alpha)
:=\begin{pmatrix}\mathcal G_{+,z}(\alpha)\\\mathcal G_{-,z}(\alpha)\end{pmatrix}.
\label{eq:Laplace-column}
\end{equation}
The notation $\bm{\mathcal G}_{-z}$ means substitution of $-z$ in this same formal transform.
\end{definition}

When a convergent integral representation is available, this is the
stationary-phase expansion
\begin{equation}
\frac{e^{-U_i/z}}{\sqrt{-2\pi z}}
\int_{\gamma_i(z)}e^{x/z}\alpha
\sim\mathcal G_{i,z}(\alpha).
\label{eq:Laplace-integral}
\end{equation}
Here $\gamma_i(z)$ is an admissible steepest-descent thimble through
$p_i$. Its orientation and the square-root branch are fixed by
$\int e^{\zeta_i^2/(2z)}d\zeta_i=\sqrt{-2\pi z}$.
The substitution $z\mapsto-z$ changes the decay sector and, when
needed, the thimble.
We use \eqref{eq:normalized-Gaussian-expansion} as the definition, so no global integral representation is required.

For a residue-free one-form with at most a double pole, we extend
the definition by formal regularization:
\begin{equation}
\mathcal G_{i,z}\!\left(c_{-2}\frac{d\zeta_i}{\zeta_i^2}
+\alpha_{\mathrm{hol}}\right)
:=\frac{c_{-2}}z+\mathcal G_{i,z}(\alpha_{\mathrm{hol}}).
\label{eq:Laplace-regularization}
\end{equation}

\begin{lemma}
\label{lem:Laplace-properties}
The normalized transforms have the following properties.
\begin{enumerate}[(i)]
\item The zero-order moment and odd-order moments satisfy
\begin{equation}
\mathcal G_{i,z}(d\zeta_i)=1,\qquad
\mathcal G_{i,z}(\zeta_i^{2m+1}d\zeta_i)=0\quad(m\geq0).
\label{eq:Laplace-moments}
\end{equation}
For any holomorphic $y=\sum_{m\geq0}y_m^i\zeta_i^m$,
\begin{equation}
z\mathcal G_{i,-z}(dy)
=\sum_{k\geq0}(2k+1)!!\,y_{2k+1}^i z^{k+1}.
\label{eq:Laplace-odd-series}
\end{equation}
In particular, $\mathcal G_{i,z}(dH(x))=0$ for every holomorphic
germ $H$ at $U_i$.
\item Formal integration by parts gives
\begin{equation}
\mathcal G_{i,z}\!\left(dF+\frac{F\,dx}{z}\right)=0,
\label{eq:Gaussian-integration-by-parts}
\end{equation}
and therefore
\begin{equation}
\bm{\mathcal G}_{-z}\!\left(dF-\frac{F\,dx}{z}\right)=0.
\label{eq:negative-Gaussian-integration-by-parts}
\end{equation}
The identity holds when $F$ has at most
a simple pole at the relevant critical point.
\item Suppose that $x$ depends on primary parameters $t^\beta$,
and that $\alpha$ is holomorphic near the moving critical point
and independent of these parameters $t^\beta$ in the local Airy coordinate. Then,
\begin{equation}
e^{-U_i/z}z\partial_\beta
\left(e^{U_i/z}\mathcal G_{i,z}(\alpha)\right)
=\mathcal G_{i,z}\bigl((\partial_\beta x)\alpha\bigr).
\label{eq:Laplace-primary-derivative}
\end{equation}
Here $\partial_\beta x=\partial_\beta U_i$.
\item The transform is unchanged by adding a constant to $x$,
and is covariant under a change of surface coordinate. 
\end{enumerate}
\end{lemma}
\begin{proof}
All assertions follow from the definition in the induced local Airy coordinates.
\end{proof}

For the $y^{\circ}$ functions, define the periods by
\begin{equation}
f_i(z):=e^{U_i/z}\mathcal G_{i,z}(dy^\circ).
\label{eq:stationary-phase-periods}
\end{equation}
The unnormalized periods $I_i$ in
Lemma~\ref{lemma:dimensional-reduction} use $y^\circ dx$ instead of
$dy^\circ$. By formal integration by parts (or for admissible cycles
with vanishing boundary terms),
\[
\int e^{x/z}y^\circ dx=-z\int e^{x/z}dy^\circ,
\qquad
\frac{I_i(z)}{\sqrt{-2\pi z}}\sim-z f_i(z).
\]

Let $\Xi_0^i(p)$ be the 1-form 
$$
  \left.\Xi_0^i(p)=\frac{dp}{d\zeta_i}\right|_{\zeta_i=0}\frac{dp}{(p-p_i)^2}.
$$
We consider the one-dimensional oscillating integrals
\begin{align}\label{eq:1-dim-oscillating-integral}
    R^i_{\ j}(z)=z\,\mathcal G_{j,z}(\Xi^i_0),
\end{align}
It is related to the $R$-matrix associated with the Eynard-Orantin topological recursion (see Section \ref{sec:TR}). Lemma \ref{lemma:dimensional-reduction}
reduces the 3-dimensional period integrals in Equation \eqref{eq:3-cycle-integrals} to $\mathcal G_{j,z}(\Xi^i_0)$. Furthermore, we relate the A-model $R$-matrix to one-dimensional oscillating
integrals by \eqref{eq:A-period-R-identification}; see Proposition \ref{prop:geometric-R-matrices}.
\begin{proposition}
\label{prop:geometric-R-matrices}
The stationary phase asymptotics of the oscillating integrals in \eqref{eq:1-dim-oscillating-integral}
associated with the first two pairs 
\eqref{eq:oscillating-integral-phases-p1} and \eqref{eq:oscillating-integral-phases-ad}, determine the A-model matrices $R^{\mathbb P^1}(z)$ and
$R^{X^{\mathrm{ad}}}(z)$, respectively. If
$\widetilde R^{X^{\mathrm{di}}}(z)$ denotes the reduced $R$-matrix obtained from
the third pair \eqref{eq:oscillating-integral-phases-di} by one-dimensional oscillating integrals in \eqref{eq:1-dim-oscillating-integral}, then
\begin{equation}
R^{X^{\mathrm{di}}}(z)
=
\exp\!\left[
-2\sum_{m\geq1}
\frac{B_{2m}}{2m(2m-1)}z^{2m-1}
\right]
\widetilde R^{X^{\mathrm{di}}}(z),
\label{eq:di-geometric-R-correction}
\end{equation}
where $B_{2m}$ are the Bernoulli numbers.
\end{proposition}

\begin{proof}
The identification of the A-model $R$-matrix $R^{\bP^1}(z)$ and the oscillating integrals associated to \eqref{eq:oscillating-integral-phases-p1}
is a classical result in \cite{dunin-barkowskiIdentificationGiventalFormula2014}; see also for example \cite[Proposition 3.10]{fangEynardOrantinRecursionEquivariant2017}.

The A-model $R$-matrix of resolved conifold $X$ is identified with the 3-dimensional period integral $\check{R}(u,z)$ in \eqref{eq:A-period-R-identification}.
Lemma \ref{lemma:dimensional-reduction} reduces 3-dimensional period integrals to oscillating integrals on Lefschetz thimble.
One can verify that
$$
    -\frac{\Xi^i_0}{\sqrt{\Delta_i}} = d_p\left[\Big(B_i(u)\partial_{u^2}x\big|_p + C_i(u)\Big)\frac{dy^\circ}{dx}\right].
$$
Then
\begin{align*}
    &\quad\int_{\Ga_i}e^{\frac{W^\bT}{z}}\overline{V_{j}(u)}\Omega = \Big(zB_j(u)\frac{\partial}{\partial u^2}+C_j(u)\Big)\int_{\Ga_i}e^{\frac{W^\bT}{z}}\Omega
    \\
    &= \sqrt{ab}(-z)^{-(a+b)/z}\Gamma(-\frac{a}{z})\Gamma(-\frac{b}{z})\Big(zB_j(u)\frac{\partial}{\partial u^2}+C_j(u)\Big)\int_{\gamma_i}e^{\frac{x}{z}}dy^\circ 
    \\
    &= \sqrt{ab}(-z)^{-(a+b)/z}\Gamma(-\frac{a}{z})\Gamma(-\frac{b}{z})\int_{\gamma_i}e^{\frac{x}{z}}\left[\Big(B_j(u)\partial_{u^2}x\big|_p + C_j(u)\Big)\frac{dy^\circ}{dx}\right]dx
    \\
    &= -\sqrt{ab}(-z)^{-(a+b)/z+1}\Gamma(-\frac{a}{z})\Gamma(-\frac{b}{z})\int_{\gamma_i}e^{\frac{x}{z}} \frac{\Xi^j_0}{\sqrt{\Delta_j}}.
\end{align*}

Then we have:
$$
    \check{R}^i_{\ j}(z) = \exp\!\left[-
\sum_{m\ge1}
\frac{B_{2m}}
{2m(2m-1)}
\bigl(a^{1-2m}+b^{1-2m}\bigr)z^{2m-1}
\right]R^i_{\ j}(z).
$$
For $(a,b)=(1,-1)$, the two terms in
$a^{1-2m}+b^{1-2m}$ cancel. The remaining differential is
$-\mathrm i\,dp/p$, and at $\mu=0$ the remaining phase is
$x^{X^{\mathrm{ad}}}$. Hence its stationary phase expansion gives
$R^{X^{\mathrm{ad}}}$ directly. For $(a,b)=(1,1)$, the differential
is $dp/p$ and the two fiber contributions in
Lemma~\ref{lemma:dimensional-reduction} coincide. At $\mu=0$ the phase is
$x^{X^{\mathrm{di}}}$, giving
\eqref{eq:di-geometric-R-correction}.
\end{proof}
\begin{remark}
  The specialization of torus weight to $\mu=0$ is taken at a locus of $\{u^1,u^2\}$, referred to as the semisimple locus, such that the critical points and Hessians of $x(p)$ in \eqref{eq:a-b-spectral-curve} remain
  nondegenerate there. On such a locus, the stationary phase expansion of the oscillating integrals associated to \eqref{eq:oscillating-integral-phases-ad} and \eqref{eq:oscillating-integral-phases-di} has regular limits.
  We begin with $\mu\neq0$ in order to identify the A-model $R$-matrix and period integrals \eqref{eq:A-period-R-identification} in the large radius limits.
\end{remark}

  \begin{remark}
    The recipe of determining $R$-matrix by the Landau-Ginzburg mirror and dimensional reduction is a well-established method in 
    \cite{fangRemodelingConjectureToric2019,liuAllgenusOpenclosedMirror2020} for semi-projective toric Calabi-Yau 3-orbifolds, 
under the  assumptions that these toric orbifolds are equipped with a Calabi-Yau torus 
    action with finitely many 0-dimensional and 1-dimensional orbits. 
    These assumptions are also needed in their work to establish the topological recursions 
    as predicted by the BKMP Conjecture.
    In our cases, the diagonal action does not preserve the Calabi-Yau form on the resolved conifold;
    the fixed locus of $X^{\rm di}, X^{\rm ad}$ is 1-dimensional. So we rederive the spectral curve and B-model $R$-matrix following a similar strategy. We also present a separate proof of the topological recursions in the next section.
  \end{remark}

The logarithm of the scalar factor in
\eqref{eq:di-geometric-R-correction} begins with
\[
-\frac z6+\frac{z^3}{180}-\frac{z^5}{630}+\cdots.
\]
In particular,
\begin{equation}
R_1^{X^{\mathrm{di}}}
=\widetilde R_1^{X^{\mathrm{di}}}-\frac16 I.
\label{eq:di-first-R-correction}
\end{equation}

For later use, we record the first coefficients of the three geometrically normalized \(R\)-matrices,
\[
R(z)=I+R_1z+O(z^2),
\]
in the normalized canonical basis ordered by \(+,-\).
Expanding the stationary-phase formulas to first order gives
\begin{equation}
R_1^{\mathbb P^1}
=
\begin{pmatrix}
\displaystyle-\frac{1}{16e^{v^2/2}}
&
\displaystyle
-\frac{1}{
4\sqrt{\Delta_+^{\mathbb P^1}}
 \sqrt{\Delta_-^{\mathbb P^1}}}
\\[8pt]
\displaystyle
-\frac{1}{
4\sqrt{\Delta_+^{\mathbb P^1}}
 \sqrt{\Delta_-^{\mathbb P^1}}}
&
\displaystyle\frac{1}{16e^{v^2/2}}
\end{pmatrix},
\label{eq:P1-first-R}
\end{equation}
and
\begin{equation}
R_1^{X^{\mathrm{ad}}}
=
\begin{pmatrix}
\displaystyle
\frac{\mathrm i(1-2e^{u^2})}
{16e^{u^2/2}\sqrt{1-e^{u^2}}}
&
\displaystyle
\frac{1}{
4(1-e^{u^2})
\sqrt{\Delta_+^{X^{\mathrm{ad}}}}
\sqrt{\Delta_-^{X^{\mathrm{ad}}}}}
\\[10pt]
\displaystyle
\frac{1}{
4(1-e^{u^2})
\sqrt{\Delta_+^{X^{\mathrm{ad}}}}
\sqrt{\Delta_-^{X^{\mathrm{ad}}}}}
&
\displaystyle
-\frac{\mathrm i(1-2e^{u^2})}
{16e^{u^2/2}\sqrt{1-e^{u^2}}}
\end{pmatrix}.
\label{eq:ad-first-R}
\end{equation}
For the diagonal theory, including the scalar correction in
\eqref{eq:di-geometric-R-correction}, we obtain
\begin{equation}
R_1^{X^{\mathrm{di}}}
=
\begin{pmatrix}
\displaystyle-\frac{1+e^{u^2}}{16e^{u^2/2}}
&
\displaystyle
-\frac{1}{
4\sqrt{\Delta_+^{X^{\mathrm{di}}}}
 \sqrt{\Delta_-^{X^{\mathrm{di}}}}}
\\[8pt]
\displaystyle
-\frac{1}{
4\sqrt{\Delta_+^{X^{\mathrm{di}}}}
 \sqrt{\Delta_-^{X^{\mathrm{di}}}}}
&
\displaystyle\frac{1+e^{u^2}}{16e^{u^2/2}}
\end{pmatrix}
-\frac16 I.
\label{eq:di-first-R}
\end{equation}
Here the square roots of \(\Delta_\pm\) are the same as those
used in \eqref{eq:Phi-matrices}.
In particular, all three matrices satisfy \(R_1^T=R_1\),
as required by \eqref{eq:R-fundamental-solution}.

\subsection{The Hodge twist in the diagonal theory}
\label{subsec:Hodge-correction}
In Proposition \ref{prop:geometric-R-matrices}, 
we have identified the $R$-matrices $R^{\bP^1}$ and $R^{X^{\rm ad}}$ with the $R$-matrices  
determined by the spectral curves \eqref{eq:oscillating-integral-phases-p1}
and \eqref{eq:oscillating-integral-phases-ad} respectively,
but for the case of $R^{X^{\rm di}}$,
there is an extra factor on the right-hand side of \eqref{eq:di-geometric-R-correction}. 
This means the cohomological field theory determined by the spectral curve  \eqref{eq:oscillating-integral-phases-di}
is not exactly the equivariant Gromov-Witten theory of the resolved conifold
with diagonal circle action.
Instead,
it corresponds to the equivariant Gromov-Witten invariants twisted by suitable Hodge bundles.

Let $\mathbb E_g=\pi_*\omega_\pi$ be the rank-$g$ Hodge bundle on
$\overline{\mathcal M}_{g,n}(\mathbb P^1,d)$, and write
\begin{equation}
\lambda_j=c_j(\mathbb E_g),
\qquad
\Lambda_g^\vee(w)
=\sum_{j=0}^g(-1)^j\lambda_jw^{g-j},
\qquad \lambda_0=1.
\label{eq:Hodge-polynomial}
\end{equation}
The class $\Lambda_g^\vee(1)$ has degree-zero term $1$, so its
inverse is well defined in cohomology.

We would like to consider the twisted descendant Gromov-Witten invariants of $X^{\mathrm{di}}$ by the class
$R^\bullet\pi_*\mathcal{O}=\mathcal{O}-\mathbb E_g^\vee$,
$$
    e_T\big(R^\bullet\pi_*(\mathcal{O}_1\oplus\mathcal{O}_1)\big)=\Lambda_g^\vee(1)^{-2}.
$$

Define the twisted diagonal descendant invariants by
\begin{equation}
\begin{aligned}
&\left\langle
\prod_{i=1}^n\tau_{k_i}(\phi_{\alpha_i})
\right\rangle^{X^{\mathrm{di}},\mathrm{tw}}_{g,d}
\\
&=
\Biggl[
\int_{[\overline{\mathcal M}_{g,n}(\mathbb P^1,d)]^{\mathrm{vir}}}
\Lambda_g^\vee(1)^{-2}
e_T\!\left(
-R^\bullet\pi_*f^*\mathcal O_{\mathbb P^1}(-1)_1
\right)^2
\prod_{i=1}^n\psi_i^{k_i}
\operatorname{ev}_i^*(\phi_{\alpha_i})
\Biggr]_{\lambda=1}.
\end{aligned}
\label{eq:corrected-diagonal-correlators}
\end{equation}
Here the superscript $\mathrm{tw}$ records the Hodge insertion.
The corresponding generating functions are
\begin{equation}
\begin{aligned}
F_g^{X^{\mathrm{di}},\mathrm{tw}}(t)
&=
\sum_{d,n\geq0}\frac{1}{n!}
t^{\alpha_1,k_1}\cdots t^{\alpha_n,k_n}
\left\langle
\prod_{i=1}^n\tau_{k_i}(\phi_{\alpha_i})
\right\rangle^{X^{\mathrm{di}},\mathrm{tw}}_{g,d},
\\
\mathcal D_{X^{\mathrm{di}}}^{\mathrm{tw}}(t;\hbar)
&=
\exp\!\left(
\sum_{g\geq0}\hbar^{g-1}
F_g^{X^{\mathrm{di}},\mathrm{tw}}(t)
\right).
\end{aligned}
\label{eq:corrected-diagonal-potentials}
\end{equation}

The Hodge insertion is compatible with the gluing and unit axioms.
The Hodge bundle also pulls back under forgetting a marked point. Consequently, multiplication by
$\Lambda_g^\vee(1)^{-2}$ defines a cohomological field theory with
the same metric and unit as the diagonal theory. These properties
of the Hodge bundle, and their applications in Gromov--Witten theory, are reviewed in \cite{faberHodgeIntegralsGromovWitten2000}.

\begin{proposition}
\label{prop:Hodge-twisted-R}
The twisted diagonal theory has the same genus-zero descendant
potential, Frobenius structure, and $S$-calibration as
$X^{\mathrm{di}}$. Its $R$-matrix is
\begin{equation}
R^{X^{\mathrm{di}},\mathrm{tw}}(z)
=
\exp\!\left[
2\sum_{m\geq1}\frac{B_{2m}}{2m(2m-1)}z^{2m-1}
\right]R^{X^{\mathrm{di}}}(z)
=\widetilde R^{X^{\mathrm{di}}}(z).
\label{eq:Hodge-twisted-R}
\end{equation}
\end{proposition}

\begin{proof}
Since $\mathbb E_0$ has rank zero, $\Lambda_0^\vee(1)=1$.
Thus \eqref{eq:corrected-diagonal-correlators} agrees with the untwisted diagonal invariant in genus zero. In particular,
\begin{equation}
F_0^{X^{\mathrm{di}},\mathrm{tw}}(t)
=F_0^{X^{\mathrm{di}}}(t),
\qquad
S^{X^{\mathrm{di}},\mathrm{tw}}(z)
=S^{X^{\mathrm{di}}}(z).
\label{eq:Hodge-genus-zero-invariance}
\end{equation}
The Frobenius structure is therefore unchanged as well.

To compute the effect on the $R$-matrix, let $x_1,\ldots,x_g$ be
the Chern roots of $\mathbb E_g$. Then
\begin{equation}
\begin{aligned}
\log\Lambda_g^\vee(1)^{-2}
&=-2\sum_{j=1}^g\log(1-x_j)
\\
&=2\sum_{k\geq1}(k-1)!\operatorname{ch}_k(\mathbb E_g)
\\
&=2\sum_{m\geq1}(2m-2)!
\operatorname{ch}_{2m-1}(\mathbb E_g).
\end{aligned}
\label{eq:Hodge-Chern-character}
\end{equation}
The last equality uses Mumford's relation, which implies
$\operatorname{ch}_{2m}(\mathbb E_g)=0$ for $m\geq1$; see
\cite[Section 1]{faberHodgeIntegralsGromovWitten2000}.

Givental's Hodge integral formula \cite[Theorem 4.1]{giventalGromovWittenInvariantsQuantization2001}
expresses the insertion of
$$\exp(c\,\operatorname{ch}_{2m-1}(\mathbb E_g))$$ in the integral by the quantized
scalar symplectic transformation for the partition function
\[
\exp\left(
c\,\frac{B_{2m}}{(2m)!}\widehat{z^{2m-1}I}
\right).
\]
This is the quantized form of the Hodge integral equations of
Faber--Pandharipande \cite{faberHodgeIntegralsGromovWitten2000}. Their derivation uses Mumford's
Grothendieck--Riemann--Roch formula and the splitting axioms, and
applies to the localized equivariant theory considered here.
The action is taken in the dilaton-shifted variables, so the
translation preserving the unit is included.

Applying this formula to \eqref{eq:Hodge-Chern-character} multiplies
the $R$-matrix by
\begin{equation}
\exp\!\left[
\sum_{m\geq1}
\frac{B_{2m}}{(2m)!}\,2(2m-2)!\,z^{2m-1}
\right]I
=
\exp\!\left[
2\sum_{m\geq1}
\frac{B_{2m}}{2m(2m-1)}z^{2m-1}
\right]I.
\label{eq:Hodge-R-factor}
\end{equation}
This is the inverse of the scalar factor in
\eqref{eq:di-geometric-R-correction}. Proposition~\ref{prop:geometric-R-matrices}
therefore gives \eqref{eq:Hodge-twisted-R}.
\end{proof}

\begin{remark}
By Givental's Hodge integral formula
\cite[Theorem~4.1]{giventalGromovWittenInvariantsQuantization2001},
\eqref{eq:Hodge-Chern-character} gives the following relations
between the total descendant potentials:
\begin{equation}
\begin{aligned}
\mathcal D_{X^{\mathrm{di}}}^{\mathrm{tw}}(t;\hbar)
&=
\exp\!\left(
2\sum_{m\geq1}
\frac{B_{2m}}{2m(2m-1)}
\widehat{z^{2m-1}I}
\right)
\mathcal D_{X^{\mathrm{di}}}(t;\hbar),
\\
\mathcal D_{X^{\mathrm{di}}}(t;\hbar)
&=
\exp\!\left(
-2\sum_{m\geq1}
\frac{B_{2m}}{2m(2m-1)}
\widehat{z^{2m-1}I}
\right)
\mathcal D_{X^{\mathrm{di}}}^{\mathrm{tw}}(t;\hbar).
\end{aligned}
\label{eq:Hodge-partition-transformation}
\end{equation}
Here the hat denotes Givental's quantization with respect to
the pairing $\eta^{X^{\mathrm{di}}}$ and the dilaton-shifted
coordinates $t^{\alpha,k}-\delta^\alpha_1\delta_{k1}$;
see \cite[Section~2]{giventalGromovWittenInvariantsQuantization2001}.
\end{remark}

It follows that the identity
\eqref{eq:di-calibration-transformation} and the diagonal descendant
transformation in \eqref{eq:descendant-transformations} apply
unchanged to the twisted theory. 
Consequently, the genus-zero diagonal identity proved by LYZZ also holds for the twisted theory.
The Hodge insertion modifies its
higher-genus invariants so that its $R$-matrix is exactly the matrix
obtained from the logarithmic phase in
\eqref{eq:oscillating-integral-phases}. These are the diagonal
invariants that will enter the topological recursion comparison in
the next section.


\section{Eynard--Orantin topological recursion and DOSS correspondence}
\label{sec:TR}

In this section, we relate the Gromov--Witten theories of Section~\ref{sec:geometry}
in the $\bP^1$ and the antidiagonal resolved conifold case, and 
the twisted  Gromov--Witten theory in the diagonal case as in the end of the preceding section
to Eynard--Orantin topological recursion. We first recall the
recursion and specify our conventions. We then describe the DOSS
correspondence \cite{dunin-barkowskiIdentificationGiventalFormula2014} and the extension in \cite{xiongEYNARDORANTINTOPOLOGICALRECURSION} that interprets
changes of $y$ in terms of ancestor backgrounds. Finally, we apply
this correspondence to the spectral curves of $\mathbb P^1$ and
the resolved conifold, and compare their recursion multidifferentials.

\subsection{Eynard--Orantin topological recursion}
\label{subsec:EO}

We recall the Eynard--Orantin recursion
\cite{eynardInvariantsAlgebraicCurves2007} in the sphere setting
needed below. The spectral data consist of a Riemann surface
$\Sigma$, functions $x,y$, and a Bergman kernel $B$.
All spectral curves considered in this subsection are parametrized by
a coordinate $p$ on the Riemann sphere, whose Bergman kernel is
\begin{equation}
B(p,q)=\frac{dp\,dq}{(p-q)^2}.
\label{eq:Bergman-kernel}
\end{equation}
It is symmetric, holomorphic away from the diagonal, and has a double pole along the diagonal with leading coefficient $1$.
Equivalently, in any local coordinate $\zeta$,
\[
B(p,q)
=\left[
\frac{1}{(\zeta(p)-\zeta(q))^2}
+\text{a holomorphic function}
\right]d\zeta(p)\,d\zeta(q).
\]

Assume that $dx$ has finitely many simple zeros $p_i$, that $x$ and
$y$ are holomorphic near each $p_i$, and that $dy(p_i)\ne0$.
Write $U_i=x(p_i)$. The local Airy coordinate $\zeta_i$ and the local
involution $\iota_i$ are characterized by
\begin{equation}
x=U_i+\frac{\zeta_i^2}{2},
\qquad
\zeta_i(\iota_i(p))=-\zeta_i(p).
\label{eq:Airy-coordinates}
\end{equation}
Thus $x(\iota_i(p))=x(p)$ and $\iota_i(p_i)=p_i$.
For the functions in \eqref{eq:oscillating-integral-phases}, these
are the critical values and local coordinates used in
Subsection~\ref{subsec:R-comparison}, with critical point labels $i\in\{+,-\}$. The logarithmic functions are
understood through their holomorphic germs at the critical points.

The initial differentials are
\begin{equation}
\omega_{0,1}(p)=y(p)\,dx(p),
\qquad
\omega_{0,2}(p,q)=B(p,q).
\label{eq:EO-initial-differentials}
\end{equation}
The recursion kernel at $p_i$ is
\begin{equation}
K_i(p,q)
=-
\frac{\displaystyle\int_{\iota_i(q)}^q B(p,\cdot)}
{2\bigl(y(q)-y(\iota_i(q))\bigr)\,dx(q)}.
\label{eq:EO-kernel}
\end{equation}
The integral is taken in the local neighborhood of $p_i$.
For $g\geq0$, $n\geq1$, and $2g-2+n>0$, the stable differentials
are defined recursively by
\begin{equation}
\begin{aligned}
\omega_{g,n}(p_1,\ldots,p_n)
={}&\sum_i\operatorname*{Res}_{q=p_i}K_i(p_1,q)
\Bigg[
\omega_{g-1,n+1}(q,\iota_i(q),p_2,\ldots,p_n)
\\
&\quad+
\sum_{\substack{g_1+g_2=g\\I\sqcup J=\{2,\ldots,n\}}}^{\prime}
\omega_{g_1,|I|+1}(q,p_I)
\omega_{g_2,|J|+1}(\iota_i(q),p_J)
\Bigg].
\end{aligned}
\label{eq:EO}
\end{equation}
Here $p_I$ denotes the collection of variables with indices in $I$.
The primed sum $\sum^\prime$ runs over ordered decompositions and excludes
$(g_1,|I|)=(0,0)$ and $(g_2,|J|)=(0,0)$; factors of type
$(0,2)$ are included. The first term is omitted when $g=0$.

Although \eqref{eq:EO} singles out $p_1$, the resulting
$\omega_{g,n}$ is symmetric in all its arguments. 
It is a meromorphic differential in each variable, with poles only at the critical points of $x$ and zero residues there \cite{eynardInvariantsAlgebraicCurves2007}.
Our sign convention in \eqref{eq:EO-kernel} gives, for the Airy
curve $x=p^2/2$, $y=p$,
\[
\omega_{0,3}(p_1,p_2,p_3)
=\frac{dp_1\,dp_2\,dp_3}{p_1^2p_2^2p_3^2}.
\]

For stable differentials, the dependence on $y$ enters through
its odd part in the local coordinate \eqref{eq:Airy-coordinates} \cite{eynardInvariantsSpectralCurves2014,xiongEYNARDORANTINTOPOLOGICALRECURSION}:
\begin{equation}
\frac{y(\zeta_i)-y(-\zeta_i)}2
=\sum_{k\geq0}y_{2k+1}^{\,i}\zeta_i^{2k+1},
\qquad y_1^{\,i}\ne0.
\label{eq:EO-odd-y}
\end{equation}
Also, for a nonzero constant $c$,
\begin{equation}
y\longmapsto cy
\qquad\Longrightarrow\qquad
\omega_{g,n}\longmapsto c^{\,2-2g-n}\omega_{g,n},
\qquad n\geq1,\quad 2g-2+n>0.
\label{eq:EO-y-scaling}
\end{equation}


\subsection{The DOSS correspondence and variation of \texorpdfstring{$y$}{y}}
\label{subsec:DOSS}

We recall the correspondence of Dunin-Barkowski--Orantin--Shadrin--Spitz (DOSS) \cite{dunin-barkowskiIdentificationGiventalFormula2014}
in the formulation of Xiong--Zhou
\cite[Section 5]{xiongEYNARDORANTINTOPOLOGICALRECURSION}.
We use the local coordinates and recursion convention of Section~\ref{subsec:EO}, together with the convention $\Phi R(z)e^{U/z}$ of Section~\ref{subsec:R-matrices}.

Recall the cohomological field theory introduced in Section \ref{subsec:cohft}. Let $v$ denote the fixed semisimple primary point of the chosen theory. The ancestor potentials are denoted by $\overline F_g(\bar t;v)$. All other Frobenius quantities and one-forms
are evaluated at the same primary point $v$ unless indicated otherwise. 
Let $\Omega_{g,n;v}$ be the cohomological field theory
after the primary shift to $v$.

In the following paragraph, Greek indices refer to the flat basis $\partial_\alpha$, and Latin indices
refer to the normalized canonical basis. The latter are raised and
lowered using $\delta_{ij}$.

Expand the Bergman kernel in the local Airy coordinates $\zeta_i$ and $\zeta'_j$ at the critical
points $p_i$ and $p_j$, respectively:
\begin{equation}
B(\zeta_i,\zeta'_j)
=\left[
\frac{\delta_{ij}}{(\zeta_i-\zeta'_j)^2}
+\sum_{k,l\geq0}B^{ij}_{k,l}\zeta_i^k(\zeta'_j)^l
\right]d\zeta_i\,d\zeta'_j.
\label{eq:B-local}
\end{equation}
For correspondence with CohFTs, we need some additional properties of the coefficients $B^{ij}_{k,l}$ \cite{dunin-barkowskiIdentificationGiventalFormula2014,xiongEYNARDORANTINTOPOLOGICALRECURSION,xiongSpectralCurvesSemisimple2026}.

The compatibility with the $R$-matrix is expressed by
\begin{equation}
\sum_{k,l\geq0}(2k-1)!!(2l-1)!!B^{ij}_{2k,2l}z^kw^l
=
\left[
\frac{I-R(z)^{-1}R(w)^{-T}}{z+w}
\right]^{ij}.
\label{eq:B-R}
\end{equation}
In particular, setting one of the variables to zero gives
\begin{equation}
(R_{k+1})^i{}_j
=(-1)^k(2k-1)!!B^{ij}_{0,2k},
\qquad k\geq0.
\label{eq:R-boundary}
\end{equation}

We also impose the compatibility condition
\begin{equation}
(l+1)B^{ij}_{2r,l+2}
+(2r+1)B^{ij}_{2r+2,l}
+\delta_{ab}B^{ai}_{0,2r}B^{bj}_{0,l}=0,
\qquad r,l\geq0.
\label{eq:B-full-compatibility}
\end{equation}
This is the local condition of
\cite[Section 5.1]{xiongEYNARDORANTINTOPOLOGICALRECURSION}.
It includes odd values of $l$, whereas \eqref{eq:B-R} only concerns
the even coefficients in both variables. These conditions will be
verified for the three spectral curves in Subsection~\ref{subsec:three-curves}.

Define the one-forms $\Xi_k^i$ and the ancestor one-forms $W_k^i$ by
\begin{equation}
\Xi_k^i(p)
=(2k-1)!!\,[\zeta_i^{2k}]
\frac{B(\zeta_i,p)}{d\zeta_i},
\qquad
W_k^i(p)
=\sum_{\ell=0}^k(R_\ell)^i{}_j\,\Xi_{k-\ell}^j(p).
\label{eq:two-form-bases}
\end{equation}
The one-forms $\Xi_k^i$ are extracted directly from the Bergman kernel;
the ancestor one-forms $W_k^i$ incorporate the change of variables in
Givental's formula and correspond to canonical ancestor parameters.
Under \eqref{eq:B-full-compatibility}, they satisfy
\begin{equation}
\begin{aligned}
-d\!\left(\frac{\Xi_k^i}{dx}\right)
&=\Xi_{k+1}^i
+(2k-1)!!\,\delta_{ab}B^{ai}_{0,2k}\Xi_0^b,
\\
W_0^i&=\Xi_0^i,
\qquad
W_{k+1}^i=-d\!\left(\frac{W_k^i}{dx}\right).
\end{aligned}
\label{eq:ancestor-form-recursions}
\end{equation}
The first identity follows by expanding
\eqref{eq:B-full-compatibility} in the second local coordinate;
the second follows from \eqref{eq:R-boundary} and the symplectic condition of the $R$-matrix.

At the fixed primary point $v$, we denote by $y^\circ$ a reference
function whose odd part satisfies the unit compatibility condition
\begin{equation}
\begin{pmatrix}
\displaystyle\sum_{k\geq0}(2k+1)!!\,y_{2k+1}^{\circ,+}z^k\\[4pt]
\displaystyle\sum_{k\geq0}(2k+1)!!\,y_{2k+1}^{\circ,-}z^k
\end{pmatrix}
=R(z)^{-1}\Phi^{-1}\begin{pmatrix}1\\0\end{pmatrix}.
\label{eq:y-unit}
\end{equation}
Equivalently, by \eqref{eq:Laplace-odd-series},
\begin{equation}
\bm{\mathcal G}_{-z}(dy^\circ)
=R(z)^{-1}\Phi^{-1}\begin{pmatrix}1\\0\end{pmatrix},
\qquad y_1^{\circ,i}=\Delta_i^{-1/2}.
\label{eq:reference-y-period}
\end{equation}
The condition determines the local odd coefficients of $y^\circ$. This reference choice represents zero
ancestor background at the fixed primary point $v$; see below.

\begin{theorem}[DOSS correspondence]
\label{thm:DOSS}
Suppose that \eqref{eq:B-R}, \eqref{eq:B-full-compatibility}, and
\eqref{eq:y-unit} hold for the $R$-matrix of the chosen theory.
Then the recursion differentials of $(x,y^\circ,B)$ satisfy, for $n\geq1$ and $2g-2+n>0$,
\begin{equation}
\begin{aligned}
\omega_{g,n}[y^\circ](p_1,\ldots,p_n)
={}&\sum_{k_1,\ldots,k_n\geq0}
\left\langle
\prod_{r=1}^n\bar\tau_{k_r}(\partial_{\alpha_r})
\right\rangle_{g,n;v}
\prod_{r=1}^n
\Phi^{\alpha_r}{}_{i_r}\,W_{k_r}^{i_r}(p_r).
\end{aligned}
\label{eq:ancestor-TR}
\end{equation}
\end{theorem}

This is \cite[Theorem 4.1]{dunin-barkowskiIdentificationGiventalFormula2014}; see also \cite[Section 5.6]{xiongEYNARDORANTINTOPOLOGICALRECURSION}.
The factors of $\Phi$ express the transformation from the normalized canonical basis to the flat basis. 
The sum in \eqref{eq:ancestor-TR} is finite, since a nonzero ancestor correlator requires $k_1+\cdots+k_n\leq3g-3+n$.

We next replace the reference function $y^\circ$ by a general function $y$, keeping the primary point $v$, $x$, $B$, and the
identification with $R$ fixed.

For a general holomorphic $y$ with $y_1^{\,i}\ne0$, its odd
coefficients determine the column series
\begin{equation}
\mathbf y(z)
=\sum_{k\geq0}(2k+1)!!
\begin{pmatrix}y_{2k+1}^{\,+}\\y_{2k+1}^{\,-}\end{pmatrix}
z^{k+1}.
\label{eq:odd-y-series}
\end{equation}
By \eqref{eq:Laplace-odd-series}, this is precisely
\begin{equation}
\mathbf y(z)=z\,\bm{\mathcal G}_{-z}(dy).
\label{eq:effective-y-periods}
\end{equation}
In particular, the reference series is
$\mathbf y^\circ(z)=zR(z)^{-1}\Phi^{-1}(1,0)^T$.

Define the flat ancestor background point $a$ by
\begin{equation}
\begin{pmatrix}
\displaystyle\sum_{k\geq0}a^{1,k}z^k
\\[4pt]
\displaystyle\sum_{k\geq0}a^{2,k}z^k
\end{pmatrix}
=z\begin{pmatrix}1\\0\end{pmatrix}
-\Phi R(z)\mathbf y(z)
=\Phi R(z)\bigl(\mathbf y^\circ(z)-\mathbf y(z)\bigr).
\label{eq:y-background}
\end{equation}
In particular,
\begin{equation}
a^{\alpha,0}=0,
\qquad
a^{\alpha,1}=\delta^\alpha_1-\Phi^\alpha{}_i y_1^{\,i}.
\label{eq:y-background-leading}
\end{equation}
Then $a=0$ holds exactly when $y$ and $y^\circ$ have the same local odd parts. Varying $y$ changes the non-primary ancestor variables $\bar{t}^{\alpha,k},\ k\geq 1$ while
the chosen semisimple primary point $v$ remains fixed.

\begin{proposition}[\cite{xiongEYNARDORANTINTOPOLOGICALRECURSION}]
\label{prop:y-variation}
Assume \eqref{eq:B-R} and \eqref{eq:B-full-compatibility}, and let
$y$ be a general holomorphic function with $y^i_1\ne 0$. For $n\geq1$ and $2g-2+n>0$, the differentials
produced by topological recursion satisfy
\begin{equation}
\begin{aligned}
\omega_{g,n}[y](p_1,\ldots,p_n)
={}&\sum_{k_1,\ldots,k_n\geq0}
\left.
\frac{\partial^n\overline F_g(\bar t;v)}
{\partial\bar t^{\alpha_1,k_1}\cdots
 \partial\bar t^{\alpha_n,k_n}}
\right|_{\bar t=a}
\prod_{r=1}^n
\Phi^{\alpha_r}{}_{i_r}\,W_{k_r}^{i_r}(p_r),
\end{aligned}
\label{eq:shifted-ancestor-TR}
\end{equation}
where $a$ is given by \eqref{eq:y-background}. The shifted evaluation is interpreted as a formal Taylor expansion around $a=0$.
\end{proposition}

The shifted evaluation point $a$ is often called the background of the topological recursion in physical context.

At the partition-function level, the formula of Xiong--Zhou also contains the genus-one normalization
\[
\left(\prod_{i\in\{+,-\}}
\frac{y_1^{\,i}}{y_1^{\circ,i}}
\right)^{1/24}.
\]
Its contribution to the genus-one vacuum term will be treated in Subsection~\ref{subsec:vacuum-correspondence}. 
In the comparison of the next section, we will keep track of both the ancestor forms \eqref{eq:two-form-bases} and the background point \eqref{eq:y-background}.

\subsection{The three spectral curves}
\label{subsec:three-curves}

Equip each pair $(x,y^\circ)$ in \eqref{eq:oscillating-integral-phases}
with the Bergman kernel $B(p,q)$ in \eqref{eq:Bergman-kernel}. For
$\mathbb P^1$, this is the spectral curve of Norbury--Scott
\cite{norburyGromovWittenInvariants2014,dunin-barkowskiDubrovinsSuperpotentialGlobal2019}, whose relation to
Gromov--Witten theory was established in
\cite{dunin-barkowskiIdentificationGiventalFormula2014,borotLoopEquationsGromovWitten2019}.
We verify the compatibility conditions of Section~\ref{subsec:DOSS}
for all three curves, using the independently normalized
$R$-matrices of Subsection~\ref{subsec:R-comparison}.

The critical points are given in \eqref{eq:phase-critical-points}.
In the local Airy coordinates \eqref{eq:Airy-coordinates}, the first
coefficient one-forms of the Bergman kernel are
\begin{equation}
\Xi_0^i(p)=W_0^i(p)
=\left.\frac{dp}{d\zeta_i}\right|_{\zeta_i=0}
\frac{dp}{(p-p_i)^2}.
\label{eq:zero-forms}
\end{equation}
The Airy coordinate derivatives are
\begin{equation}
\left(\left.\frac{dp}{d\zeta_i}\right|_{\zeta_i=0}\right)^2
=\frac{1}{x''(p_i)}
=
\begin{cases}
\displaystyle\frac{p_i}{2},&\mathbb P^1,\\[5pt]
\displaystyle\frac{p_i(1+p_i)}{2},&X^{\mathrm{ad}},\\[5pt]
\displaystyle\frac{p_i(1+p_i)^2}{2},&X^{\mathrm{di}}.
\end{cases}
\label{eq:three-coordinate-derivatives}
\end{equation}

For each pair $(x,y^\circ)$ in
\eqref{eq:oscillating-integral-phases}, let $R(z)$ denote the
normalized stationary-phase matrix constructed in
\eqref{eq:1-dim-oscillating-integral}. Thus
\[
z\eta^{-1}df_j=e^{U_j/z}\Phi R(z)_{\bullet j},
\]
where $f_j$ is the scalar period
\eqref{eq:stationary-phase-periods}, and $d$ in this formula
denotes differentiation with respect to the primary variables.

\begin{lemma}
\label{lem:three-curve-compatibility}
Each of the three spectral curves $(x,y^\circ,B)$ satisfies \eqref{eq:B-R}, \eqref{eq:B-full-compatibility}, and \eqref{eq:y-unit} with the corresponding matrix $R(z)$.

Moreover, for any flat tangent vector
$h=h^\alpha\partial_\alpha$, we have
\begin{equation}
\bm{\mathcal G}_{-z}\!\left(
(h^\alpha\partial_\alpha x)\,dy^\circ
\right)
=R(z)^{-1}\Phi^{-1}h.
\label{eq:flat-negative-periods}
\end{equation}
Here $h$ on the right denotes its flat-coordinate column.
\end{lemma}

\begin{proof}
Put $r=e^{v^2}$ in the $\mathbb P^1$ case and $r=e^{u^2}$
in the resolved conifold cases, and write
\[
k_i=\left.\frac{dp}{d\zeta_i}\right|_{\zeta_i=0}.
\]

We first verify the differential compatibility condition
\eqref{eq:B-full-compatibility}. The derivatives
of the phase functions are
\begin{equation}
\begin{aligned}
\frac{1}{(x^{\mathbb P^1})'(p)}
&=\frac{p^2}{p^2-r},
\\
\frac{1}{(x^{X^{\mathrm{ad}}})'(p)}
&=\frac{p(1+p)(p+r)}{p^2+2rp+r},
\\
\frac{1}{(x^{X^{\mathrm{di}}})'(p)}
&=\frac{p(1+p)(p+r)}{p^2-r}.
\end{aligned}
\label{eq:three-reciprocal-derivatives}
\end{equation}
In each case, $F(p)=1/x'(p)$ has a polynomial part of degree
at most one and simple poles at $p_+$ and $p_-$, with
residue $k_i^2=1/x''(p_i)$ at $p_i$.
We claim that
\begin{equation}
d_p\!\left(\frac{B(p,q)}{dx(p)}\right)
+d_q\!\left(\frac{B(p,q)}{dx(q)}\right)
=-\delta_{ij}\Xi_0^i(p)\Xi_0^j(q).
\label{eq:B-differential}
\end{equation}
Indeed, the coefficient of $dp\,dq$ on the left is
\begin{equation}
\frac{F'(p)+F'(q)}{(p-q)^2}
-\frac{2(F(p)-F(q))}{(p-q)^3}.
\label{eq:B-rational-calculation}
\end{equation}
The polynomial part of $F$ contributes zero, while each
simple pole term $k_i^2/(p-p_i)$ contributes
\[
-\frac{k_i^2}{(p-p_i)^2(q-p_i)^2}.
\]
This proves \eqref{eq:B-differential} by
\eqref{eq:zero-forms}. Expanding at $p_i,p_j$ and comparing
the regular coefficients of
$\zeta_i^{2k}(\zeta'_j)^l$ gives
\eqref{eq:B-full-compatibility} for all $k,l\geq0$.

We next prove the equation~\eqref{eq:flat-negative-periods} and the condition~\eqref{eq:y-unit}.
Since $dy^\circ$ is independent of the primary variables when $p$ is fixed, differentiation of the stationary-phase periods gives
\[
e^{-U_j/z}z\partial_\alpha f_j
=
\mathcal G_{j,z}\!\left(
(\partial_\alpha x)\,dy^\circ
\right).
\]
Using the defining relation for $R(z)$ and
$\Phi^{-1}=\Phi^T\eta$, we obtain
\begin{equation}
\begin{aligned}
\mathcal G_{j,z}\!\left(
(h^\alpha\partial_\alpha x)\,dy^\circ
\right)
&=e^{-U_j/z}z h^\alpha\partial_\alpha f_j\\
&=h^T\eta\Phi R(z)_{\bullet j}\\
&=\bigl[R(z)^T\Phi^{-1}h\bigr]_j.
\end{aligned}
\label{eq:flat-positive-periods}
\end{equation}
Replacing $z$ by $-z$ and using
$R(-z)^T=R(z)^{-1}$ proves
\eqref{eq:flat-negative-periods}.
In particular, taking $h=\partial_1$ and using
$\partial_1x=1$ gives
\begin{equation}
\mathcal G_{i,z}(dy^\circ)
=
\left[
R(z)^T\Phi^{-1}
\begin{pmatrix}1\\0\end{pmatrix}
\right]_i.
\label{eq:unit-period-projection}
\end{equation}
Substituting $-z$ and applying
\eqref{eq:Laplace-odd-series} proves
\eqref{eq:y-unit}.

It remains to verify \eqref{eq:B-R}.
Taking $h=\Phi^\alpha{}_i\partial_\alpha$ in
\eqref{eq:flat-positive-periods} gives
\begin{equation}
R^i{}_j(z)
=
\mathcal G_{j,z}\!\left(
\bigl(\Phi^1{}_i+\Phi^2{}_i\partial_2x\bigr)\,dy^\circ
\right).
\label{eq:canonical-period-projection}
\end{equation}
Set
\[
g_i(p)=
\bigl(\Phi^1{}_i+\Phi^2{}_i\partial_2x(p)\bigr)
\frac{dy^\circ(p)}{dx(p)}.
\]
Substitution of the $\Phi$ from
Section~\ref{subsec:Frobenius} gives
\begin{equation}
\begin{aligned}
g_i^{\mathbb P^1}(p)
&=\frac{\Phi^1{}_i p+r\Phi^2{}_i}{p^2-r}
=\frac{k_i}{p-p_i},
\\[3pt]
g_i^{X^{\mathrm{ad}}}(p)
&=-\mathrm i\,
\frac{(1+p)\bigl(\Phi^1{}_i(p+r)-r\Phi^2{}_i\bigr)}
{p^2+2rp+r}
=-\mathrm i\,\Phi^1{}_i+\frac{k_i}{p-p_i},
\\[3pt]
g_i^{X^{\mathrm{di}}}(p)
&=\frac{(1+p)\bigl(\Phi^1{}_i(p+r)+r\Phi^2{}_i\bigr)}
{p^2-r}
=\Phi^1{}_i+\frac{k_i}{p-p_i}.
\end{aligned}
\label{eq:three-canonical-projection-functions}
\end{equation}
The quantities $\Phi$, $p_i$, and $k_i$ in each row
belong to that row's theory.
It follows that $dg_i=-\Xi_0^i$. Formal integration by parts
\eqref{eq:Gaussian-integration-by-parts} therefore gives
\begin{equation}
R^i{}_j(z)
=\mathcal G_{j,z}(g_i\,dx)
=z\,\mathcal G_{j,z}(\Xi_0^i).
\label{eq:R-period-form}
\end{equation}
The local expansion
\[
\Xi_0^i
=
\left(
\frac{\delta_{ij}}{\zeta_j^2}
+\sum_{l\geq0}B^{ij}_{0,l}\zeta_j^l
\right)d\zeta_j
\]
and the regularization \eqref{eq:Laplace-regularization}
then imply
\begin{equation}
R^i{}_j(z)
=
\delta^i_j+
\sum_{m\geq0}(-1)^m(2m-1)!!
B^{ij}_{0,2m}z^{m+1}.
\label{eq:geometric-R-from-B}
\end{equation}
Thus the boundary relation \eqref{eq:R-boundary} holds.

Taking the even coefficients in both variables in
\eqref{eq:B-differential}, and using
\eqref{eq:geometric-R-from-B}, gives
\begin{equation}
\sum_{k,l\geq0}(-1)^{k+l}(2k-1)!!(2l-1)!!
B^{ij}_{2k,2l}z^kw^l
=
\left[
\frac{R(z)^T R(w)-I}{z+w}
\right]^{ij}.
\label{eq:three-curves-B-R-signed}
\end{equation}
Replacing $z,w$ by $-z,-w$ and using
$R(-z)^T=R(z)^{-1}$ proves \eqref{eq:B-R}.
\end{proof}

With this lemma, one has the proposition:
\begin{proposition}
\label{prop:three-reconstructions}
The spectral curves $(x,y^\circ,B)$ associated with
$\mathbb P^1$, $X^{\mathrm{ad}}$, and $X^{\mathrm{di}}$
have the respective $R$-matrices
\begin{equation}
R^{\mathbb P^1}(z),
\qquad
R^{X^{\mathrm{ad}}}(z),
\qquad
\widetilde R^{X^{\mathrm{di}}}(z)
=R^{X^{\mathrm{di}},\mathrm{tw}}(z).
\label{eq:three-spectral-R-matrices}
\end{equation}
For $n\geq1$ and $2g-2+n>0$, their topological recursion
multidifferentials recover, through \eqref{eq:ancestor-TR},
the ancestor correlators of $\mathbb P^1$ at $v$,
the anti-diagonal resolved conifold at $u$, and the Hodge-twisted diagonal theory at $u$, respectively.
\end{proposition}

The next subsection compares the recursion kernels after replacing $y^\circ$ by the shifted functions $y_{\mathrm{sh}}$. 
The backgrounds $a$ corresponding to these shifted curves will be given in Section~\ref{sec:proof}.

\subsection{A common recursion sphere and change of spectral curves}
\label{subsec:common-kernel}


We now place the spectral curves of
Proposition~\ref{prop:three-reconstructions} on a common sphere
and modify the resolved conifold $y$-functions to compare their
recursion kernels.

On the Riemann sphere with coordinate $\ell$, set
\begin{equation}
\mathsf{x}(\ell)=v^1+\ell+\frac{e^{v^2}}{\ell},
\qquad
B(\ell_1,\ell_2)
=\frac{d\ell_1\,d\ell_2}{(\ell_1-\ell_2)^2}.
\label{eq:common-W}
\end{equation}
The critical points and critical values of $\mathsf{x}$ are
\begin{equation}
\ell_\pm=\pm e^{v^2/2},
\qquad
\mathsf{x}(\ell_\pm)=v^1\pm2e^{v^2/2}.
\label{eq:common-critical-points}
\end{equation}
We work on the locus $(v^1)^2\ne4e^{v^2}$, so that $\log \mathsf{x}$
is holomorphic near both critical points. Since
$d\log \mathsf{x}=d\mathsf{x}/\mathsf{x}$, the functions $\mathsf{x}$ and $\log \mathsf{x}$ have the same
simple critical points. Their local involutions are 
\begin{equation}
\ell\longmapsto\frac{e^{v^2}}{\ell}.
\label{eq:common-involution}
\end{equation}

Recall the Frobenius coordinate transformations
\eqref{eq:intro-Frobenius-coordinate-transformations}:
\begin{equation*}
\begin{aligned}
X^{\mathrm{di}}:\qquad
v^1&=e^{u^1}(1+e^{u^2}),
&v^2&=2u^1+u^2,
\\
X^{\mathrm{ad}}:\qquad
v^1&=\mathrm i\,e^{u^1}(1-2e^{u^2}),
&v^2&=2u^1+u^2+\log(1-e^{u^2}).
\end{aligned}
\end{equation*}
The corresponding changes of sphere coordinate satisfy
\begin{equation*}
\begin{aligned}
X^{\mathrm{di}}:\qquad
\ell&=e^{u^1}p,
&\mathsf{x}(\ell)&=e^{x^{X^{\mathrm{di}}}(u;p)},
\\
X^{\mathrm{ad}}:\qquad
\ell&=\mathrm i\,e^{u^1}(p+e^{u^2}),
&\mathsf{x}(\ell)&=\mathrm i\,e^{x^{X^{\mathrm{ad}}}(u;p)}.
\end{aligned}
\end{equation*}
These affine changes preserve $B$.
We express all curves in the common coordinate $\ell$
and suppress additive constants in $x$.
Figure~\ref{fig:ancestor-shifts-and-comparison} summarizes the relations between these spectral curves.

By Proposition~\ref{prop:y-variation}, the vertical arrows in
Figure~\ref{fig:ancestor-shifts-and-comparison} represent shifts of the ancestor
background at fixed $x$ and $B$.

At each critical point $\ell_i$, we have
\[
\left.\frac{d y_{\mathrm{sh}}^{\mathrm{di}}}{d\ell}\right|_{\ell_i}
=
-\left.\frac{d y_{\mathrm{sh}}^{\mathrm{ad}}}{d\ell}\right|_{\ell_i}
=
\frac{\mathsf{x}(\ell_i)}{\ell_i}\ne0.
\]
Thus the three curves
$\mathcal C^{\mathbb P^1}$,
$\mathcal C_{\mathrm{sh}}^{\mathrm{di}}$, and
$\mathcal C_{\mathrm{sh}}^{\mathrm{ad}}$
satisfy the regularity assumptions of Section~\ref{subsec:EO}.
We denote their recursion differentials by
$\omega^{\mathbb P^1}_{g,n}$, $\omega^+_{g,n}$, and
$\omega^-_{g,n}$, respectively.

\begin{proposition}
\label{prop:kernel-comparison}
For every $g\geq0$, $n\geq1$, and $2g-2+n>0$, one has
\begin{equation}
\begin{aligned}
\omega^+_{g,n}(\ell_1,\ldots,\ell_n)
&=\omega^{\mathbb P^1}_{g,n}(\ell_1,\ldots,\ell_n),
\\
\omega^-_{g,n}(\ell_1,\ldots,\ell_n)
&=(-1)^n\omega^{\mathbb P^1}_{g,n}(\ell_1,\ldots,\ell_n).
\end{aligned}
\label{eq:stable-differentials}
\end{equation}
\end{proposition}

\begin{proof}
Write $\bar\ell=e^{v^2}/\ell$ near either critical point.
Then $\mathsf{x}(\bar\ell)=\mathsf{x}(\ell)$, and hence
\begin{equation}
\begin{aligned}
&\bigl(\mathsf{x}(\ell)\log\ell
      -\mathsf{x}(\bar\ell)\log\bar\ell\bigr)\,d\log \mathsf{x}(\ell)
\\
&\qquad
=\mathsf{x}(\ell)\bigl(\log\ell-\log\bar\ell\bigr)
  \frac{d\mathsf{x}(\ell)}{\mathsf{x}(\ell)}
\\
&\qquad
=\bigl(\log\ell-\log\bar\ell\bigr)\,d\mathsf{x}(\ell).
\end{aligned}
\label{eq:kernel-cancellation}
\end{equation}
Thus the denominators in \eqref{eq:EO-kernel} agree for $\mathcal C^{\mathbb P^1}$ and $\mathcal C_{\mathrm{sh}}^{\mathrm{di}}$.
Their numerators also agree, because both the Bergman
kernel and the local involution are the same. They therefore have identical recursion kernels and the same initial differential $\omega_{0,2}=B$.
Induction in \eqref{eq:EO} proves the first
identity in \eqref{eq:stable-differentials}.

The curve $\mathcal C_{\mathrm{sh}}^{\mathrm{ad}}$ is obtained
from $\mathcal C_{\mathrm{sh}}^{\mathrm{di}}$ by replacing
$y_{\mathrm{sh}}^{\mathrm{di}}$ with
$y_{\mathrm{sh}}^{\mathrm{ad}}=-y_{\mathrm{sh}}^{\mathrm{di}}$. By \eqref{eq:EO-y-scaling}, its stable differentials
are multiplied by
\(
(-1)^{2-2g-n}=(-1)^n,
\)
which proves the second identity.
\end{proof}

\begin{remark}
Throughout the comparison on the common sphere, we label all
canonical data by $\ell_\pm$. Since
$\ell(p_\pm^{X^{\mathrm{ad}}})=\ell_\mp$, we interchange the
original anti-diagonal canonical labels and retain the same notation.
\end{remark}

Proposition~\ref{prop:kernel-comparison} gives an equality of stable differentials on the common sphere. 
To translate Proposition~\ref{prop:kernel-comparison} into
descendant identities, we must identify the ancestor backgrounds
represented by $y_{\mathrm{sh}}^{\mathrm{di}}$ and
$y_{\mathrm{sh}}^{\mathrm{ad}}$ and compare the ancestor one-forms
under $\mathsf{x}\mapsto\log \mathsf{x}$.
These two steps are carried out in Section~\ref{sec:proof}.

\section{Proof of the correspondences}
\label{sec:proof}

In this section, we prove the anti-diagonal correspondence and its Hodge-twisted
diagonal analogue by combining the calibration identities of
Section~\ref{sec:geometry} with the topological recursion comparison
of Section~\ref{sec:TR}. 
The genus-zero parts of the LYZZ correspondences were proved in \cite{liuEquivariantGromovWitten2022}.
We therefore concentrate on positive genus. 

We first determine the ancestor backgrounds induced by the descendant transformations and show that, under Proposition~\ref{prop:y-variation}, they are represented by the shifted $y$-functions introduced in Section~\ref{subsec:common-kernel}.
We then define a family of linear functionals $\lambda_{\alpha,k}$ on the span of the $\mathbb P^1$ ancestor one-forms on the common sphere. Applying these functionals extracts descendant correlators from the multidifferentials of the $\mathbb P^1$ spectral curve. 
We next compute the action of these functionals on the multidifferentials of the shifted spectral curves of resolved conifold. The results recover derivatives of the transformed resolved conifold descendant potentials.
Combining these with the comparison of recursion differentials in Proposition~\ref{prop:kernel-comparison} proves the correspondences with at least one marked point in the stable range.
Subsection~\ref{subsec:vacuum-correspondence} deals with the remaining vacuum terms.
 
Throughout this section, the diagonal theory is understood to include
the Hodge twist introduced in Section~\ref{subsec:Hodge-correction}.

\subsection{Descendant and ancestor backgrounds}
\label{subsec:backgrounds}





\begin{figure}[htbp]
\centering
\begin{tikzpicture}[
  >=Latex,
  font=\small,
  stage/.style={
    align=center,
    text width=2.7cm,
    inner sep=4pt
  },
  flow/.style={
    ->,
    line width=0.55pt,
    shorten <=3pt,
    shorten >=3pt
  },
  dependence/.style={flow, dashed},
  lab/.style={
    font=\scriptsize,
    align=center,
    inner sep=3pt
  }
]
  \node[stage] (s) at (0,0)
    {{\footnotesize $\mathbb P^1$ descendant}\\[6pt]
     $s_v^{\alpha,k}=v^\alpha\delta_{k0}$};

  \node[stage] (t) at (4.7,0)
    {{\footnotesize $X^\nu$ descendant}\\[6pt]
     $t_\nu(s_v)$};

  \node[stage] (a) at (9.4,0)
    {{\footnotesize $X^\nu$ ancestor}\\[6pt]
     $a_\nu(z)$};

  \node[stage] (v) at (0,-2.7)
    {{\footnotesize $\mathbb P^1$ primary}\\[6pt]
     $v$};

  \node[stage] (u) at (4.7,-2.7)
    {{\footnotesize $X^\nu$ primary}\\[6pt]
     $u$};

  \draw[flow] (v.north) --
    node[lab,right] {
      embedding\\
      \eqref{eq:source-centers}
    }
    (s.south);

  \draw[flow] (s.east) --
    node[lab,above] {LYZZ}
    node[lab,below] {\eqref{eq:descendant-transformations}}
    (t.west);

  \draw[flow] (t.east) --
    node[lab,above] {$S^{X^\nu}(u;z)$}
    node[lab,below] {\eqref{eq:target-ancestor-inputs}}
    (a.west);

  \draw[flow] (v.east) --
    node[lab,above] {coordinate change}
    node[lab,below] {
      \eqref{eq:di-coordinate-transformation}
      or \eqref{eq:ad-coordinate-transformation}
    }
    (u.west);

  \draw[dependence] (u.east) to[out=0,in=-90]
    node[lab,below right,pos=0.45] {fixed base point}
    (a.south);
\end{tikzpicture}

\caption{Background transformations for
$\nu=\mathrm{di},\mathrm{ad}$.
The dashed arrow records the dependence of the ancestor
transformation on the fixed primary point $u$.}
\label{fig:background-transformations}
\end{figure}
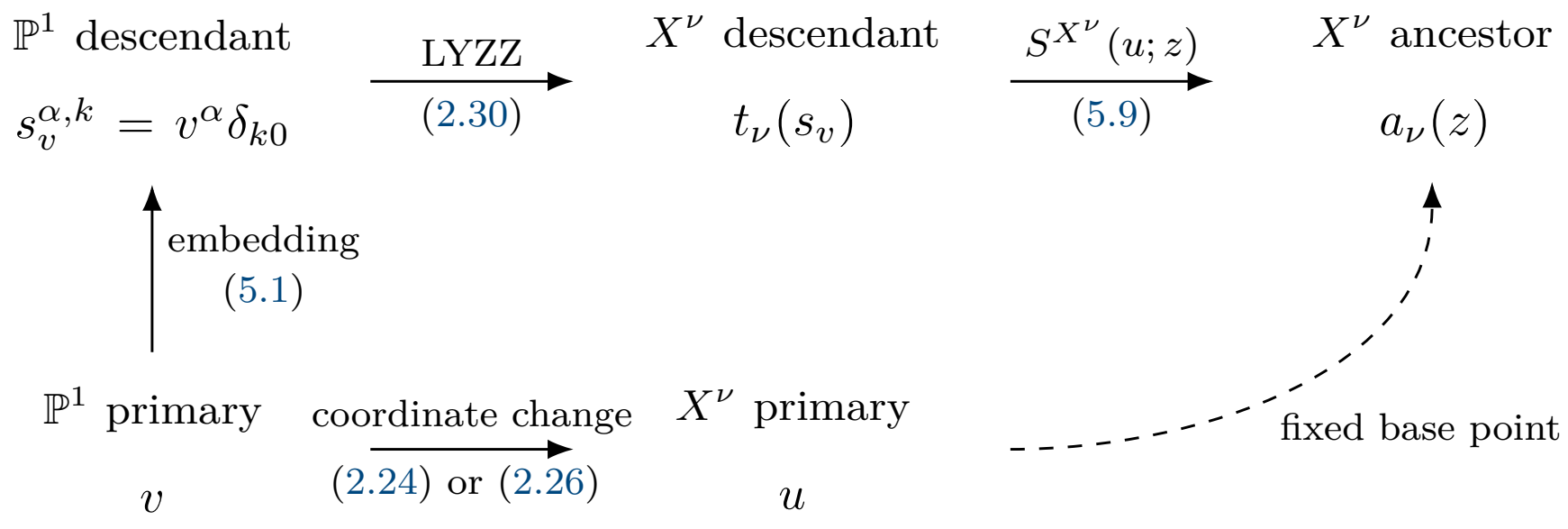

Fix a semisimple primary point $v$ with $(v^1)^2\ne4e^{v^2}$.
For each comparison, choose the corresponding point $u$ on a local
branch of \eqref{eq:di-coordinate-transformation} or
\eqref{eq:ad-coordinate-transformation}. Define the pure primary
descendant point $s_v$ by
\begin{equation}
s_v^{\alpha,k}=v^\alpha\delta_{k0}.
\label{eq:source-centers}
\end{equation}
Figure~\ref{fig:background-transformations} shows the transformations between each backgrounds.

Recall that $\overline F_g^{\mathbb P^1}(\bar t;v)$ is the
ancestor potential of the CohFT shifted to $v$.
The original $\mathbb P^1$ spectral curve determines its
derivatives at $\bar t=0$.
Under the ancestor--descendant transformation, this background
corresponds to the descendant point
$s^{\alpha,k}_v=v^\alpha\delta_{k0}$, as verified below.

We now verify that $s_v$ corresponds to zero ancestor background at the fixed primary point $v$.
Writing
\[
s(z)=\sum_{k\geq0}s^{\alpha,k}z^k\partial_\alpha,
\]
the ancestor--descendant transformation is
\begin{equation}
a^{\mathbb P^1}(v;z)
=
\bigl[
S^{\mathbb P^1}(v;z)
\bigl(s(z)-z\partial_1\bigr)
\bigr]_+
+z\partial_1,
\label{eq:source-ancestor-input}
\end{equation}
where $[\;\cdot\;]_+$ retains the nonnegative powers of $z$.
This is the change of variables in Givental's
ancestor--descendant formula
\cite{giventalGromovWittenInvariantsQuantization2001,
kontsevichRelationsCorrelatorsTopological1998}.
The calibration \eqref{eq:P1-calibration} gives
$\theta_{1,1}^{\mathbb P^1}(v)=v^1v^2$, and hence
\[
S_1^{\mathbb P^1}(v)\partial_1
=
(\eta^{\mathbb P^1})^{\alpha\beta}
\partial_\beta\theta_{1,1}^{\mathbb P^1}(v)\,
\partial_\alpha
=
v^\alpha\partial_\alpha.
\]
When $s_v(z)=v^\alpha\partial_\alpha$, it follows that
\begin{equation}
\left.a^{\mathbb P^1}(v;z)\right|_{s=s_v}
=
v^\alpha\partial_\alpha
-S_1^{\mathbb P^1}(v)\partial_1
=0.
\label{eq:source-zero-ancestor}
\end{equation}
More generally,
\begin{equation}
a^{\mathbb P^1}(v;z)
=
\bigl[
S^{\mathbb P^1}(v;z)
\bigl(s(z)-s_v(z)\bigr)
\bigr]_+.
\label{eq:source-centered-ancestor}
\end{equation}

Since $a^{\mathbb P^1}(v;z)$ is affine in $s$, we have
\[
\frac{\partial a^{\beta,m}}{\partial s^{\alpha,k}}
=
\bigl(S_{k-m}^{\mathbb P^1}(v)\bigr)^\beta{}_\alpha,
\qquad
S_j^{\mathbb P^1}=0\quad(j<0),
\]
where $a^{\beta,m}$ denotes the coefficient of
$z^m\partial_\beta$ in $a^{\mathbb P^1}(v;z)$.
Together with $a^{\mathbb P^1}(v;z)|_{s=s_v}=0$,
the ancestor--descendant formula therefore gives
\begin{equation}
\begin{aligned}
&\left.
\frac{\partial^n F_g^{\mathbb P^1}(s)}
{\partial s^{\alpha_1,k_1}\cdots
 \partial s^{\alpha_n,k_n}}
\right|_{s=s_v}
\\
&\quad=
\sum_{m_1,\ldots,m_n\geq0}
\left.
\frac{\partial^n
\overline F_g^{\mathbb P^1}(\bar t;v)}
{\partial\bar t^{\beta_1,m_1}\cdots
 \partial\bar t^{\beta_n,m_n}}
\right|_{\bar t=0}
\prod_{r=1}^n
\bigl(S_{k_r-m_r}^{\mathbb P^1}(v)\bigr)
^{\beta_r}{}_{\alpha_r},
\end{aligned}
\label{eq:source-ancestor-descendant-derivatives}
\end{equation}
for $n\geq1$ and $2g-2+n>0$, with repeated $\beta_r$ indices
summed. 
The primary point $v$ is held fixed throughout.

For the resolved conifold theory, at a fixed primary point $u$, write
\[
t(z)=\sum_{k\geq0}t^{\alpha,k}z^k\partial_\alpha.
\]
Substitution of $s_v$ into \eqref{eq:descendant-transformations} gives
\begin{equation}
\begin{aligned}
t_{\mathrm{di}}(s_v;z)-z\partial_1
&=v^1\partial_1
+\frac{v^2\partial_2-\partial_1}{1-z}
-\frac{2z}{(1-z)^2}\partial_2,\\
t_{\mathrm{ad}}(s_v;z)-z\partial_1
&=-\mathrm i v^1(\partial_1+\partial_2)
+\frac{v^2\partial_2-(\partial_1+\partial_2)}{1-z}
-\frac{2z}{(1-z)^2}\partial_2.
\end{aligned}
\label{eq:mapped-centers}
\end{equation}
In particular, the primary components are
\[
\begin{aligned}
\bigl(t_{\mathrm{di}}^{1,0}(s_v),t_{\mathrm{di}}^{2,0}(s_v)\bigr)
&=(v^1-1,\,v^2),\\
\bigl(t_{\mathrm{ad}}^{1,0}(s_v),t_{\mathrm{ad}}^{2,0}(s_v)\bigr)
&=(-\mathrm i v^1-1,\,-\mathrm i v^1+v^2-1).
\end{aligned}
\]
The constant terms arise from the dilaton shift in
\eqref{eq:descendant-transformations}.
In the diagonal case the primary points are related by
\begin{equation}
v^1=e^{u^1}(1+e^{u^2}),\qquad v^2=2u^1+u^2,
\label{eq:di-fixed-primary}
\end{equation}
and in the anti-diagonal case by
\begin{equation}
v^1=\mathrm i e^{u^1}(1-2e^{u^2}),\qquad
v^2=2u^1+u^2+\log(1-e^{u^2}).
\label{eq:ad-fixed-primary}
\end{equation}
Write $r=e^{u^2}$ in each case.

At fixed $u$, the descendant backgrounds in \eqref{eq:mapped-centers}
give the ancestor inputs
\begin{equation}
\begin{aligned}
a_{\mathrm{di}}(z;v)
&=\bigl[
S^{X^{\mathrm{di}}}(u(v);z)
\bigl(t_{\mathrm{di}}(s_v;z)-z\partial_1\bigr)
\bigr]_++z\partial_1,
\\
a_{\mathrm{ad}}(z;v)
&=\bigl[
S^{X^{\mathrm{ad}}}(u(v);z)
\bigl(t_{\mathrm{ad}}(s_v;z)-z\partial_1\bigr)
\bigr]_++z\partial_1.
\end{aligned}
\label{eq:target-ancestor-inputs}
\end{equation}

\begin{lemma}
\label{lem:calibration-center-sums}
At the corresponding primary points $u$ and $v$ one has
\begin{equation}
\begin{aligned}
\sum_{j\geq0}S_j^{X^{\mathrm{di}}}(u)\partial_2
&=e^{u^1}\bigl(r\partial_1+(1-r)\partial_2\bigr),\\
\sum_{j\geq0}S_j^{X^{\mathrm{di}}}(u)
  (\partial_1+2j\partial_2)
&=v^1\partial_1
  +v^2e^{u^1}\bigl(r\partial_1+(1-r)\partial_2\bigr),\\[3pt]
\sum_{j\geq0}S_j^{X^{\mathrm{ad}}}(u)\partial_2
&=e^{u^1}\bigl(-r\partial_1+(1-r)\partial_2\bigr),\\
\sum_{j\geq0}S_j^{X^{\mathrm{ad}}}(u)
  \bigl(\partial_1+(2j+1)\partial_2\bigr)
&=-\mathrm i v^1(\partial_1+\partial_2)
  +v^2e^{u^1}\bigl(-r\partial_1+(1-r)\partial_2\bigr).
\end{aligned}
\label{eq:calibration-center-sums}
\end{equation}
\end{lemma}

\begin{proof}
We evaluate the sums through the functions $\theta$ of Section~\ref{subsec:S-calibration}. Comparing
coefficients in \eqref{eq:scalar-matrix-calibration} gives
\[
(S_j)^\alpha{}_\beta
=
\eta^{\alpha\gamma}\partial_\gamma\theta_{\beta,j}.
\]
Consequently,
\[
\sum_{j\geq0}S_j(u)\partial_\beta
=
\left.
\eta^{\alpha\gamma}
\partial_\gamma\theta_\beta(u;\nu)
\right|_{\nu=1}\partial_\alpha.
\]

We first prove the first and third identities. Since
$(-1)_d=0$ for $d\geq2$, the scalar calibrations
\eqref{eq:di-calibration} and \eqref{eq:ad-calibration} give
\[
\begin{aligned}
\theta_2^{X^{\mathrm{di}}}(u;1)
&=e^{u^1}(1+e^{u^2})-1,
\\
\theta_2^{X^{\mathrm{ad}}}(u;1)
&=-e^{u^1}(1-e^{u^2})+1.
\end{aligned}
\]
Computing the gradients and applying the inverse metrics yields
\[
\begin{aligned}
\sum_{j\geq0}S_j^{X^{\mathrm{di}}}(u)\partial_2
&=
e^{u^1+u^2}\partial_1
+e^{u^1}(1-e^{u^2})\partial_2,
\\
\sum_{j\geq0}S_j^{X^{\mathrm{ad}}}(u)\partial_2
&=
-e^{u^1+u^2}\partial_1
+e^{u^1}(1-e^{u^2})\partial_2.
\end{aligned}
\]
For the remaining identities, introduce within this proof
\[
\begin{aligned}
\Theta_{\mathrm{di}}(u)
&=
\sum_{j\geq0}
\left(
\theta_{1,j}^{X^{\mathrm{di}}}(u)
+2j\theta_{2,j}^{X^{\mathrm{di}}}(u)
\right)
\\
&=
\left[
\theta_1^{X^{\mathrm{di}}}(u;\nu)
+2\nu\frac{\partial}{\partial\nu}
\theta_2^{X^{\mathrm{di}}}(u;\nu)
\right]_{\nu=1},
\\[4pt]
\Theta_{\mathrm{ad}}(u)
&=
\sum_{j\geq0}
\left(
\theta_{1,j}^{X^{\mathrm{ad}}}(u)
+(2j+1)\theta_{2,j}^{X^{\mathrm{ad}}}(u)
\right)
\\
&=
\left[
\theta_1^{X^{\mathrm{ad}}}(u;\nu)
+\left(2\nu\frac{\partial}{\partial\nu}+1\right)
\theta_2^{X^{\mathrm{ad}}}(u;\nu)
\right]_{\nu=1}.
\end{aligned}
\]
Thus,
\begin{equation}
\begin{aligned}
\sum_{j\geq0}S_j^{X^{\mathrm{di}}}(u)
\bigl(\partial_1+2j\partial_2\bigr)
&=
(\eta^{X^{\mathrm{di}}})^{\alpha\gamma}
\partial_\gamma\Theta_{\mathrm{di}}(u)\,
\partial_\alpha,
\\
\sum_{j\geq0}S_j^{X^{\mathrm{ad}}}(u)
\bigl(\partial_1+(2j+1)\partial_2\bigr)
&=
(\eta^{X^{\mathrm{ad}}})^{\alpha\gamma}
\partial_\gamma\Theta_{\mathrm{ad}}(u)\,
\partial_\alpha.
\end{aligned}
\label{eq:S-sum}
\end{equation}

For the diagonal case, the coefficients
\eqref{eq:di-b-coefficients} satisfy
\[
b_{1,1}^{\beta,j,\mathrm{di}}
=
\delta^\beta_1+2j\delta^\beta_2.
\]
Thus \eqref{eq:di-calibration-transformation}, with
$(\alpha,k)=(1,1)$, gives
\[
\Theta_{\mathrm{di}}(u)
=
\theta_{1,1}^{\mathbb P^1}(v(u))
=
v^1(u)v^2(u),
\]
where the equality follows from
\eqref{eq:P1-calibration}. The diagonal coordinate transformation
satisfies
\[
\frac{\partial(v^1,v^2)}{\partial(u^1,u^2)}=
\begin{pmatrix}
\partial_1v^1&\partial_1v^2\\
\partial_2v^1&\partial_2v^2
\end{pmatrix}
=
\begin{pmatrix}
e^{u^1}(1+e^{u^2})&2\\
e^{u^1+u^2}&1
\end{pmatrix}.
\]
By the chain rule,
\[
\begin{pmatrix}
\partial_1\Theta_{\mathrm{di}}\\
\partial_2\Theta_{\mathrm{di}}
\end{pmatrix}
=
\begin{pmatrix}
e^{u^1}(1+e^{u^2})&2\\
e^{u^1+u^2}&1
\end{pmatrix}
\begin{pmatrix}v^2\\v^1\end{pmatrix}.
\]
Applying the inverse diagonal metric to this column and using
\eqref{eq:S-sum} gives the second identity in \eqref{eq:calibration-center-sums}.

For the anti-diagonal case,
\eqref{eq:ad-b-coefficients} gives
\[
b_{1,1}^{\beta,j,\mathrm{ad}}
=
\delta^\beta_1+(2j+1)\delta^\beta_2.
\]
The gradient identity
\eqref{eq:ad-calibration-transformation}, again with
$(\alpha,k)=(1,1)$, therefore becomes
\[
\begin{pmatrix}v^2\\v^1\end{pmatrix}
=
-
\begin{pmatrix}
e^{-u^1}&-e^{-u^1}\\
\mathrm i e^{u^2}&\mathrm i(1-e^{u^2})
\end{pmatrix}
\begin{pmatrix}
\partial_1\Theta_{\mathrm{ad}}\\
\partial_2\Theta_{\mathrm{ad}}
\end{pmatrix}.
\]
Solving for the two derivatives gives
\[
\partial_1\Theta_{\mathrm{ad}}=\mathrm i v^1-e^{u^1}(1-r)v^2,
\qquad
\partial_2\Theta_{\mathrm{ad}}=\mathrm i v^1+e^{u^1}r v^2.
\]
Applying the inverse anti-diagonal metric gives the last identity
in \eqref{eq:calibration-center-sums}.
\end{proof}

Substituting this lemma into \eqref{eq:target-ancestor-inputs} gives
\begin{proposition}
At the target primary point $u$, the ancestor inputs corresponding to
$t_\nu(s_v)$ are
\begin{equation}
\begin{aligned}
a_{\mathrm{di}}(z;v)
&=z\partial_1-\frac{v^1z}{1-z}\partial_1
  -\frac{2ze^{u^1}}{(1-z)^2}
   \bigl(r\partial_1+(1-r)\partial_2\bigr),\\
a_{\mathrm{ad}}(z;v)
&=z\partial_1+\frac{\mathrm i v^1z}{1-z}(\partial_1+\partial_2)
  -\frac{2ze^{u^1}}{(1-z)^2}
   \bigl(-r\partial_1+(1-r)\partial_2\bigr).
\end{aligned}
\label{eq:ancestor-backgrounds}
\end{equation}
In particular, $a_{\mathrm{di}}^{\alpha,0}=a_{\mathrm{ad}}^{\alpha,0}=0$.
\end{proposition}

\begin{proof}
For the coefficient calculation write
\[
(c_{\mathrm{di}},h_{\mathrm{di}})=(v^1,\partial_1),
\qquad
(c_{\mathrm{ad}},h_{\mathrm{ad}})
=(-\mathrm i v^1,\partial_1+\partial_2).
\]
By \eqref{eq:mapped-centers}, the coefficient of $z^k$ in
$t_\nu(s_v;z)-z\partial_1$ is
$\delta_{k0}c_\nu h_\nu+v^2\partial_2-h_\nu-2k\partial_2$ for $k\geq0$.
Since $S_0=I$, taking the coefficient of $z^m$ in
\eqref{eq:target-ancestor-inputs} gives, for $m\geq0$,
\[
\begin{aligned}
a_\nu^{\alpha,m}\partial_\alpha
&=\delta_{m1}\partial_1+\delta_{m0}c_\nu h_\nu
 +\sum_{j\geq0}S_j^{X^\nu}(u)
   \bigl((v^2-2m-2j)\partial_2-h_\nu\bigr)\\
&=\delta_{m1}\partial_1-(1-\delta_{m0})c_\nu h_\nu
 -2m\sum_{j\geq0}S_j^{X^\nu}(u)\partial_2,
\end{aligned}
\]
where the second equality follows from \eqref{eq:calibration-center-sums}.
Summing over $m$ proves \eqref{eq:ancestor-backgrounds}.
\end{proof}

\subsection{Shifted $y$-functions and transformed backgrounds}
\label{subsec:background-curves}

At each fixed target primary point $u$,
Proposition~\ref{prop:y-variation} associates an ancestor background to a choice of $y$ through
\eqref{eq:y-background}.
We show that the shifted functions
$y_{\mathrm{sh}}^{\mathrm{di}}$ and
$y_{\mathrm{sh}}^{\mathrm{ad}}$ introduced in
Subsection~\ref{subsec:common-kernel} give precisely
$a_{\mathrm{di}}$ and $a_{\mathrm{ad}}$, respectively.

\begin{proposition}
\label{prop:effective-background}
For corresponding primary points $v$ and $u$ with $(v^1)^2\ne4e^{v^2}$,
the ancestor backgrounds in
\eqref{eq:ancestor-backgrounds} are represented, through
Proposition~\ref{prop:y-variation}, by the curves
$(\log \mathsf{x},y_{\mathrm{sh}}^\nu,B)$ with
\begin{equation}
\mathsf{x}=v^1+\ell+\frac{e^{v^2}}\ell,\qquad
y_{\mathrm{sh}}^{\mathrm{di}}=\mathsf{x}\log\ell,\qquad
y_{\mathrm{sh}}^{\mathrm{ad}}=-\mathsf{x}\log\ell.
\label{eq:actual-effective-y}
\end{equation}
Here $\nu=\mathrm{di},\mathrm{ad}$. Both curves have Bergman kernel
\[
B(\ell_1,\ell_2)
=\frac{d\ell_1\,d\ell_2}{(\ell_1-\ell_2)^2}.
\]
\end{proposition}

\begin{proof}
By \eqref{eq:y-background}, it suffices to verify
\begin{equation}
\mathbf y_{\mathrm{sh}}^\nu(z)
=
R(z)^{-1}\Phi^{-1}
\left[
z\begin{pmatrix}1\\0\end{pmatrix}
-
\begin{pmatrix}
a_\nu^1(z)\\a_\nu^2(z)
\end{pmatrix}
\right],
\qquad \nu=\mathrm{di},\mathrm{ad}.
\label{eq:effective-odd-series}
\end{equation}
All Frobenius data are evaluated at $u$, and $r=e^{u^2}$.

By \eqref{eq:effective-y-periods},
\[
\mathbf y_{\mathrm{sh}}^\nu(z)
=
z\bm{\mathcal G}_{-z}(dy_{\mathrm{sh}}^\nu).
\]
We compute these transforms using
\eqref{eq:negative-Gaussian-integration-by-parts}
and \eqref{eq:flat-negative-periods}.

In the diagonal case, $\ell=e^{u^1}p$ and
\eqref{eq:di-fixed-primary} give
\begin{equation}
\mathsf{x}=e^{u^1}\frac{(1+p)(p+r)}{p}=e^{x^{X^{\mathrm{di}}}}.
\label{eq:di-common-projection}
\end{equation}
Thus \(x^{X^{\mathrm{di}}}=\log \mathsf{x}\), and the coordinate
change preserves the Bergman kernel.
The original differential and the primary derivatives satisfy
\[
dy^\circ=\frac{dp}{p},
\qquad
\partial_1x=1,
\qquad
\partial_2x=\frac{r}{p+r}.
\]
In particular,
\[
e^{u^1}\bigl(r\partial_1x+(1-r)\partial_2x\bigr)
=e^{u^1}\frac{r(p+1)}{p+r}.
\]
A direct calculation gives
\begin{equation}
\begin{aligned}
&d(\mathsf{x}\log p)
-\frac{v^1dy^\circ}{1-z}
-\frac{2e^{u^1}r(p+1)}
{(p+r)(1-z)^2}\,dy^\circ
\\
&\quad=
\left(d-z^{-1}dx\right)
\left[
-\frac{z\mathsf{x}\log p}{1-z}
-\frac{ze^{u^1}(p-r)(p+1)}
{p(1-z)^2}
\right].
\end{aligned}
\label{eq:di-background-certificate}
\end{equation}
The transform \(\bm{\mathcal G}_{-z}\) annihilates the
right-hand side by
\eqref{eq:negative-Gaussian-integration-by-parts}.
Applying \eqref{eq:flat-negative-periods}
\begin{equation*}
\bm{\mathcal G}_{-z}\!\left(
(h^\alpha\partial_\alpha x)\,dy^\circ
\right)
=R(z)^{-1}\Phi^{-1}h
\end{equation*} 
to the remaining terms therefore yields
\begin{equation}
\begin{aligned}
\bm{\mathcal G}_{-z}\bigl(d(\mathsf{x}\log p)\bigr)
&=
R(z)^{-1}\Phi^{-1}
\left[
\frac{v^1}{1-z}\begin{pmatrix}1\\0\end{pmatrix}
+
\frac{2e^{u^1}}{(1-z)^2}\begin{pmatrix}r\\1-r\end{pmatrix}
\right]
\\
&=
R(z)^{-1}\Phi^{-1}
\left[
\begin{pmatrix}1\\0\end{pmatrix}
-\frac1z
\begin{pmatrix}
a_{\mathrm{di}}^1(z)\\a_{\mathrm{di}}^2(z)
\end{pmatrix}
\right],
\end{aligned}
\label{eq:di-effective-period}
\end{equation}
where the last equality follows from
\eqref{eq:ancestor-backgrounds}.
Since \(\ell=e^{u^1}p\), the difference
\(\mathsf{x}\log\ell-\mathsf{x}\log p\) is a constant multiple of \(\mathsf{x}=e^x\).
Consequently, \(\mathsf{x}\log\ell\) and \(\mathsf{x}\log p\) determine the same ancestor background.
Multiplying \eqref{eq:di-effective-period} by \(z\) proves
\eqref{eq:effective-odd-series} in the diagonal case.

In the anti-diagonal case, $\ell=\mathrm i e^{u^1}(p+r)$ and
\eqref{eq:ad-fixed-primary} give
\begin{equation}
\mathsf{x}=\mathrm i e^{u^1}\frac{p(1+p)}{p+r}
=\mathrm i\,e^{x^{X^{\mathrm{ad}}}}.
\label{eq:ad-common-projection}
\end{equation}
Hence \(d\log \mathsf{x}=dx^{X^{\mathrm{ad}}}\).
The additive constant relating the two $x$-functions does not change the local Airy coordinates or the normalized Laplace transform, and the affine coordinate change preserves the Bergman kernel.

In the original coordinate \(p\), we have
\[
dy^\circ=-\mathrm i\,\frac{dp}{p},
\qquad
dx=
\left(
\frac1p+\frac1{p+1}-\frac1{p+r}
\right)dp,
\qquad
\partial_2x=-\frac{r}{p+r}.
\]
It follows that
\[
\begin{aligned}
dy^\circ+\mathrm i\,dx
&=
-\mathrm i
\left(
\frac{dp}{p+r}-\frac{dp}{p+1}
\right),
\\
(\partial_2x)\,dy^\circ
&=
-\mathrm i
\left(
\frac{dp}{p+r}-\frac{dp}{p}
\right).
\end{aligned}
\]
Since \(\bm{\mathcal G}_{-z}(dx)=0\) by
Lemma~\ref{lem:Laplace-properties}, the period identity
\eqref{eq:flat-negative-periods} gives
\begin{equation}
\begin{aligned}
\bm{\mathcal G}_{-z}
\bigl(dy^\circ+\mathrm i\,dx\bigr)
&=
R(z)^{-1}\Phi^{-1}
\begin{pmatrix}1\\0\end{pmatrix},
\\
\bm{\mathcal G}_{-z}
\bigl((\partial_2x)\,dy^\circ\bigr)
&=
R(z)^{-1}\Phi^{-1}
\begin{pmatrix}0\\1\end{pmatrix}.
\end{aligned}
\label{eq:ad-insertion-periods}
\end{equation}

Substituting the anti-diagonal background from
\eqref{eq:ancestor-backgrounds}, we obtain the differential
identity
\begin{equation}
\begin{aligned}
&-\frac{\mathsf{x}}{1-z}\frac{dp}{p+r}
-\left(1-\frac{a_{\mathrm{ad}}^1(z)}z\right)
\bigl(dy^\circ+\mathrm i\,dx\bigr)
+\frac{a_{\mathrm{ad}}^2(z)}z
(\partial_2x)\,dy^\circ
\\
&\quad=
\left(d-z^{-1}dx\right)
\left[
\frac{\mathrm i ze^{u^1}}{(1-z)^2}
\left(
p+r+\frac{r(1-r)}{p+r}
+(1-2r)(2-z)
\right)
\right].
\end{aligned}
\label{eq:ad-background-certificate}
\end{equation}
Moreover, since
\(d\log\ell=dp/(p+r)\) and \(d\mathsf{x}=\mathsf{x}\,dx\),
\begin{equation}
d(-\mathsf{x}\log\ell)
+\frac{\mathsf{x}}{1-z}\frac{dp}{p+r}
=
\left(d-z^{-1}dx\right)
\left[
\frac{z\mathsf{x}\log\ell}{1-z}
\right].
\label{eq:ad-effective-logarithm}
\end{equation}
Adding these two identities and applying
\(\bm{\mathcal G}_{-z}\), the right-hand sides vanish by
\eqref{eq:negative-Gaussian-integration-by-parts}.
Using \eqref{eq:ad-insertion-periods}, we conclude that
\begin{equation}
\bm{\mathcal G}_{-z}\bigl(d(-\mathsf{x}\log\ell)\bigr)
=
R(z)^{-1}\Phi^{-1}
\left[
\begin{pmatrix}1\\0\end{pmatrix}
-\frac1z
\begin{pmatrix}
a_{\mathrm{ad}}^1(z)\\a_{\mathrm{ad}}^2(z)
\end{pmatrix}
\right].
\label{eq:ad-effective-period}
\end{equation}
Multiplying by \(z\) and using
\eqref{eq:effective-y-periods} proves
\eqref{eq:effective-odd-series} in the anti-diagonal case.
The proposition now follows from
Proposition~\ref{prop:y-variation}.
\end{proof}

To translate the comparison of their recursion differentials into identities between descendant correlators, we next compare one-forms with respect to the descendant variables.

\subsection{Descendant extractions}
\label{subsec:insertions}

On the common sphere, write
\[
\mathsf{x}(\ell)=v^1+\ell+\frac{e^{v^2}}{\ell},
\qquad
\ell_\pm=\pm e^{v^2/2}.
\]
The ancestor one-forms $W^{\mathbb P^1,i}_k$ of the $\mathbb P^1$ curve defined in \eqref{eq:ancestor-form-recursions} satisfy
\begin{equation}
W_0^{\mathbb P^1,i}(\ell)
=\sqrt{\frac{\ell_i}{2}}\,
\frac{d\ell}{(\ell-\ell_i)^2},
\qquad
W_{j+1}^{\mathbb P^1,i}
=-d\!\left(\frac{W_j^{\mathbb P^1,i}}{d\mathsf{x}}\right).
\label{eq:source-W}
\end{equation}

Every finite collection of these forms is linearly independent.
Indeed, $W_j^{\mathbb P^1,i}$ has a pole of order $2j+2$ at $\ell_i$, whereas its pole at the other critical point has order at most $2j$; see \cite{zhouTopologicalRecursionVirasoro2024,xiongEYNARDORANTINTOPOLOGICALRECURSION} for more details about this. Comparing the highest pole orders successively proves the assertion.

For a descendant label $(\alpha,k)$, define a linear functional on their linear span space by
\begin{equation}
\lambda_{\alpha,k}\bigl(W_j^{\mathbb P^1,i}\bigr)
=\bigl((\Phi^{\mathbb P^1}(v))^{-1}\bigr)^i{}_\beta
\bigl(S_{k-j}^{\mathbb P^1}(v)\bigr)^\beta{}_\alpha,
\qquad
S_m^{\mathbb P^1}=0\quad(m<0).
\label{eq:source-functional}
\end{equation}
This is the functional associated with a descendant-ancestor transformation; see also \cite{borotLoopEquationsGromovWitten2019}.
To see this, differentiation of the ancestor transformation
\eqref{eq:source-ancestor-input} in the direction $s^{\alpha,k}$ gives
the flat ancestor coefficient
\[
\frac{\partial a^{\beta,j}}{\partial s^{\alpha,k}}
=\bigl(S_{k-j}^{\mathbb P^1}(v)\bigr)^\beta{}_\alpha.
\]
Multiplication by $(\Phi^{\mathbb P^1})^{-1}$ expresses this
coefficient in the normalized canonical basis, giving
\eqref{eq:source-functional}.

For a meromorphic one-form $\omega$ on the common sphere, define
\begin{equation}
L\omega=-d\!\left(\mathsf{x}\frac{\omega}{d\mathsf{x}}\right).
\label{eq:logarithmic-operator}
\end{equation}
This is the recursion operator for ancestor forms when the $x$
is $\log \mathsf{x}$.

Let $W_m^{+,i}$ and $W_m^{-,i}$ denote the ancestor one-forms of the shifted diagonal and anti-diagonal curves, respectively.
The local Airy coordinates associated with $\mathsf{x}$ and $\log \mathsf{x}$ satisfy
\[
\zeta_{\mathbb P^1}
=\sqrt{\mathsf{x}(\ell_i)}\,\zeta_X+O(\zeta_X^3).
\]
Consequently, the initial forms and their differential recursions give
\begin{equation}
\begin{aligned}
W_m^{+,i}
&=\sqrt{\mathsf{x}(\ell_i)}\,L^mW_0^{\mathbb P^1,i},
\\
W_m^{-,i}
&=\sqrt{\mathsf{x}(\ell_i)}\,L^mW_0^{\mathbb P^1,i},
\end{aligned}
\label{eq:target-W}
\end{equation}
with the primary points ($v$ and $u$) corresponding to the relevant transformations \eqref{eq:di-coordinate-transformation} and \eqref{eq:ad-coordinate-transformation}.

Fix $v$ with $(v^1)^2\ne4e^{v^2}$ and, for each comparison,
choose the corresponding point $u$ on the chosen local branch.
We hold $v$ and $u$ fixed when differentiating with respect to
the descendant variables $s^{\alpha,k}$.

The differentials $\omega^+$ and $\omega^-$ use the shifted functions $y_{\mathrm{sh}}^{\mathrm{di}}$ and $y_{\mathrm{sh}}^{\mathrm{ad}}$, respectively. 
By Proposition~\ref{prop:y-variation} and the identification of the backgrounds in Proposition~\ref{prop:effective-background}, for $n\geq1$ and $2g-2+n>0$, we have
\begin{equation}
\begin{aligned}
\omega^+_{g,n}(\ell_1,\ldots,\ell_n)
={}&\sum_{m_1,\ldots,m_n\geq0}
\left.
\frac{\partial^n\overline F_g^{X^{\mathrm{di}},\mathrm{tw}}(\bar t;u)}
{\partial\bar t^{\alpha_1,m_1}\cdots
 \partial\bar t^{\alpha_n,m_n}}
\right|_{\bar t=a_{\mathrm{di}}}
\prod_{r=1}^n
(\Phi^{X^{\mathrm{di}}}(u))^{\alpha_r}{}_{i_r}
\,W_{m_r}^{+,i_r}(\ell_r),
\\[6pt]
\omega^-_{g,n}(\ell_1,\ldots,\ell_n)
={}&\sum_{m_1,\ldots,m_n\geq0}
\left.
\frac{\partial^n\overline F_g^{X^{\mathrm{ad}}}(\bar t;u)}
{\partial\bar t^{\alpha_1,m_1}\cdots
 \partial\bar t^{\alpha_n,m_n}}
\right|_{\bar t=a_{\mathrm{ad}}}
\prod_{r=1}^n
(\Phi^{X^{\mathrm{ad}}}(u))^{\alpha_r}{}_{i_r}
\,W_{m_r}^{-,i_r}(\ell_r).
\end{aligned}
\label{eq:modified-conifold-ancestor-TR}
\end{equation}
Here $a_{\mathrm{di}}$ and $a_{\mathrm{ad}}$ are given by \eqref{eq:ancestor-backgrounds}.



The ancestor backgrounds in these formulas have already been determined. To compare the descendant correspondence, we must identify the action
of $\lambda_{\alpha,k}$ on $W_m^{+,i}$ and $W_m^{-,i}$.
By \eqref{eq:target-W}, this reduces to computing the dual action
of $L$.

\begin{lemma}
\label{lem:logarithmic-action}
The functionals in \eqref{eq:source-functional} satisfy
\begin{equation}
\lambda_{2,k}L=(k+1)\lambda_{2,k},
\qquad
\lambda_{1,k}L=k\lambda_{1,k}+2\lambda_{2,k-1},
\label{eq:Jordan-action}
\end{equation}
where $k\geq0$ and $\lambda_{2,-1}=0$.
\end{lemma}

\begin{proof}
In this proof, put
\[
\kappa_i=\sqrt{\frac{\ell_i}{2}},
\qquad \mathsf{x}_i=\mathsf{x}(\ell_i).
\]
We first compute the action of $L$ on the $\mathbb P^1$ ancestor one-forms:
\begin{equation}
LW_j^{\mathbb P^1,i}
=\mathsf{x}_iW_{j+1}^{\mathbb P^1,i}
+\left(j+\frac12\right)W_j^{\mathbb P^1,i}
+\frac{\kappa_i}{2\kappa_{-i}}W_j^{\mathbb P^1,-i}.
\label{eq:L-on-W}
\end{equation}
Here $-i$ means another critical point, i.e., for $i=\pm$, $-i=\mp$.

For $j=0$, this follows from
\[
\frac{(\mathsf{x}-\mathsf{x}_i)W_0^{\mathbb P^1,i}}{d\mathsf{x}}
=\frac{\kappa_i\ell}{\ell^2-\ell_i^2}
=\frac{\kappa_i}{2}
\left(\frac1{\ell-\ell_i}+\frac1{\ell+\ell_i}\right)
\]
by applying $-\d$. 
If we write
\[
D\omega:=-d\!\left(\frac{\omega}{d\mathsf{x}}\right),
\]
then $[L,D]=D$. Applying this commutator repeatedly proves
\eqref{eq:L-on-W} for every $j$.

On the other hand, consider the $\mathbb{P}^1$ Euler field
\[
E=v^1\partial_{v^1}+2\partial_{v^2}.
\]
The explicit calibration series in Section~\ref{subsec:S-calibration}
satisfy
\begin{equation}
\begin{aligned}
E\theta_{1,m}^{\mathbb P^1}
&=m\theta_{1,m}^{\mathbb P^1}
+2\theta_{2,m-1}^{\mathbb P^1}+2\delta_{m0},
\\
E\theta_{2,m}^{\mathbb P^1}
&=(m+1)\theta_{2,m}^{\mathbb P^1},
\end{aligned}
\qquad
\theta_{2,-1}^{\mathbb P^1}=0.
\label{eq:source-Euler}
\end{equation}
Indeed, one can obtain from \eqref{eq:P1-calibration} the coefficients
\begin{equation}
\theta^{\mathbb{P}^1}_{1,m}
=
\sum_{d=0}^{\lfloor m/2\rfloor}
(v^2-2\sum_{i=1}^{d}\frac{1}{i})
\frac{e^{dv^2}(v^1)^{m-2d}}{(d!)^2(m-2d)!},
\label{coeff-p1-1}
\end{equation} 
and
\begin{equation}
\theta_{2,m}^{\mathbb{P}^1}
=
\sum_{d=0}^{\lfloor(m+1)/2\rfloor}
\frac{e^{dv^2}(v^1)^{m+1-2d}}{(d!)^2(m+1-2d)!}.
\end{equation}
The normalized canonical vector field \eqref{eq:normalized-idempotent-matrix} is
\[
f_i:=\Phi^\gamma{}_i\partial_\gamma=\kappa_i
\left(\partial_{v^1}+\ell_i^{-1}\partial_{v^2}\right),
\]
and $[E,f_i]=-\tfrac12f_i$. Write
\[
A_{\alpha,m}^i
=\bigl((\Phi^{\mathbb P^1})^{-1}\bigr)^i{}_\beta
\bigl(S_m^{\mathbb P^1}\bigr)^\beta{}_\alpha
=f_i\cdot\theta_{\alpha,m}^{\mathbb P^1}.
\]
Equation \eqref{eq:source-Euler} then gives
\begin{equation}
\begin{aligned}
EA_{1,m}^i
&=\left(m-\frac12\right)A_{1,m}^i+2A_{2,m-1}^i,
\\
EA_{2,m}^i
&=\left(m+\frac12\right)A_{2,m}^i.
\end{aligned}
\label{eq:canonical-calibration-homogeneity}
\end{equation}

Meanwhile, direct differentiation of the inverse canonical transform gives
\[
E\bigl((\Phi^{\mathbb P^1})^{-1}\bigr)^i{}_\beta
=E(\Phi^{\gamma}{}_{i}\eta_{\gamma\beta})
=\frac{\kappa_i}{2\kappa_{-i}}
\bigl((\Phi^{\mathbb P^1})^{-1}\bigr)^{-i}{}_\beta.
\]
And the quantum differential equation implies
\[
ES_m^{\mathbb P^1}=(E\star)S_{m-1}^{\mathbb P^1}.
\]
Quantum multiplication by $E$ has eigenvalue $\mathsf{x}_i$ in the $i$-th canonical direction. Indeed, in the flat basis, the quantum product is
\[
E\star=
\begin{pmatrix}
v^1&2e^{v^2}\\
2&v^1
\end{pmatrix}.
\]
Using $e^{v^2}=\ell_i^2$, one has
\[
\begin{aligned}
\left(\frac{\kappa_i}{\ell_i},\,\kappa_i\right)(E\star)
&=\left(
\frac{\kappa_i}{\ell_i}(v^1+2\ell_i),\,
\kappa_i(v^1+2\ell_i)
\right)\\
&=\mathsf{x}_i\left(\frac{\kappa_i}{\ell_i},\,\kappa_i\right).
\end{aligned}
\]

Therefore,
\begin{equation}
EA_{\alpha,m}^i
=\mathsf{x}_iA_{\alpha,m-1}^i
+\frac{\kappa_i}{2\kappa_{-i}}A_{\alpha,m}^{-i}.
\label{eq:canonical-calibration-connection}
\end{equation}

By \eqref{eq:source-functional} and \eqref{eq:L-on-W},
\[
\begin{aligned}
\lambda_{\alpha,k}\bigl(LW_j^{\mathbb P^1,i}\bigr)
&=\mathsf{x}_iA_{\alpha,k-j-1}^i
+\left(j+\frac12\right)A_{\alpha,k-j}^i
+\frac{\kappa_i}{2\kappa_{-i}}A_{\alpha,k-j}^{-i}
\\
&=EA_{\alpha,k-j}^i
+\left(j+\frac12\right)A_{\alpha,k-j}^i.
\end{aligned}
\]
Substituting \eqref{eq:canonical-calibration-homogeneity} proves
both identities in \eqref{eq:Jordan-action}.
\end{proof}

Iteration of the lemma yields
\begin{equation}
\begin{aligned}
\lambda_{2,k}L^m
&=(k+1)^m\lambda_{2,k},
&&k,m\geq0,
\\
\lambda_{1,k}L^m
&=k^m\lambda_{1,k}
+2mk^{m-1}\lambda_{2,k-1},
&&k,m\geq1,
\\
\lambda_{1,0}L^m&=0,
&&m\geq1.
\end{aligned}
\label{eq:Jordan-powers}
\end{equation}
These are precisely the powers that occur in the LYZZ coefficients
\eqref{eq:di-b-coefficients} and \eqref{eq:ad-b-coefficients}.

We now compare the resolved conifold ancestor invariants with $\mathbb{P}^1$ descendant invariants explicitly. Keep the primary point fixed and substitute $t=t_{\mathrm{di}}(s)$
or $t=t_{\mathrm{ad}}(s)$ into \eqref{eq:target-ancestor-inputs}.
At fixed $u$, the target ancestor background is
\[
a_\nu(z;s)
=\bigl[S^{X^\nu}(u;z)(t_\nu(s;z)-z\partial_1)\bigr]_+
+z\partial_1.
\]
At $s=s_v$ this gives the background in \eqref{eq:ancestor-backgrounds}. Differentiating at fixed $u$ and using \eqref{eq:descendant-transformations} gives
\begin{equation}
\begin{aligned}
\frac{\partial a_{\mathrm{di}}^{\beta,m}}
{\partial s^{\alpha,k}}
&=\sum_{j\geq0}
\bigl(S_j^{X^{\mathrm{di}}}\bigr)^\beta{}_\gamma
b_{\alpha,k}^{\gamma,m+j,\mathrm{di}},
\\
\frac{\partial a_{\mathrm{ad}}^{\beta,m}}
{\partial s^{\alpha,k}}
&=\sum_{j\geq0}
\bigl(S_j^{X^{\mathrm{ad}}}\bigr)^\beta{}_\gamma
b_{\alpha,k}^{\gamma,m+j,\mathrm{ad}}.
\end{aligned}
\label{eq:target-ancestor-variations}
\end{equation}

Thus one has the proposition below.
\begin{proposition}
\label{prop:external-identity}
For every $\alpha\in\{1,2\}$ and $k,m\geq0$,
\begin{equation}
\begin{aligned}
\lambda_{\alpha,k}\bigl(W_m^{+,i}\bigr)
&=\bigl((\Phi^{X^{\mathrm{di}}})^{-1}\bigr)^i{}_\beta
\frac{\partial a_{\mathrm{di}}^{\beta,m}}
{\partial s^{\alpha,k}},
\\
\lambda_{\alpha,k}\bigl(W_m^{-,i}\bigr)
&=-\bigl((\Phi^{X^{\mathrm{ad}}})^{-1}\bigr)^i{}_\beta
\frac{\partial a_{\mathrm{ad}}^{\beta,m}}
{\partial s^{\alpha,k}}.
\end{aligned}
\label{eq:external-identity}
\end{equation}
\end{proposition}

\begin{proof}
We first prove the identities for $m=0$. In each theory,
\[
(\Phi^{-1})^i{}_\beta(S_j)^\beta{}_\gamma
=
\Phi^\delta{}_i\,\partial_\delta\theta_{\gamma,j}.
\]
Thus we compare the differentiated calibration identities
in the normalized canonical basis.

In the diagonal case, the coordinate transformation gives
\[
\ell_i=e^{u^1}p_i,
\qquad
p_i^2=e^{u^2},
\qquad
\mathsf{x}(\ell_i)=e^{u^1}(1+p_i)^2.
\]
The corresponding rows of the transition matrix are
\[
\bigl((\Phi^{X^{\mathrm{di}}})^T\bigr)^{i}{}_{-}
=
\frac{(p_i,1-p_i)}{\sqrt{2p_i}},
\qquad
\bigl((\Phi^{\mathbb P^1})^T\bigr)^{i}{}_{-}
=
\frac{(\ell_i,1)}{\sqrt{2\ell_i}}.
\]
The transpose of the Jacobian of
\eqref{eq:di-coordinate-transformation} is
\[
\frac{\partial(v^1,v^2)}{\partial(u^1,u^2)}=
\begin{pmatrix}
e^{u^1}(1+e^{u^2})&2\\
e^{u^1+u^2}&1
\end{pmatrix}.
\]
Consequently,
\begin{equation*}
\begin{aligned}
&\frac{(p_i,1-p_i)}{\sqrt{2p_i}}
\begin{pmatrix}
e^{u^1}(1+e^{u^2})&2\\
e^{u^1+u^2}&1
\end{pmatrix}
\\
&\qquad=
\frac{1+p_i}{\sqrt{2p_i}}
(e^{u^1}p_i,1)
=
\sqrt{\mathsf{x}(\ell_i)}
\frac{(\ell_i,1)}{\sqrt{2\ell_i}}.
\end{aligned}
\end{equation*}
That is, 
\begin{equation}
(\Phi^{X^{\mathrm{di}}})^\delta{}_i
\frac{\partial v^\rho}{\partial u^\delta}
=
\sqrt{\mathsf{x}(\ell_i)}\,
(\Phi^{\mathbb P^1})^\rho{}_i.
\label{eq:di-initial-leg}
\end{equation}

Differentiating \eqref{eq:di-calibration-transformation}
with respect to $u^\delta$ gives
\begin{equation}
\frac{\partial v^\rho}{\partial u^\delta}
\frac{\partial\theta_{\alpha,k}^{\mathbb P^1}}{\partial v^\rho}
=
\sum_{j\geq0}
b_{\alpha,k}^{\gamma,j,\mathrm{di}}
\frac{\partial\theta_{\gamma,j}^{X^{\mathrm{di}}}}
{\partial u^\delta}.
\end{equation}
Using \eqref{eq:source-functional} and
\eqref{eq:scalar-matrix-calibration}, we therefore obtain
\begin{equation}
\begin{aligned}
&\sqrt{\mathsf{x}(\ell_i)}\,
\lambda_{\alpha,k}
\bigl(W_0^{\mathbb P^1,i}\bigr)
=
\sqrt{\mathsf{x}(\ell_i)}\,
(\Phi^{\mathbb P^1})^\rho{}_i
\frac{\partial\theta_{\alpha,k}^{\mathbb P^1}}
{\partial v^\rho}
=
(\Phi^{X^{\mathrm{di}}})^\delta{}_i
\frac{\partial v^\rho}{\partial u^\delta}
\frac{\partial\theta_{\alpha,k}^{\mathbb P^1}}
{\partial v^\rho}
\\
&\quad=
(\Phi^{X^{\mathrm{di}}})^\delta{}_i
\sum_{j\geq0}
b_{\alpha,k}^{\gamma,j,\mathrm{di}}
\frac{\partial\theta_{\gamma,j}^{X^{\mathrm{di}}}}
{\partial u^\delta}
=
\bigl((\Phi^{X^{\mathrm{di}}})^{-1}\bigr)^i{}_\beta
\sum_{j\geq0}
\bigl(S_j^{X^{\mathrm{di}}}\bigr)^\beta{}_\gamma
b_{\alpha,k}^{\gamma,j,\mathrm{di}}.
\end{aligned}
\end{equation}
By \eqref{eq:target-W} and
\eqref{eq:target-ancestor-variations}, this proves the
diagonal identity in \eqref{eq:external-identity} for $m=0$.

In the anti-diagonal case, the equation of critical point and
the coordinate transformation give
\[
p_i^2+2e^{u^2}p_i+e^{u^2}=0,
\qquad
\ell_i=\mathrm i\,e^{u^1}(p_i+e^{u^2}),
\]
and
\[
\mathsf{x}(\ell_i)
=
\frac{\mathrm i\,e^{u^1}(1+p_i)^2}{1-e^{u^2}}.
\]
The transition matrix satisfies
\[
\bigl((\Phi^{X^{\mathrm{ad}}})^T\bigr)^{i}{}_{-}
=
-\frac{
\left(\dfrac{p_i+e^{u^2}}{1-e^{u^2}},\,1\right)
}{
\sqrt{-\dfrac{2(p_i+e^{u^2})}{1-e^{u^2}}}
}.
\]
With the square-root choices fixed above, one also has
\[
\frac{\sqrt{\mathsf{x}(\ell_i)}}{\sqrt{2\ell_i}}
=
\frac{\mathrm i(1+p_i)}
{
(1-e^{u^2})
\sqrt{-\dfrac{2(p_i+e^{u^2})}{1-e^{u^2}}}
}.
\]
These formulas imply
\begin{equation}
\sqrt{\mathsf{x}(\ell_i)}
\frac{(\ell_i,1)}{\sqrt{2\ell_i}}
\begin{pmatrix}
e^{-u^1}&-e^{-u^1}\\
\mathrm i\,e^{u^2}&\mathrm i(1-e^{u^2})
\end{pmatrix}
=
\bigl((\Phi^{X^{\mathrm{ad}}})^T\bigr)^{i}{}_{-}.
\end{equation}
Indeed, multiplication of the row vector by the displayed
matrix gives
\[
(\ell_i,1)
\begin{pmatrix}
e^{-u^1}&-e^{-u^1}\\
\mathrm i\,e^{u^2}&\mathrm i(1-e^{u^2})
\end{pmatrix}
=
\mathrm i\,
(p_i+2e^{u^2},\,1-p_i-2e^{u^2}),
\]
and the required simplification follows from
$p_i^2+2e^{u^2}p_i+e^{u^2}=0$.
Thus,
\begin{equation}
\sqrt{\mathsf{x}(\ell_i)}\,
\bigl((\Phi^{\mathbb P^1})^T\bigr)^{i}{}_{-}
\begin{pmatrix}
e^{-u^1}&-e^{-u^1}\\
\mathrm i e^{u^2}&\mathrm i(1-e^{u^2})
\end{pmatrix}
=
\bigl((\Phi^{X^{\mathrm{ad}}})^T\bigr)^{i}{}_{-}.
\label{eq:ad-initial-leg}
\end{equation}
It is consistent with the transformation between $\theta$ functions \eqref{eq:ad-calibration-transformation}:
\begin{equation*}
\begin{pmatrix}
\partial_{v^1}\theta_{\alpha,k}^{\mathbb P^1}\\
\partial_{v^2}\theta_{\alpha,k}^{\mathbb P^1}
\end{pmatrix}_{v=v(u)}
=
-
\begin{pmatrix}
e^{-u^1}&-e^{-u^1}\\
\mathrm i e^{u^2}&\mathrm i(1-e^{u^2})
\end{pmatrix}
\sum_{j\ge0}
b_{\alpha,k}^{\gamma,j,\mathrm{ad}}
\begin{pmatrix}
\partial_{u^1}\theta_{\gamma,j}^{X^{\mathrm{ad}}}\\
\partial_{u^2}\theta_{\gamma,j}^{X^{\mathrm{ad}}}
\end{pmatrix}.
\end{equation*}

Applying these transformations, we obtain
\begin{align*}
&\sqrt{\mathsf{x}(\ell_i)}\,
\lambda_{\alpha,k}\bigl(W_0^{\mathbb P^1,i}\bigr)\\
&\quad=
-\sqrt{\mathsf{x}(\ell_i)}\,
\bigl((\Phi^{\mathbb P^1})^T\bigr)^{i}{}_{-}
\begin{pmatrix}
e^{-u^1}&-e^{-u^1}\\
\mathrm i e^{u^2}&\mathrm i(1-e^{u^2})
\end{pmatrix}
\sum_{j\ge0}
b_{\alpha,k}^{\gamma,j,\mathrm{ad}}
\begin{pmatrix}
\partial_{u^1}\theta_{\gamma,j}^{X^{\mathrm{ad}}}\\
\partial_{u^2}\theta_{\gamma,j}^{X^{\mathrm{ad}}}
\end{pmatrix}\\
&\quad=
-(\Phi^{X^{\mathrm{ad}}})^\delta{}_i
\sum_{j\ge0}
b_{\alpha,k}^{\gamma,j,\mathrm{ad}}\,
\partial_{u^\delta}\theta_{\gamma,j}^{X^{\mathrm{ad}}}\\
&\quad=
-\bigl((\Phi^{X^{\mathrm{ad}}})^{-1}\bigr)^i{}_\beta
\sum_{j\ge0}
(S_j^{X^{\mathrm{ad}}})^\beta{}_\gamma
b_{\alpha,k}^{\gamma,j,\mathrm{ad}}.
\end{align*}
The last equality follows again from
\eqref{eq:scalar-matrix-calibration} and
\eqref{eq:canonical-frame}.
Together with \eqref{eq:target-W} and
\eqref{eq:target-ancestor-variations}, this proves
the anti-diagonal identity in \eqref{eq:external-identity}
for \(m=0\).

We next compare the higher ancestor indices.
The explicit coefficients
\eqref{eq:di-b-coefficients} and
\eqref{eq:ad-b-coefficients} satisfy
\begin{equation}
\begin{aligned}
b_{2,k}^{\gamma,m+j,\mathrm{di}}
&=(k+1)^m b_{2,k}^{\gamma,j,\mathrm{di}},
\\
b_{2,k}^{\gamma,m+j,\mathrm{ad}}
&=(k+1)^m b_{2,k}^{\gamma,j,\mathrm{ad}},
\\
b_{1,k}^{\gamma,m+j,\mathrm{di}}
&=k^m b_{1,k}^{\gamma,j,\mathrm{di}}
+2mk^{m-1}b_{2,k-1}^{\gamma,j,\mathrm{di}},
\\
b_{1,k}^{\gamma,m+j,\mathrm{ad}}
&=k^m b_{1,k}^{\gamma,j,\mathrm{ad}}
+2mk^{m-1}b_{2,k-1}^{\gamma,j,\mathrm{ad}}.
\end{aligned}
\label{eq:b-height-recurrences}
\end{equation}
Here $m\geq1$ and $j\geq0$; the first two identities hold
for $k\geq0$, and the last two for $k\geq1$.

On the other hand, \eqref{eq:target-W} and
\eqref{eq:Jordan-powers} give
\[
\begin{aligned}
\lambda_{2,k}\bigl(W_m^{\pm,i}\bigr)
&=
\sqrt{\mathsf{x}(\ell_i)}\,(k+1)^m
\lambda_{2,k}\bigl(W_0^{\mathbb P^1,i}\bigr),
\\
\lambda_{1,k}\bigl(W_m^{\pm,i}\bigr)
&=
\sqrt{\mathsf{x}(\ell_i)}
\left[
k^m\lambda_{1,k}\bigl(W_0^{\mathbb P^1,i}\bigr)
+
2mk^{m-1}
\lambda_{2,k-1}\bigl(W_0^{\mathbb P^1,i}\bigr)
\right],
\qquad k\geq1.
\end{aligned}
\]
Substituting the identities already proved for $m=0$,
then applying \eqref{eq:b-height-recurrences}, yields
\[
\begin{aligned}
\lambda_{\alpha,k}\bigl(W_m^{+,i}\bigr)
&=
\bigl((\Phi^{X^{\mathrm{di}}})^{-1}\bigr)^i{}_\beta
\sum_{j\geq0}
\bigl(S_j^{X^{\mathrm{di}}}\bigr)^\beta{}_\gamma
b_{\alpha,k}^{\gamma,m+j,\mathrm{di}},
\\
\lambda_{\alpha,k}\bigl(W_m^{-,i}\bigr)
&=
-\bigl((\Phi^{X^{\mathrm{ad}}})^{-1}\bigr)^i{}_\beta
\sum_{j\geq0}
\bigl(S_j^{X^{\mathrm{ad}}}\bigr)^\beta{}_\gamma
b_{\alpha,k}^{\gamma,m+j,\mathrm{ad}}.
\end{aligned}
\]
By \eqref{eq:target-ancestor-variations}, these are precisely
the required identities.

Finally, for $(\alpha,k)=(1,0)$ and $m\geq1$,
\eqref{eq:Jordan-powers} gives
$\lambda_{1,0}(W_m^{\pm,i})=0$.
The corresponding right-hand sides also vanish, since
$b_{1,0}^{\gamma,m+j,\mathrm{di}}
=b_{1,0}^{\gamma,m+j,\mathrm{ad}}=0$
for $m+j\geq1$.
This completes the proof.
\end{proof}

The proposition describes the action of those functionals on the common recursion
differentials in different perspectives for $\mathbb{P}^1$ and resolved conifold. 
This gives the comparison needed to pass from the equality of differentials to the descendant correspondences.

\subsection{The stable range correspondence}
\label{subsec:stable-marked}
We now prove the stable correspondence.
We evaluate the derivatives of
\[
F_g^{X^{\mathrm{di}},\mathrm{tw}}\bigl(t_{\mathrm{di}}(s)\bigr),
\qquad
F_g^{X^{\mathrm{ad}}}\bigl(t_{\mathrm{ad}}(s)\bigr)
\]
at the descendant point $s_v^{\alpha,k}=v^\alpha\delta_{k0}$.
In each target theory we use the ancestor--descendant formula at the
corresponding fixed primary point $u$. The ancestor backgrounds are given
by \eqref{eq:target-ancestor-inputs}, and their transformation coefficients are in \eqref{eq:target-ancestor-variations}.

\begin{theorem}
\label{thm:stable-marked}
Let $v$ be any primary point on the chosen regular domain
$(v^1)^2\ne4e^{v^2}$, and set $s_v^{\alpha,k}=v^\alpha\delta_{k0}$.
For every $g\geq0$, $n\geq1$ with $2g-2+n>0$, and every
collection of descendant labels $(\alpha_r,k_r)$, one has
\begin{equation}
\begin{aligned}
&\left.
\frac{\partial^n F_g^{\mathbb P^1}(s)}
{\partial s^{\alpha_1,k_1}\cdots
 \partial s^{\alpha_n,k_n}}
\right|_{s=s_v}=
\left.
\frac{\partial^n}
{\partial s^{\alpha_1,k_1}\cdots
 \partial s^{\alpha_n,k_n}}
F_g^{X^{\mathrm{di}},\mathrm{tw}}
\bigl(t_{\mathrm{di}}(s)\bigr)
\right|_{s=s_v},
\\[6pt]
&\left.
\frac{\partial^n F_g^{\mathbb P^1}(s)}
{\partial s^{\alpha_1,k_1}\cdots
 \partial s^{\alpha_n,k_n}}
\right|_{s=s_v}=
\left.
\frac{\partial^n}
{\partial s^{\alpha_1,k_1}\cdots
 \partial s^{\alpha_n,k_n}}
F_g^{X^{\mathrm{ad}}}
\bigl(t_{\mathrm{ad}}(s)\bigr)
\right|_{s=s_v}.
\end{aligned}
\label{eq:marked-correspondence}
\end{equation}

\end{theorem}

\begin{proof}
Apply the tensor product $\bigotimes_{r=1}^n\lambda_{\alpha_r,k_r}$ to the stable recursion differentials on the common sphere.

For the $\mathbb{P}^1$ curve, the DOSS correspondence and the definition of $\lambda_{\alpha,k}$ give
\begin{equation}
\begin{aligned}
&\left(
\bigotimes_{r=1}^n\lambda_{\alpha_r,k_r}
\right)\omega_{g,n}^{\mathbb P^1}
\\
&\quad=
\sum_{m_1,\ldots,m_n\geq0}
\left\langle
\prod_{r=1}^n\bar\tau_{m_r}(\partial_{\beta_r})
\right\rangle^{\mathbb P^1}_{g,n;v}
\prod_{r=1}^n
\bigl(S_{k_r-m_r}^{\mathbb P^1}(v)\bigr)^{\beta_r}{}_{\alpha_r}
\\
&\quad=
\left.
\frac{\partial^n F_g^{\mathbb P^1}(s)}
{\partial s^{\alpha_1,k_1}\cdots
 \partial s^{\alpha_n,k_n}}
\right|_{s=s_v}.
\end{aligned}
\label{eq:source-marked-contraction}
\end{equation}
The ancestor correlators have zero ancestor background at $v$,
and $S^{\mathbb P^1}$ is evaluated at the same point.
The last equality follows from
\eqref{eq:source-ancestor-descendant-derivatives}
and the definition \eqref{eq:ancestor-invariants}
of the ancestor potential.

For the resolved conifold curves, use the ancestor potential expressions \eqref{eq:modified-conifold-ancestor-TR}. 
The comparison in Proposition~\ref{prop:external-identity} yields
\begin{equation}
\begin{aligned}
&\left(
\bigotimes_{r=1}^n\lambda_{\alpha_r,k_r}
\right)\omega_{g,n}^{+}
=
\left.
\frac{\partial^n}
{\partial s^{\alpha_1,k_1}\cdots
 \partial s^{\alpha_n,k_n}}
F_g^{X^{\mathrm{di}},\mathrm{tw}}
\bigl(t_{\mathrm{di}}(s)\bigr)
\right|_{s=s_v},
\\[6pt]
&\left(
\bigotimes_{r=1}^n\lambda_{\alpha_r,k_r}
\right)\omega_{g,n}^{-}=
(-1)^n
\left.
\frac{\partial^n}
{\partial s^{\alpha_1,k_1}\cdots
 \partial s^{\alpha_n,k_n}}
F_g^{X^{\mathrm{ad}}}
\bigl(t_{\mathrm{ad}}(s)\bigr)
\right|_{s=s_v}.
\end{aligned}
\label{eq:target-marked-contractions}
\end{equation}
Indeed, the ancestor argument is affine in $s$ at fixed $u$.
The chain rule therefore gives one factor
\(
\partial a^{\beta_r,m_r}/\partial s^{\alpha_r,k_r}
\)
for each distinguished differential, exactly as in
\eqref{eq:external-identity}. In the anti-diagonal case, each
identity contributes a negative sign.

All contractions in \eqref{eq:source-marked-contraction} and
\eqref{eq:target-marked-contractions} are finite: the contributing
ancestor indices satisfy
\[
m_1+\cdots+m_n\leq3g-3+n.
\]

Finally, Proposition~\ref{prop:kernel-comparison} gives
\[
\omega_{g,n}^{+}=\omega_{g,n}^{\mathbb P^1},
\qquad
\omega_{g,n}^{-}=(-1)^n\omega_{g,n}^{\mathbb P^1}.
\]
Combining these identities with \eqref{eq:source-marked-contraction}
and \eqref{eq:target-marked-contractions} proves both
correspondences. 
\end{proof}

In positive genus the theorem determines every Taylor coefficient of
positive order at $s_v$. The remaining constants are determined below.

\subsection{Vacuum terms and the full correspondences}
\label{subsec:vacuum-correspondence}

The genus-zero identities were proved by LYZZ
\cite{liuEquivariantGromovWitten2022}.
Since the Hodge twist is trivial in genus zero, these identities also apply to the twisted diagonal theory.
We now determine the positive-genus constants left undetermined
by the marked correspondence.

\begin{theorem}[LYZZ conjecture with Hodge twist]
\label{thm:full-correspondences}
For every $g\geq0$, one has
\begin{equation}
\begin{aligned}
F_g^{\mathbb P^1}(s)
&=
F_g^{X^{\mathrm{di}},\mathrm{tw}}
\bigl(t_{\mathrm{di}}(s)\bigr)
+\delta_{g0}A_{\mathrm{di}}(s),
\\
F_g^{\mathbb P^1}(s)
&=
F_g^{X^{\mathrm{ad}}}
\bigl(t_{\mathrm{ad}}(s)\bigr)
+\delta_{g0}A_{\mathrm{ad}}(s),
\end{aligned}
\label{eq:full-correspondences}
\end{equation}
where
\begin{equation}
\begin{aligned}
A_{\mathrm{di}}(s)
&=
\sum_{k,\ell\geq0}
\frac{s^{1,k}s^{2,\ell}}
{(k+\ell+1)k!\ell!}
-
\sum_{\ell\geq0}
\frac{s^{2,\ell}}{(\ell+2)\ell!},
\\
A_{\mathrm{ad}}(s)
&=
\sum_{k,\ell\geq0}
\frac{\mathrm i^{\,k+\ell+1}s^{1,k}s^{2,\ell}}
{(k+\ell+1)k!\ell!}
+
\sum_{\ell\geq0}
\frac{\mathrm i^{\,\ell}s^{2,\ell}}
{(\ell+2)\ell!}.
\end{aligned}
\label{eq:genus-zero-correction-polynomials}
\end{equation}
\end{theorem}

\begin{proof}

For $g\geq1$, Theorem~\ref{thm:stable-marked} shows that each
difference in \eqref{eq:full-correspondences} is constant. 
We determine these constants separately for $g\geq2$ and $g=1$.

For $g\geq2$, the dilaton equation gives
\begin{equation}
\left(
\frac{\partial}{\partial t^{1,1}}
-
t^{\alpha,k}\frac{\partial}{\partial t^{\alpha,k}}
\right)F_g(t)
=(2g-2)F_g(t)
\label{eq:vacuum-dilaton}
\end{equation}
for each of the three theories. 

Each descendant transformation in
\eqref{eq:descendant-transformations} is linear in the dilaton-shifted variables. 
Applying \eqref{eq:vacuum-dilaton} to the two constant differences
therefore gives
\[
\begin{aligned}
(2g-2)\left[
F_g^{\mathbb P^1}(s)
-
F_g^{X^{\mathrm{di}},\mathrm{tw}}
\bigl(t_{\mathrm{di}}(s)\bigr)
\right]&=0,
\\
(2g-2)\left[
F_g^{\mathbb P^1}(s)
-
F_g^{X^{\mathrm{ad}}}
\bigl(t_{\mathrm{ad}}(s)\bigr)
\right]&=0.
\end{aligned}
\]
This proves both identities for $g\geq2$.

In genus one, write $F_1(v)$ or $F_1(u)$ for the primary restriction of the corresponding potential. Givental's genus-one formula is \cite{giventalSemisimpleFrobeniusStructures2001}
\begin{equation}
dF_1
=
\frac12\operatorname{Tr}(R_1\,dU)
+
\frac1{48}d\log(\Delta_+\Delta_-).
\label{eq:genus-one-primary-differential}
\end{equation}
Substituting the data of Sections~\ref{subsec:Frobenius} and \ref{subsec:R-comparison} gives
\begin{equation}
\begin{aligned}
F_1^{\mathbb P^1}(v)
&=-\frac{v^2}{24},
\\
F_1^{X^{\mathrm{ad}}}(u)
&=-\frac{u^2}{24}
-\frac1{12}\log(1-e^{u^2}),
\\
F_1^{X^{\mathrm{di}},\mathrm{tw}}(u)
&=-\frac{u^2}{24}
+\frac1{12}\log(1-e^{u^2}).
\end{aligned}
\label{eq:genus-one-primary-potentials}
\end{equation}
The integration constants are fixed by the Gromov--Witten normalization: after subtracting the degree-zero contribution, the primary potentials have zero constant term in their large-radius expansions.

It remains to include the ancestor backgrounds. At a fixed
semisimple primary point $v$, let $a^{\alpha,0}=0$. A genus-one ancestor vacuum
term with $b\geq1$ insertions is integrated over
$\overline{\mathcal M}_{1,b}$, which has complex dimension $b$.
Since every background insertion has positive ancestor degree,
all insertions must have ancestor degree one, and only the
degree-zero part of the cohomological field theory contributes.
Using
\[
\int_{\overline{\mathcal M}_{1,b}}
\bar\psi_1\cdots\bar\psi_b
=\frac{(b-1)!}{24}
\]
and the genus-one topological field theory in the normalized
canonical basis, we obtain
\begin{equation}
\begin{aligned}
\overline F_1(a;v)
&=
\frac1{24}
\sum_{i\in\{+,-\}}
\sum_{b\geq1}\frac1b
\left(
\frac{(\Phi^{-1})^i{}_\alpha a^{\alpha,1}}
{\Delta_i^{-1/2}}
\right)^b
\\
&=
-\frac1{24}
\sum_{i\in\{+,-\}}
\log\left(
\frac{y_1^{\,i}}{y_1^{\circ,i}}
\right),
\end{aligned}
\label{eq:genus-one-ancestor-vacuum}
\end{equation}
where
$
y_1^{\,i}
=
y_1^{\circ,i}
-
(\Phi^{-1})^i{}_\alpha a^{\alpha,1}$ and $y_1^{\circ,i}=\Delta_i^{-1/2}$. 
At the corresponding target point $u$, the genus-one
ancestor--descendant formula \cite{giventalGromovWittenInvariantsQuantization2001} gives
\begin{equation}
F_1^{X}\bigl(t_\nu(s_v)\bigr)
=F_1^{X}(u)+\overline F_1^{X}(a_\nu;u),
\label{eq:genus-one-vacuum-transformation}
\end{equation}
where $X$ denotes the respective theory, including the Hodge twist
in the diagonal case.

At a critical point $p_i$ of $x$, the ratio of the first odd
coefficients is $(dy_{\mathrm{sh}}^\nu/dy^\circ)(p_i)$.
Set $\mathsf{x}_i=\mathsf{x}(\ell_i)$ and $r=e^{u^2}$.
In the diagonal case, $dy^\circ=dp/p$, $\ell=e^{u^1}p$, and the
ratio equals $\mathsf{x}_i$. In the anti-diagonal case,
$dy^\circ=-\mathrm i\,dp/p$ and $\ell=\mathrm i e^{u^1}(p+r)$,
so the ratio equals $\mathsf{x}_ip_i/(\mathrm i(p_i+r))$.
The two anti-diagonal critical points satisfy
$p_+p_-=r$ and $(p_++r)(p_-+r)=r(1-r)$ by
\eqref{eq:phase-critical-points}. Hence
\begin{equation}
\begin{aligned}
\prod_{i\in\{+,-\}}\frac{y_{\mathrm{sh},1}^{\mathrm{di},i}}
     {y_1^{\circ,i}}
&=\mathsf{x}_+\mathsf{x}_-=e^{2u^1}(1-r)^2,\\
\prod_{i\in\{+,-\}}\frac{y_{\mathrm{sh},1}^{\mathrm{ad},i}}
     {y_1^{\circ,i}}
&=-\frac{\mathsf{x}_+\mathsf{x}_-}{1-r}=\frac{e^{2u^1}}{1-r}.
\end{aligned}
\label{eq:genus-one-slope-products}
\end{equation}
In the last equalities we used $\mathsf{x}_+\mathsf{x}_-=(v^1)^2-4e^{v^2}$
and the respective primary coordinate transformations.
Taking the logarithms on the branches fixed by the large-radius
normalization, \eqref{eq:genus-one-primary-potentials}--\eqref{eq:genus-one-slope-products}
give
\begin{equation}
\begin{aligned}
F_1^{X^{\mathrm{di}},\mathrm{tw}}\bigl(t_{\mathrm{di}}(s_v)\bigr)
&=-\frac{2u^1+u^2}{24}=-\frac{v^2}{24},\\
F_1^{X^{\mathrm{ad}}}\bigl(t_{\mathrm{ad}}(s_v)\bigr)
&=-\frac{2u^1+u^2+\log(1-r)}{24}=-\frac{v^2}{24}.
\end{aligned}
\label{eq:genus-one-target-vacua}
\end{equation}
Both values equal $F_1^{\mathbb P^1}(s_v)$, so both genus-one
constants vanish. This completes the proof.
\end{proof}

\section{Conclusion}
\label{sec:outlook}
In this paper, we derive transformations of A-model descendant Gromov-Witten potentials between $\bP^1$ and the resolved
conifold by comparing Eynard-Orantin topological recursion on the B-model side. The spectral curve for $\bP^1$ is 
from Norbury-Scott. The spectral curves for $X^{\rm di}, X^{\rm ad}$ are from the dimensional reduction of Landau-Ginzburg
mirror. For the comparison of spectral curves, we fix the data $(x,B)$ for spectral curves, and shift the $y$-coordinate.

The shifted spectral curves of $X^{\rm di}, X^{\rm ad}$ have same recursion kernel as the spectral curve for $\bP^1$. 
The shift of $y$-coordinate gives the effect as dilaton shift in CohFTs, 
which is directly written down in the ancestor background. The $S$-calibration gives the ancestor-descendant transformation,
so that we derive the transformation in A-model descendant invariants.


Several aspects of this method deserve further development.
An important direction is the resolved conifold with more general torus actions, where the comparison allow transformations that do not preserve the standard Bergman kernel.
For a local change of spectral parameter $\ell=f(p)$, the
difference between the two standard sphere kernels is
\[
\Delta B(p_1,p_2)
=
\left(
\frac{f'(p_1)f'(p_2)}{(f(p_1)-f(p_2))^2}
-\frac{1}{(p_1-p_2)^2}
\right)dp_1\,dp_2.
\]
This difference vanishes for the affine transformations used
in this paper, but is generally nonzero.
Extending the comparison therefore requires tracking its effect
on the recursion differentials and the associated ancestor and
descendant potentials. 

The available constructions of spectral curves for Frobenius manifolds and CohFTs
\cite{milanovEynardOrantinRecursion2015,fangRemodelingConjectureToric2019,dunin-barkowskiDubrovinsSuperpotentialGlobal2019,
xiongSpectralCurvesSemisimple2026}
provide further opportunities to study transitions between theories through transformations of their spectral data.
A related example is the transition between the generalized
Br\'ezin--Gross--Witten theory and the Witten--Kontsevich model
\cite{xiongEmergentGeometryGeneralized2026}.
There, matching the deformed spectral curves in renormalized
coordinates identifies the stable recursion differentials and
relates the positive-genus free energies with the appropriate normalization corrections.
It would be interesting to develop these examples into a systematic approach to correspondences between theories arising from different enumerative problems.

\vspace{.1in}
{\bf Declarations}

\vspace{.1in}
{\bf Data availability statement}
Data sharing is not applicable to this article as no datasets
were generated or analysed during the current study.

\vspace{.1in}
{\bf AI Assistance}
The first-named and second-named authors had a beneficial discussion with AI during the development of the mathematical constructions presented in this paper. AI tools were also used to polish the language of the manuscript. 
The authors take responsibility for all results of the paper.

\vspace{.1in}
{\bf Conflict of interest}
On behalf of all authors, the corresponding author states that there is no conflict of interest.

\providecommand{\bysame}{\leavevmode\hbox to3em{\hrulefill}\thinspace}
\providecommand{\MR}{\relax\ifhmode\unskip\space\fi MR }
\providecommand{\MRhref}[2]{%
  \href{http://www.ams.org/mathscinet-getitem?mr=#1}{#2}
}
\providecommand{\href}[2]{#2}

\end{document}